\documentclass[11pt]{amsart}
\usepackage{amssymb, paralist, xspace, graphicx, url, amscd, euscript, mathrsfs,stmaryrd,epic,eepic,color}
\usepackage[all]{xy}
\SelectTips{cm}{}
\usepackage[notcite,notref]{showkeys}

\numberwithin{equation}{section}

\numberwithin{subsection}{section}

\allowdisplaybreaks[1]

\newenvironment{enumeratei}
{\begin{enumerate}[\upshape (i)]}
{\end{enumerate}}

\newtheorem*{namedtheorem}{\theoremname}
\newcommand{\theoremname}{testing}

\newtheorem{theorem}[subsection]{Theorem}
\newtheorem{proposition}[subsection]{Proposition}
\newtheorem{proposition-definition}[subsection]
{Proposition-Definition}

\newtheorem{lemma}[subsection]{Lemma}

\theoremstyle{definition}
\newtheorem{definition}[subsection]{Definition}

\newtheorem{remark}[subsection]{Remark}
\newtheorem{question}[subsection]{Question}

\theoremstyle{remark}

\newcommand\cB{\mathcal{B}}
\newcommand\cC{\mathcal{C}}
\newcommand\cD{\mathcal{D}}
\newcommand\cE{\mathcal{E}}

\newcommand\cG{\mathcal{G}}

\newcommand\cI{\mathcal{I}}

\newcommand\cK{\mathcal{K}}

\newcommand\cM{\mathcal{M}}

\newcommand\cO{\mathcal{O}}

\newcommand\cX{\mathcal{X}}
\newcommand\cY{\mathcal{Y}}
\newcommand\cZ{\mathcal{Z}}

\newcommand\ocM{\overline{\mathcal{M}}}
\newcommand\ocX{\overline{\mathcal{X}}}

\newcommand\FF{\mathbb{F}}
\newcommand\GG{\mathbb{G}}

\newcommand\LL{\mathbb{L}}

\newcommand\PP{\mathbb{P}}
\newcommand\QQ{\mathbb{Q}}

\newcommand\ZZ{\mathbb{Z}}

\newcommand\bL{\mathbf{L}}

\newcommand\rmm{\mathrm{m}}

\newcommand\rmo{\mathrm{o}}

\newcommand\fB{\mathfrak{B}}
\newcommand\fC{\mathfrak{C}}

\newcommand\fM{\mathfrak{M}}

\newcommand\arr{\ifinner\to\else\longrightarrow\fi}
\newcommand\larr{\longrightarrow}

\renewcommand\H{\operatorname{H}}

\newcommand\eqdef{\overset{\mathrm{\scriptscriptstyle def}} =}

\def\displaytimes_#1{\mathrel{\mathop{\times}\limits_{#1}}}

\def\displayotimes_#1{\mathrel{\mathop{\bigotimes}\limits_{#1}}}

\renewcommand\hom{\operatorname{\underline {Hom}}}

\newcommand\Aut{\operatorname{Aut}}

\newcommand\pic{\operatorname{Pic}}

\newcommand\spec{\operatorname{Spec}}

\newdir{ >}{{}*!/-5pt/@{>}}

\newcommand\double{\rightrightarrows}

\newcommand\doublelong[2]{\mathbin{\xymatrix{{}\ar@<3pt>[r]^{#1}
\ar@<-3pt>[r]_{#2}&}}}

\newcommand{\underhom}
{\mathop{\underline{\mathrm{Hom}}}\nolimits}

\newlength{\ignora}

\newcommand{\der}{\operatorname{Der}}

\newcommand{\B}{\mathop{\cB}\nolimits}

\newcommand{\lr}{linearly reductive\xspace}

\newcommand{\mmu}{\boldsymbol{\mu}}

\newcommand{\bmu}{{\boldsymbol{\mu}}}

\newcommand{\fppf}{_{\mathrm{fppf}}}

\newcommand{\tw}{{\mathrm{tw}}}
\newcommand{\sm}{{\mathrm{sm}}}
\newcommand{\sh}{{\mathrm{sh}}}

\newcommand{\lci}{local complete intersection\xspace}

\newcommand{\thickslash}{\mathbin{\!\!\pmb{\fatslash}}}
\renewcommand{\setminus}{\smallsetminus}

\theoremstyle{plain}
\newtheorem{thm}[subsection]{Theorem}
\newtheorem*{thm*}{Theorem}
\newtheorem{cor}[subsection]{Corollary}
\newtheorem{lem}[subsection]{Lemma}
\newtheorem{prop}[subsection]{Proposition}
\newcommand{\mc}{\mathcal}
\newcommand{\Sp}{\text{\rm Spec}}
\theoremstyle{definition}
\newtheorem{pg}[subsection]{}
\newtheorem{rem}[subsection]{Remark}
\newtheorem{defn}[subsection]{Definition}
\newcommand{\Z}{\mathbb{Z}}
\newcommand{\scr}[1]{\mathbf{\EuScript{#1}}}
\newcommand{\lsm}[1]{\mathcal{M}_{#1}}
\newcommand{\mls}{\mc }
\newcommand{\Hom}{\underline {\text{\rm Hom}}}
\newcommand{\Sec}{\underline {\text{\rm Sec}}}
\newcommand{\rigSec}{{\text{\rm Sec}}}

\numberwithin{equation}{subsection}

\begin{document}

\title[Twisted yet tame]{Twisted stable maps to tame Artin stacks}

\author[Abramovich]{Dan Abramovich}

\author[Olsson]{Martin Olsson}

\author[Vistoli]{Angelo Vistoli}

\address[Abramovich]{Department of Mathematics\\
Brown University\\
Box 1917\\
Providence, RI 02912\\
U.S.A.}
\email{abrmovic@math.brown.edu}

\address[Olsson]{Department of Mathematics \#3840\\
University of California\\
Berkeley, CA 94720-3840\\
U.S.A.}
\email{molsson@math.berkeley.edu}

\address[Vistoli]{Scuola Normale Superiore\\ Piazza dei Cavalieri 7\\
56126 Pisa\\ Italy}
\email{angelo.vistoli@sns.it}

\thanks{Vistoli supported in part by the PRIN Project ``Geometria
sulle variet\`a algebriche'', financed by MIUR. Olsson partially supported by NSF grant DMS-0555827 and an Alfred P. Sloan fellowship. Abramovich support in part by NSF grants DMS-0301695 and DMS-0603284.}
\date{\today}

\maketitle

\tableofcontents

\section{Introduction}

\subsection{Context and statement of main result} This paper is a continuation of \cite{AOV}, where the basic theory of tame Artin stacks is developed. Our main goal here is the construction of an appropriate analogue of Kontsevich's  space of stable maps in the case where the target is a tame Artin stack. When the target is a tame Deligne--Mumford stack, the theory was developed in \cite{A-V}, and found a number of applications. The theory  for arbitrary tame Artin stacks developed here is very similar, but it is necessary to overcome a number of technical hurdles and to generalize a few questions of foundation. However the method of construction we use here is very different from the ad-hoc method of \cite{A-V}, and more natural: we rely on the second author's general results in \cite{logtwisted,homstack,bounded}, as extended in the appendices.

The main result is the following.  Let $\mc M$ be a finitely presented tame algebraic stack over a scheme $S$ with finite inertia, and fix integers $g$ and $n$.  As we explain in section \ref{section4}, there is a natural notion of twisted curve in arbitrary characteristic, defined using the theory of tame stacks.  Using this notion we can then consider certain twisted stable maps
$$
\mc C\rightarrow \mc M,
$$
from $n$-marked genus $g$  twisted curves $\mc C$ to $\mc M$ (see section \ref{Section:3} for precise definitions), and obtain a stack $\cK_{g, n}(\cM)$ classifying twisted stable maps to $\mc M$.  The main result of the paper is the following:
\begin{thm*}[{Theorem \ref{Th:3.1} in the text}]  Let $M$ denote the coarse moduli space of $\mc M$.
The stack $\cK_{g,n}(\cM)$ of twisted stable maps
from $n$-pointed genus $g$ curves into $\cM$  is a locally finitely presented
algebraic stack over $S$, which is proper and quasi-finite over the stack of stable
maps into $M$. If $M$ is projective over a field, the substack $\cK_{g,n}(\cM,\beta)$   of twisted stable maps
from $n$-pointed genus $g$ curves into $\cM$ with target class $\beta$ is proper, and is open and closed in $\cK_{g,n}(\cM)$.
\end{thm*}

\subsection{The paper is organized as follows}

Section \ref{section4} is devoted to a number of basic properties of twisted curves in arbitrary characteristic, where, unlike the situation of \cite{A-V}, they may fail to be Deligne--Mumford stacks. In particular we show in Proposition \ref{Prop:2.3} that the notion of twisted curve is stable under small deformations and can be tested on geometric fibers.

In Section \ref{interlude} we collect together some facts about relative moduli spaces which will be used in the following section.

In Section  \ref{Section:3} we define twisted stable maps into a tame stack $\cX$ and show in Theorem \ref{Th:3.1} that they form an Artin stack, which is proper and quasi-finite over the Kontsevich space of the coarse moduli space $X$ of $\cX$. In particular, when $X$ is projective, the stacks of twisted stable maps of $\cX$ with fixed numerical data admit projective coarse moduli spaces. The construction proceeds rather naturally from the following:
\begin{enumerate}
\item  the existence of the stack of twisted curves, which was shown in \cite{logtwisted} for Deligne--Mumford twisted curves and proved in our case in Appendix \ref{AppendixA},  and
\item The existence and finiteness properties of $\Hom$-stacks, proved in \cite{homstack} and \cite{bounded} in many cases and extended in our case in Appendix \ref{AppendixB}.
\end{enumerate}
Some extra care is needed for proving the quasi-finite claim. Properness relies on Lemma \ref{Lem:purity}, a suitable generalization  of the Purity Lemma of \cite{A-V}.

In section \ref{reductionsection} we concentrate on the case where   the target stack is $\cB G$, the classifying stack of a finite flat linearly reductive group scheme $G$. The main result here, Theorem \ref{flattheorem}, is that the stack $\cM_{g,n}(\cB G)$ of stable maps with target
$\cB G$ is finite and flat over the Deligne--Mumford stack $\cM_{g,n}$. This result is  known when $G$ is tame and \'etale (see \cite{ACV}) and relatively straightforward when $G$ is  a diagonalizable group scheme (or even a twist of such). However, our argument  in general takes some delicate twists and turns.

As an example for the behavior of these stacks, we consider in Section \ref{examplesection} two ways to compactify the moduli space $X(2)$ of elliptic curves with full level-2 structure. The first is as a component in $\cK_{0,4}(\cB\bmu_2)$, the stack of totally branched $\bmu_2$-covers of stable 4-pointed curves of genus 0.  This provides an opportunity to consider the cyclotomic inertia stacks and evaluation maps. The second is as a component in $\cK_{1,1}(\cB\bmu_2^2)$ parametrizing elliptic curves with a map to $\cB\bmu_2^2$, where  $X(2)$ meets other components in a geometrically appealing way. We also  find the Katz--Mazur regular model of $X(2)$ as the closure of a component of the generic fiber.

As already mentioned, the paper contains three appendices, written by the second author.  Appendix \ref{AppendixA} generalizes the main results of \cite{logtwisted}. Appendix \ref{AppendixC} contains some preparatory results needed for Appendix \ref{AppendixB} which  generalizes some of the results from \cite{homstack} and \cite{bounded} to tame stacks.  The logical order of Appendix  \ref{AppendixA} is after section \ref{section4}, whereas the Appendix \ref{AppendixB} only uses results from \cite{AOV}.

\subsection{The question of the correct generality} One may ask whether the main result, and the tools leading to it, are set in the correct generality. 

First, it appears to us that the main result, especially the crucial statement that $\cK_{g,n}(\cM) \to \cK_{g,n}(M)$ is proper, is the most general possible, at least without substantially new ideas.
Already for $\cM = \B\Z/(p)$,  the classifying stack of the group $\Z/(p)$ in characteristic $p>0$, the result fails.  For example, consider a discrete valuation ring $R$ with residue field $k$ of characteristic $p$, and let $E/R$ be an elliptic curve whose closed fiber $E_k$ is supersingular.  Let $K$ be the fraction field of $R$, and let $P_K\rightarrow E_K$ be a nontrivial $\Z/(p)$-torsor (for example we might have $K$ of characteristic $0$).  Let $P\rightarrow E$ be the normalization of $E$ in $P_K$.  Then $P\rightarrow E$ must be ramified along the entire closed fiber $E_k$.  For if not, the map $P\rightarrow E$ must be \'etale by purity, in which case we get a $\Z/(p)$-torsor $P_k\rightarrow E_k$.  Since $E_k$ is supersingular, this $\Z/(p)$-torsor is trivial, which in turn implies that $P$ and $P_K$ are also trivial.  This implies that there is no morphism of stacks $\mls E\rightarrow E$ which is an isomorphism away from the origin of $E$, for which $P_K$ extends to a $\Z/(p)$-torsor over $\mls E$.

This example in particular shows that sections \ref{reductionsection} and \ref{examplesection} do not generalize to apply for such groups. In this context there are some other ideas which can be used, notably Raynaud's constructions in \cite{Raynaud}, see also \cite{raynaudgroup}, and Romagny's effective model \cite{Romagny}, see also \cite{Tossici-actions}. The point is the following: suppose $R$ is a  
  discrete valuation ring $R$ as above, $Y \to X$ a finite cover of nodal curves, $G$ is a finite group acting on $Y$ over $X$, with $p\big\vert \,|G|$, and $X_K = Y_K /G$. Then one can try to replace $G$ by a more degenerate group scheme  $\cG/\cX$, with $\cX$ a twisted curve over $X$, in such a way that $\cG$ has a better behaved  action on $Y$. The paper \cite{raynaudgroup} suggests that under certain conditions on $G$  this would lead to a proper moduli space; this is a project the first author and Romagny plan to pursue, where results of the present paper are crucial.
  
One of the many points in the proof where the assumption that $\cM$ be tame is used is Proposition \ref{Prop:B3} on a quotient presentation of a relative moduli space. This does not generalize if $\cM$ is not tame -- see Remark \ref{relative-example}.

Back to the main theorem, one can consider the case where $\cM/S$ is finitely presented, and ask what can be said. In case $\cM$ has finite diagonal the stack $\cK_{g,n}(\cM)$ is still a locally finitely presented algebraic stack -- see Proposition \ref{Prop:existence}. This relies on Theorem \ref{athm} (i). The following questions around Theorem \ref{athm} are probably of interest both on their own right and in relation to the main theorem. The notation as in Theorem \ref{athm}, with $\cX$ flat and proper over $B$ with finite diagonal, but $\cY$ only assumed of finite presentation over $B$:

\begin{question} Under what conditions on $\cX, \cY$ is $\underhom(\cX,\cY)$ an algebraic stack locally of finite presentation over $B$? Note that the general statement of \cite{Aoki} is retracted in \cite{Aokierratum} because of a sheaf-theoretic gap in the proof, so the correct generality is yet to be pinned down.
\end{question}

\begin{question} Suppose now $\cY \to Y$ is proper. Under what further conditions  is $\underhom(\cX,\cY) \to \underhom(\cX,Y)$ of finite type?  In particular, is it sufficient that $\cY\to Y$ have finite diagonal?
\end{question}

\begin{question} Suppose further $\cY \to Y$ is proper and quasi-finite. Under what further conditions  is $\underhom(\cX,\cY) \to \underhom(\cX,Y)$  quasi-finite?  Note that since $\Aut(\alpha_p)=\GG_m$  is not quasi-finite, it is not enough for $\cX, \cY$ to have finite diagonal.  Theorem \ref{athm} (ii) shows that it suffices for $\cY$ to be tame. Is it also sufficient for $\cX$ to be  tame?
\end{question}

\subsection{Acknowledgements} Thanks to Johan de Jong for helpful comments. Thanks to Shaul Abramovich for help and advice with visualization. We are indebted to  the referee who provided a preliminary opinion on our paper, whose comments led us to try to optimize some of the results, and to the second referee for doing a very careful job. 

\section{Twisted curves}\label{section4}

\begin{defn}\label{A1} Let $S$ be a scheme. An \emph{$n$-marked twisted  curve over $S$} is a collection of data $(f\colon\mc C\rightarrow S, \{\Sigma _i\subset \mc C\}_{i=1}^n)$ as follows:
\begin{enumeratei}
\item $\mc C$ is a proper tame  stack over $S$ whose geometric fibers are connected of dimension $1$, and such that the moduli space $C$ of $\mc C$ is a nodal curve over $S$.
\item The $\Sigma _i\subset \mc C$ are closed substacks which are fppf gerbes over $S$, and whose images in $C$ are contained in the smooth locus of the morphism $C\rightarrow S$.
\item If $\mc U\subset \mc C$ denotes the complement of the $\Sigma _i$ and the singular locus of $\mc C\rightarrow S$, then $\mc U\rightarrow C$ is an open immersion.
\item For any geometric point $\bar p\rightarrow C$ mapping to a smooth point of $C$, there exists an integer $r$ such that
\begin{equation*}
\Sp (\mc O_{C, \bar p})\times _C\mc C\simeq [D^\sh
/\mu _r],
\end{equation*}
where 
$D^\sh$
%$\widetilde{\mc O_{S, f(\bar p)}[z]}$
 denotes the strict henselization of $D\colon= \Sp (\mc O_{S, f(\bar p)}[z])$ at the point $(\mathfrak{m}_{S, f(\bar p)}, z)$ and $\zeta \in \mu _r$ acts by $z\mapsto \zeta \cdot z$ (note that $r=1$ unless $\bar p$ maps to a point in the image of some $\Sigma _i$). Here $\mathfrak{m}_{S, f(\bar p)}$ denote the maximal ideal in the strict henselization $\mc O_{S, f(\bar p)}$.
\item If $\bar p\rightarrow C$ is a geometric point mapping to a node of $C$, then there exists an integer $r$ and an element $t\in \mathfrak{m}_{S, f(\bar p)}$ such that
\begin{equation*}
 \Sp (\mc O_{C, \bar p})\times _C\mc C\simeq [D^\sh
 /\mu _r],
 \end{equation*}
where $
D^\sh
$ denotes the strict henselization of $$D\colon= \Sp (\mc O_{S, f(\bar p)}[z, w]/(zw-t))$$ at the point $(\mathfrak{m}_{S, f(\bar p)}, z, w)$ and $\zeta \in \mu _r$ acts by $x\mapsto \zeta \cdot z$ and $y\mapsto \zeta ^{-1}\cdot y$.
\end{enumeratei}
\end{defn}

\begin{rem} If $\mc C\rightarrow S$ is a proper tame Artin stack which admits a collection of closed substacks $\Sigma _i\subset \mc C$ ($i=1,\dots , n$ for some $n$) such that $(\mc C, \{\Sigma _i\}_{i=1}^n)$ is an $n$-marked twisted curve, we will refer to $\mc C$ as a \emph{twisted curve}, without reference to  markings.
\end{rem}

The following proposition allows us to detect twisted curves on fibers:

\begin{prop}\label{Prop:2.3} Let $S$ be a scheme, and let $(f\colon\mc C\rightarrow S, \{\Sigma _i\}_{i=1}^n)$ be a proper flat  tame  stack $f\colon\mc C\rightarrow S$ with a collection of closed substacks $\Sigma _i\subset \mc C$ which are $S$-gerbes.  If for some geometric point $\bar x\rightarrow S$ the fiber $(\mc C_{\bar x}, \{\Sigma _{i, \bar x}\})$ is an $n$-marked twisted  curve, then there exists an open neighborhood of $\bar x$ in $S$ over which $(\mc C\rightarrow S, \{\Sigma _i\})$ is an $n$-marked twisted  curve.
\end{prop}
\begin{proof} Let $\pi \colon\mc C\rightarrow C$ be the coarse moduli space of $S$.  Since $\mc C$ is a tame  stack flat over $S$, 
 the space $C$ is flat over $S$,
 and for any morphism $S'\rightarrow S$ the base change $C\times _SS'$ is the coarse moduli space of $\mc C\times _SS'$, by \cite[Corollary 3.3]{AOV}. It follows that in some \'etale neighborhood of $\bar x\rightarrow S$ the space $C/S$ is a nodal curve over $S$ (see for example \cite[\S 1]{D-M}).  Shrinking on $S$ we may therefore assume that $C$ is a nodal curve over $S$. After further shrinking on $S$, we can also assume that the images of the $\Sigma _i$ in $C$ are contained in the smooth locus of the morphism $C\rightarrow S$ and that (\ref{A1} (iii)) holds.
To prove the proposition it then suffices to verify that conditions (iv) and (v) hold for geometric points of $\bar p\rightarrow C$ over $\bar x$.  

We may without loss of generality assume that $S$ is the spectrum of a strictly henselian local ring whose closed point is $\bar x$. Also let $\mc C_{\bar p}$ denote the fiber product
\begin{equation*}
\mc C_{\bar p}:= \mc C\times _{C, \bar p}\Sp (\mc O_{C, \bar p}).
\end{equation*}

Consider first the case when $\bar p\rightarrow C$ has image in the smooth locus of $C/S$.  Since the closed fiber $C_{\bar x}$ is a twisted curve, we can choose an isomorphism  
\begin{equation*}
\mc C_{\bar p}\times _S\Sp (k(\bar x))\simeq [\Sp (
{k(\bar x)[z]})^\sh
%\widetilde {k(\bar x)[z]})
/\mu _r],
\end{equation*}
for some integer $r\geq 1$, where the $\mu _r$-action is as in (\ref{A1} (iv)).  Let 
\begin{equation*}
P_0\rightarrow \mc C_{\bar p}\times _S\Sp (k(\bar x))
\end{equation*}
be the $\mu _r$-torsor
\begin{equation*}
\Sp 
({k(\bar x)[z]})^\sh
%(\widetilde {k(\bar x)[z]})
\rightarrow [\Sp 
({k(\bar x)[z]})^\sh
%(\widetilde {k(\bar x)[z]})
/\mu _r].
\end{equation*}
As in the proof of \cite[Proposition 3.6]{AOV} there exists a $\mu _r$-torsor $P\rightarrow \mc C_{\bar p}$ whose reduction is $P_0$.  Furthermore, $P$ is an affine scheme over $S$. Since $\mc C_{\bar p}$ is flat over $S$, the scheme $P$ is also flat over $S$.  Since its closed fiber is smooth over $S$, it follows that $P$ is a smooth $S$--scheme of relative dimension $1$ with $\mu _r$-action.  Write $P = \Sp (A)$ and fix a $\mu _r$-equivariant isomorphism
\begin{equation*}
A\otimes _{\mc O_S}k(\bar x)\simeq 
({k(\bar x)[z]})^\sh.
%\widetilde {k(\bar x)[z]}.
\end{equation*}
Since $\mu _r$ is linearly reductive we can find an element $\widetilde z\in A$ lifting $z$ such that an element $\zeta \in \mu _r$ acts on $\widetilde z$ by $\widetilde z\mapsto \zeta \cdot \widetilde r$.  Since $A$ is flat over $\mc O_S$ the induced map
\begin{equation*}
(\mc O_S[\widetilde z])^{\text{sh}}\rightarrow A
 \end{equation*}
 is an isomorphism, which shows that (iv) holds.

 The  verification of (v) requires some deformation theory. First we give an explicit argument. We then provide an argument using the cotangent complex.

 Let $\bar p\rightarrow C$ be a geometric point mapping to a node of $C$.  Then we can write
 \begin{equation*}
 \mc C_{\bar p} = [\Sp (A)/\mu _r],
 \end{equation*}
for some integer $r$, and $A$ is a flat $\mc O_{S, \bar x}$-algebra such that there exists an isomorphism
\begin{equation*}
\iota \colon A\otimes _{\mc O_{S, \bar x}}k(\bar x)\simeq 
({k(\bar x)[z,w]/(zw)})^\sh,
%\widetilde {k(\bar x)[z,w]/(zw)},
\end{equation*}
with $\mu _r$-action as in (v).  We claim that we can lift this isomorphism to a $\mu _r$-equivariant isomorphism
\begin{equation}\label{localiso}
A\simeq 
({\mc O_{S, \bar x}[\widetilde z, \widetilde w]/(\widetilde z\widetilde w-t)})^\sh
%\widetilde {\mc O_{S, \bar x}[\widetilde z, \widetilde w]/(\widetilde z\widetilde w-t)}
\end{equation}
for some $t\in \mc O_{S, \bar x}$.  For this note that $A^{\mu _r}$ is equal to $\mc O_{C, \bar p}$ and that $A$ is a finite $\mc O_{C, \bar p}$-algebra.  By a standard application of the Artin approximation theorem, in order to find \ref{localiso} it suffices to prove the analogous statements for the map on completions (with respect to $\mathfrak{m}_{C, \bar p}$)
\begin{equation*}
\widehat {\mc O}_{C, \bar p}\rightarrow \widehat A.
\end{equation*}

For this in turn we inductively find an element $t_q\in \mc O_{S, \bar x}/\mathfrak{m}_{S, \bar x}^q$ and a $\mu _r$-equivariant isomorphism
\begin{equation*}
\rho _q\colon( ({\mc O_{S, \bar x}/\mathfrak{m}_{S, \bar x}^q)[z,w]/(zw-t_q)})^\wedge \rightarrow \widehat A/\mathfrak{m}_{S, \bar x}^q\widehat A.
\end{equation*}
For $\rho _1$ we take the isomorphism induced by $\iota $ and $t_1 = 0$. For the inductive step we assume $(\rho _q, t_q)$ has been constructed and find $(\rho _{q+1}, t_{q+1})$.  For this choose first any liftings $\widetilde z, \widetilde w\in \widehat A/\mathfrak{m}_{S, \bar x}^{q+1}\widehat A$ of $z$ and $w$ such that $\zeta \in \mu _r$ acts by
\begin{equation}\label{goodaction}
\widetilde z\mapsto \zeta \cdot \widetilde z, \ \ \widetilde w\mapsto \zeta ^{-1}\cdot \widetilde w.
\end{equation}
This is possible because $\mu _r$ is linearly reductive.

Since $\widehat A$ is flat over $\mc O_{S, \bar x}$ we have an exact sequence
$$
0\rightarrow \mathfrak{m}_{S, \bar x}^q/ \mathfrak{m}_{S, \bar x}^{q+1}\otimes _{k(\bar x)} \widehat A/\mathfrak{m}_{S, \bar x}\rightarrow \widehat A/\mathfrak{m}_{S, \bar x}^{q+1}\widehat A\rightarrow \widehat A/\mathfrak{m}_{S, \bar x}^{q}\widehat A\rightarrow 0.
$$
Choosing any lifting $t_{q+1}\in \mc O_{S, \bar x}/\mathfrak{m}_{S, \bar x}^{q+1}$ of $t_q$ and consider
$$
\widetilde z\widetilde w-t_{q+1}.
$$
This is a $\mu _r$-invariant element of 
$$
\mathfrak{m}_{S, \bar x}^q/ \mathfrak{m}_{S, \bar x}^{q+1}\otimes _{k(\bar x)} \widehat A/\mathfrak{m}_{S, \bar x}\simeq \mathfrak{m}_{S, \bar x}^q/ \mathfrak{m}_{S, \bar x}^{q+1}\otimes _{k(\bar x)}( ({\mc O_{S, \bar x}/\mathfrak{m}_{S, \bar x})[z,w]/(zw)})^\wedge.
$$
It follows that after possibly changing our choice of lifting $t_{q+1}$ of $t_q$, we can write
\begin{equation*}
\widetilde z\widetilde w = t_{q+1}+\widetilde z^rg+\widetilde w^rh,
\end{equation*}
where $t_{q+1}\in \mc O_{S, \bar x}/\mathfrak{m}_{S, \bar x}^{q+1}$ and reduces to $t_{q}$, and $g, h\in \mathfrak{m}_{S, \bar x}^{q}\widehat A/\mathfrak{m}_{S, \bar x}^{q+1}\widehat A$ are $\mu _r$-invariant.   Replacing $\widetilde z$ by $\widetilde z-\widetilde w^{r-1}h$ and $\widetilde w$ by $\widetilde w-\widetilde z^{r-1}g$ (note that with these new choices the action of $\mu _r$ is still as in \ref{goodaction}) we obtain $t_{q+1}$ and $\rho _{q+1}$.
\end{proof}

\subsection{The cotangent complex and an alternative proof of Proposition~\ref{Prop:2.3}(v)} We give an alternative proof of (v)  using  a description of the cotangent complex of $\cC$, which might be of  interest on its own.  Since $\cC$
is an Artin stack, a few words are in order. 

The cotangent complex
$\LL_{\cX/\cY}$ of a morphism of 
Artin stacks $\cX\to \cY$ is defined in \cite[Chapter 17]{L-MB} as
an object in the derived category of quasi-coherent sheaves in
the lisse-\'etale site of  $\cX$. Unfortunately, as was observed by
Behrend and Gabber, this site is not functorial, which had the
potential of rendering  $\LL_{\cX}$ both incomputable and
useless. However, in \cite{Olsson-sheaves}, see especially section 8, it
is shown that the
lisse-\'etale site has surrogate properties replacing functoriality
which are in particular sufficient for dealing with  $\LL_{\cX}$, and
for deformation theory.
 In \cite{Olsson-sheaves} a cotangent complex $\LL_{\cX}$ with the appropriate properties is only defined as
an object of a carefully constructed filtered category $D'_{qcoh}(\cX_{\text {lis-\'et}})$. This issue is removed in \cite[2.2.5]{Laszlo-Olsson}, specifically the equivalence at the end of page 119 between the appropriate derived categories of quasi-coherent sheaves on the stack and on a simplicial resolution. 
In particular,    to any morphism of Artin stacks $f:X\to Y$ we can after all associate a cotangent complex $\LL_f\in D^{\le 1}_{coh}(X)$ which has good properties.

 Two key properties are the following:
\begin{enumerate}
\item given a
morphism $f\colon\cX \to \cY$ of algebraic stacks over a third $\cZ$, there
is a natural
distinguished triangle
$$\bL f^*\LL_{\cY/\cZ} \to \LL_{\cX/\cZ} \to \LL_{\cX/\cY}
\stackrel{[1]}{\to}.$$
\item Given a fiber diagram
$$\xymatrix{\cX'\ar^f[r]\ar[d] & \cX \ar[d]\\ \cY' \ar[r] & \cY}$$
with the horizontal maps flat, we have $f^*\LL_{\cX/\cY} = \LL_{\cX'/\cY'}$.
\end{enumerate}

If the base scheme is $S$,  consider $S \to \cB\GG_{m,S}$. Applying
the above to the fiber square $$\xymatrix{\GG_{m,S}\ar[r]\ar[d] & S \ar[d]\\ S \ar[r] & \cB\GG_{m,S}}$$ and the morphisms $S \to \cB\GG_{m,S} \to S$  it is easy to see that
$\LL_{\cB\GG_{m,S}/S}=\cO_{\cB\GG_{m,S}}[-1]$. (In fact for a general smooth group scheme $G$ we have $\LL_{S/\cB G} = Lie(G)^*$ and therefore $\LL_{\cB G/S} = Lie(G)^*[-1]$ considered as a $G$-equivariant object.)

Let now $S = \Sp \ k$ with $k$ separably closed,
and consider a twisted curve $\cC/S$ with coarse moduli space $\pi \colon\cC\rightarrow C$.  Let $\bar p\rightarrow C$ be a geometric point mapping to a node, and fix an isomorphism
$$
\cC_{\bar p}:= \cC\times _C\bar p\simeq [D^{\text{sh}}/\mu _r],
$$
where $D:= \Sp (k[z, w]/(zw))$.     By a standard limit argument, we can thicken this isomorphism to a diagram
$$
\xymatrix{
& \cC _U\ar[ld]_-\alpha \ar[rd]^\beta \ar[dd]^-{\pi _U}& \\
\cC \ar[dd]^-\pi && [D/\mu _r]\ar[dd]\\
& U\ar[ld]_-a\ar[rd]^-b& \\
C&& D^{\mu _r},}
$$
where  $D^{\mu _r}$ denote the coarse space of $[D/\mu _r]$ (so $D^{\mu _r}$ is equal to the spectrum of the invariants in $k[z,w]/(zw)$), $a$ and $b$ are \'etale, $U$ is affine, and the squares are cartesian.  In this setting we can calculate a versal deformation of $\mc C_U$ as follows.

First of all we have the deformation
$$
[\Sp (k[z, w, t]/(zw-t))/\mu _r]
$$
of $[D/\mu _r]$.  Since $\beta $ is \'etale and representable and by the invariance of the \'etale site under infinitesimal thickenings, this also defines a formal deformation  (i.e. compatible family of deformations over the reductions)
$$
\widehat {\mc C}_{U, t}\rightarrow \text{Spf}(k[[t]])
$$
of ${\mc C}_{U}$.  We claim that this deformation is versal.  Since this deformation is nontrivial modulo $t^2$, to verify that $\widehat {\mc C}_{U, t}$ is versal it suffices to show that the deformation functor of $\mc C_U$ is unobstructed and that the tangent space is $1$-dimensional.  For this it suffices to show that
$$
\text{Ext}^2(\LL_{\cC _U}, \mc O_{\cC}) = 0
$$
and
$$
\text{Ext}^1(\LL_{\cC _U}, \mc O_{\cC}) = k.
$$
The map $\cC_U\rightarrow [D/\mu _r]$ induces a morphism $\cC_U\rightarrow \cB\mu _r$.  Consider the composite map $\cC_U\rightarrow \cB\mu _r\rightarrow \cB\mathbb{G}_m$, and the distinguished triangle associated to the resulting diagram
$$
\cC_U\rightarrow \cB \mathbb{G}_m\rightarrow \Sp (k).
$$
  Since $\mathbb{L}_{\cB\mathbb{G}_m/k}\simeq \mc O_{\cB\mathbb{G}_m}[-1]$ this  distinguished triangle can be written as
$$
\xymatrix{
\mc O_\cC[-1]\ar[r]&\mathbb{L}_{\mc C_U}\ar[r]& \mathbb{L}_{\cC_U/\cB\mathbb{G}_m}\ar[r]^- {[1]}&.}
$$
Considering the associated long exact sequence of $\text{Ext}$-groups one sees that 
$$
\text{Ext}^i(\LL_{\cC _U}, \mc O_{\cC}) = \text{Ext}^i(\LL_{\cC _U/\cB\mathbb{G}_m}, \mc O_{\cC})
$$
for $i>0$.  Now to prove that 
$$
\text{Ext}^2(\LL_{\cC _U/\cB\mathbb{G}_m}, \mc O_{\cC}) = 0,  \ \ \text{and} \ \ \text{Ext}^1(\LL_{\cC _U/\cB\mathbb{G}_m}, \mc O_{\cC}) = k,
$$
it suffices to show that $\LL_{\cC/\cB\mathbb{G}_m}$ can be represented by a two-term complex
$$
F^{-1}\rightarrow F^0
$$
with the $F^i$ locally free sheaves of finite rank and the cokernel of the map
$$
\mls Hom(F^0, \mc O_{\cC_U})\rightarrow \mls Hom(F^{-1}, \mc O_{\cC_U})
$$
isomorphic to the structure sheaf $\mc O_{\mc Z}$ of the singular locus $\mc Z\subset \cC_U$ (with the reduced structure).  

Now since $\beta $ is \'etale, the  complex $\LL_{\cC_U/\cB\mathbb{G}_m}$ is isomorphic to $\beta ^*\LL_{[D/\mu _r]/\cB\mathbb{G}_m}$.  Therefore to construct such a presentation $F^\cdot $ of $\LL_{\cC_U/\cB\mathbb{G}_m}$, it suffices to construct a corresponding presentation of $\LL_{[D/\mu _r]/\cB\mathbb{G}_m}$.

For this consider the  fiber diagram
$$\xymatrix{X\ar[r]\ar[d] & \cC \ar[d]\\ S \ar[r] & \cB\GG_{m}}$$
Here $X$ is the surface $\Sp (k[z', w', u, u^{-1}]/(z'w'))$ endowed with the action of $\GG_m$ via $t\cdot(z',w',u) = (tz', t^{-1}w', t^r u)$. Since $X/S$ is a local complete intersection, the map $\cC\to \cB\GG_m$ is also a representable local complete intersection morphism, and therefore the natural map
$\LL_{\cC/\cB\GG_m}\rightarrow H_0(\LL_{\cC/\cB\GG_m}) = \Omega^1_{\cC/\cB\GG_m}$ is a quasi-isomorphism. Concretely,
the pullback of $\LL_{\cC/                
  \cB\GG_{m}}$ is isomorphic to $\Omega^1_{X/S}$. Now we have a standard two-term resolution 
  $$0\to \cO_X \to \cO_X^3 \to  \Omega^1_{X/S}\to 0$$ where the first map is given by $(w',z',0)$ and the second by $(dz', dw', du/u)$.  Our desired presentation of $\LL_{[D/\mu _r]/\cB\GG_m}$ is then obtained by noting that this resolution  clearly
  descends to a locally free resolution on $\cC$.

\begin{rem} Let $S$ be a scheme.  A priori the collection of twisted $n$-pointed  curves over $S$  (with morphisms only isomorphisms) forms a $2$-category.  However, by the same argument proving \cite[4.2.2]{A-V} the 2-category of twisted $n$-pointed  curves over $S$ is equivalent to a $1$-category.  In what follows we will therefore consider the category of twisted $n$-pointed  curves, whose arrows are isomorphism classes of 1-arrows in the 2-category.
\end{rem}

\begin{prop}\label{cohlemb} Let $f\colon\mc C\rightarrow S$ be a twisted curve.  Then the adjunction map $\mc O_S\rightarrow f_*\mc O_{\mc C}$ is an isomorphism, and for any quasi-coherent sheaf $\mc F$ on $\mc C$ we have $R^if_*\mc F = 0$ for $i\geq 2$.
\end{prop}
\begin{proof}
Let $\pi \colon\mc C\rightarrow C$ be the coarse space, and let $\bar f\colon C\rightarrow S$ be the structure morphism.  Since $\pi _*$ is an exact functor on the category of quasi-coherent sheaves on $\mc C$, we have for any quasi-coherent sheaf $\mc F$ on $\mc C$ an equality
$$
R^if_*\mc F = R^i\bar f_*(\pi _*\mc F).
$$
Since $\bar f\colon C\rightarrow S$ is a nodal curve this implies that $R^if_*\mc F = 0$ for $i\geq 2$, and the first statement follows from the fact that $\pi _*\mc O_{\mc C} = \mc O_C$.
\end{proof}

We conclude this section by recording some facts about the Picard functor of a twisted curve (these results will be used in section \ref{reductionsection} below).

Let $S$ be a scheme and $f\colon\mc C\rightarrow S$ a twisted curve over $S$.   Let $\scr Pic_{\mc C/S}$ denote the stack over $S$ which to any $S$-scheme $T$ associates the groupoid of line bundles on the base change $\mc C_T:= \mc C\times _ST$.  

\begin{prop} The stack $\scr Pic_{\mc C/S}$ is a smooth locally finitely presented algebraic stack over $S$, and for any object $\mc L\in \scr Pic_{\mc C/S}(T)$ (a line bundle on $\mc C_T$), the group scheme of automorphisms of $\mc L$ is canonically isomorphic to $\mathbb{G}_{m, T}$.
\end{prop}
\begin{proof} That $\scr Pic_{\mc C/S}$ is algebraic and locally finitely presented is a standard application of Artin's criterion for verifying that a stack is algebraic \cite{Artin}. See \cite[Th\'eor\`emes 1.1, 1.2, and Section 3, esp. 3.1.3]{Brochard} for details. This stack  can also be constructed  as  a $\hom$-stack,  see \cite{Aoki}.

To see that $\scr Pic_{\mc C/S}$ is smooth we apply the infinitesimal criterion. For this it suffices to consider the case where $S$ affine and we have a closed subscheme $S_0\subset S$ with square-0 ideal $I$. Writing $\mc C_0$ for  $\mc C \times_S S_0$, we have an exact sequence
$$0 \to \mc O_{\mc C}\otimes I \to \mc O_{\mc C} \to \mc O_{\mc C_0} \to 0$$ 
giving an exact sequence of cohomology
$$\scr Pic_{\mc C/S}(S) \to \scr Pic_{\mc C/S}(S_0) \to R^2f_*\mc O_{\mc C}\otimes I.$$
By Proposition \ref{cohlemb} the term on the right vanishes, giving the existence of the desired lifting of $\scr Pic_{\mc C/S}(S) \to \scr Pic_{\mc C/S}(S_0)$.

The statement about the group of automorphisms follows from the fact that the map $\mc O_S\rightarrow f_*\mc O_{\mc C}$ is an isomorphism, and the same remains true after arbitrary base change $T\rightarrow S$.
\end{proof}

Let $\text{Pic}_{\mc C/S}$ denote the rigidification of $\scr Pic_{\mc C/S}$ with respect to the group scheme $\mathbb{G}_m$ \cite[Proposition IV.2.3.18]{giraud}, see also \cite[Appendix A]{AOV}.

Let $\pi \colon\mc C\rightarrow C$ denote the coarse moduli space of $\mc C$ (so $C$ is a nodal curve over $S$).   Pulling back along $\pi$ defines a fully faithful functor
\begin{equation*}
\pi ^*\colon\scr Pic_{C/S}\rightarrow \scr Pic_{\mc C/S}
\end{equation*}
which induces a morphism
\begin{equation}\label{Picmapb}
\pi ^*\colon\text{Pic}_{C/S}\rightarrow \text{Pic}_{\mc C/S}.
\end{equation}
This morphism is a monomorphism of group schemes over $S$.   Indeed suppose given two $S$-valued points $[L_1], [L_2]\in \text{Pic}_{C/S}(S)$ defined by line bundles $L_1$ and $L_2$ on $C$ such that $\pi ^*[L_1] = \pi ^*[L_2]$.  Then after possibly replacing $S$ by an \'etale cover, the two line bundles $\pi ^*L_1$ and $\pi ^*L_2$ on $\mc C$ are isomorphic.  Since $L_i = \pi _*\pi ^*L_i$ ($i=1,2$) this implies that $L_1$ and $L_2$ are also isomorphic.

\begin{lem}\label{cokernel} The cokernel $W$ of the morphism $\pi^*$ in (\ref{Picmapb}) is an \'etale group scheme over $S$.
\end{lem}

\begin{proof}
The space $W$ is the quotient of a smooth group scheme over $S$ by the action of a smooth group scheme, hence is smooth over $S$.  To show that $W$ is \'etale we use the infinitesimal lifting criterion. If $T_0\hookrightarrow T$ is a square zero nilpotent thickening defined by a sheaf of ideals $I\subset \mc O_T$, and if $\mc L$ is a line bundle on $\mc C_T$ whose reduction $\mc L_0$ on $\mc C_{T_0}$ is obtained by pullback from a line bundle on $C_{T_0}$, then we need to show that $\mc L$ is the pullback of a line bundle on $C_T$. Using the exponential sequence 
$$ 
0 \larr I \larr \cO^\times_T \larr \cO^\times _{T_0} \larr 1
$$
this amounts to the statement that the map
\begin{equation*}
\H^i(C_{T_0}, I\mc O_{C_T})\arr \H^i(\mc C_{T_0}, I\mc O_{\mc C_T})
\end{equation*}
is an isomorphism.
\end{proof}

\begin{lem} The cokernel $W$ is quasi-finite over $S$.
\end{lem}
\begin{proof} We may assume that $S$ is affine, say $S = \spec A$. Then the twisted curve $\cC$ is defined over some finitely generated $\ZZ$-subalgebra $A_{0}$, and formation of $W$ commutes with base change; hence we may assume that $A = A_{0}$. In particular $S$ is noetherian.

We may also assume that $S$ is reduced. After passing to a stratification of $S$ and to \'etale coverings, we may assume that $S$ is connected, that $C \arr S$ is of constant topological type, that the nodes of the geometric fibers of $C \arr S$ are supported along closed subschemes $S_{1}$, \dots,~$S_{k}$ of $C$ which map isomorphically onto $S$, and the reduced inverse image $\Sigma_{i}$ of $S_{i}$ in $\cC$ is isomorphic to the classifying stack $\cB\mmu_{m_{i}}$ for certain positive integers $m_{1}$, \dots,~$m_{k}$. Given a  line bundle $\mc L$ on $\mc C$, we can view the pullback of $\mc L$ to $\Sigma_{i}$ as a line bundle $L_{i}$ on $S_{i}$, with an action of $\mmu_{m_{i}}$. Such action is given by a character $\mmu_{m_{i}} \arr \GG_{\rmm}$ of the form $t \mapsto t^{a_{i}}$ for some $a_{i} \in \ZZ/(m_{i})$. By sending $L$ into the collection $(a_{i})$ we obtain a morphism $\pic_{\cC/S} \arr \prod_i\ZZ/(m_i)$, whose kernel is easily seen to be $\pic_{C/S}$. Therefore we have a categorically injective map $W \arr \prod_i\ZZ/(m_i)$; since $W$ and $\prod_i\ZZ/(m_i)$ are \'etale over $S$, this is an open embedding. Since $\prod_i\ZZ/(m_i)$ is noetherian and quasi-finite over $S$, so is $W$.
\end{proof}
\begin{rem}
Note that, in general, $W$ is non-separated. Consider $S$ the affine line over a field $k$ with parameter $t$, and consider the blowup $X$ of $\PP^1\times S$ at the origin. The zero fiber  has a singularity with local equation $uv=t$. There is a natural action of $\bmu_r$ along the fiber, which near the singularity looks like $(u,v) \mapsto (\zeta_r u, \zeta_r^{-1}v)$. The quotient $\cC = [X/\bmu_r]$ is a twisted curve, with twisted markings at the sections at 0 and $\infty$, and a twisted node. The coarse curve $C = X/\bmu_r$ has an $A_{r-1}$  singularity  $xy = t^r$, where $x = u^r$ and $y=v^r$. Let $E\subset \cC$ be one component of the singular fiber. The invertible sheaf $O_\cC(E)$ gives a generator of $W$, but it coincides with the trivial sheaf away from the singular fiber. In fact $W$ is the group-scheme obtained by gluing $r$ copies of $S$ along the open subset $t\neq 0$. 
\end{rem}

If $N$ denotes an integer annihilating $W$ we obtain a map
\begin{equation*}
\times N\colon\text{Pic}_{\mc C/S}\rightarrow \text{Pic}_{C/S}.
\end{equation*}

\begin{definition}\label{Def:Pic^o} Let $\text{Pic}^{\rmo}_{\mc C/S}$ denote the fiber product 
\begin{equation*}
\text{Pic}^{\rmo}_{\mc C/S} = \text{Pic}_{\mc C/S}\times _{\times N, \text{Pic}_{C/S}}\text{Pic}^{\rmo}_{C/S}.
\end{equation*}
\end{definition}
The open and closed subfunctor $\text{Pic}^{\rmo}_{\mc C/S}\subset \text{Pic}_{\mc C/S}$ classifies degree $0$-line bundles on $\mc C$. In particular, $\text{Pic}^{\rmo}_{\mc C/S}$ is independent of the choice of $N$ in the above construction. Note also that any torsion point of $\text{Pic}_{\mc C/S}$ is automatically contained in $\text{Pic}^{\rmo}_{\mc C/S}$ since the cokernel of
\begin{equation*}
\text{Pic}^{\rmo}_{C/S}\rightarrow \text{Pic}_{C/S}
\end{equation*}
is torsion free.

Since $\text{Pic}^{\rmo}_{C/S}$ is a semiabelian scheme over $S$, we obtain the following:

\begin{cor}\label{3.7} The group scheme $\text{\rm Pic}_{\mc C/S}^{\rmo}$ is an extension of a quasi-finite \'etale group scheme over $S$ by a semiabelian scheme.
\end{cor}

\section{Interlude: Relative moduli spaces}\label{interlude}

In this section we gather some results on relative moduli spaces which will be used in the verification of the valuative criterion for moduli stacks of twisted stable maps in the next section.

For a morphism of algebraic stacks $f\colon\cX \to \cY$, denote by
$\cI(\cX/\cY) = \operatorname{Ker} \left(\cI(\cX) \to f^*\cI(\cY)\right) = \cX \times_{\cX \times_\cY \cX} \cX$ the relative inertia stack.

Now let $f\colon\cX \to \cY$ be a morphism of algebraic stacks, locally of finite presentation, and assume the relative inertia $\cI(\cX/\cY) \to \cX$ is finite. Recall that when $\cY$ is an algebraic space we know that the morphism factors through the coarse moduli space: 
$\cX \to X \to \cY$, see \cite{Keel-Mori, Conrad}. We claim that there is a relative construction that generalizes to our situation:

\begin{theorem}\label{Th:rel-cms}
There exists an algebraic stack $X$,  morphisms $ \cX \stackrel{\pi}{\to} X\stackrel{\overline f}{\to} \cY$, and an isomorphism $\iota :f\rightarrow \overline f\circ \pi $ satisfying the following properties:
\begin{enumerate}
\item $\overline f\colon X \to \cY$ is representable,
\item  ($ \cX \stackrel{\pi}{\to} X\stackrel{\overline f}{\to} \cY, \iota )$ is initial for diagrams satisfying (1); namely if ($\cX \stackrel{\pi'}{\to} X'\stackrel{\overline {f'}}{\to} \cY, \iota ':f\simeq \overline {f'}\circ \pi ')$ has representable $\overline {f'}$ then there is a unique set of morphisms $(h\colon X \to X', \lambda :\pi '\simeq h\circ \pi , \tau :\overline f \simeq  \overline f' \circ h),$ such that the diagram
$$
\xymatrix{
f\ar[r]^\iota \ar[d]^-{\iota '}& \overline {f}\circ \pi \ar[d]^-\tau \\
\overline {f'}\circ \pi '& \overline {f'}\circ h\circ \pi \ar[l]_-{\lambda }}
$$
commutes.
\item $\pi$ is proper and quasi-finite,
\item $\cO_X \to \pi_*\cO_\cX$ is an isomorphism, and
\item if $X'$ is a stack and $X' \to X$ is a representable flat morphism, then $\cX \times_X X' \to X'\to \cY$ satisfies properties (1)-(4) starting from $\cX \times_X X' \to \cY$.
\end{enumerate}
\end{theorem}

\begin{definition} We call $\overline f:X{\to} \cY$ the relative coarse moduli space of $f\colon\cX \to \cY$.
\end{definition}

\begin{proof}
Consider a smooth presentation $R \double U$ of $\cY$. Denote by $X_R$ the coarse moduli space of $\cX_R = \cX \times_\cY R$ and $X_U$ the coarse moduli space of $\cX_U = \cX \times_\cY U$, both exist by the assumption of finiteness of relative inertia.
 
Since the formation of coarse moduli spaces commutes with flat base change, we have that $X_R = X_U \times_U R$. It is straightforward to check that $X_R \double U_R$ is a smooth groupoid. Denote its quotient stack $X := [X_R \double U_R]$. It is again straightforward to construct the morphisms $ \cX \stackrel{\pi}{\to} X\stackrel{\overline f}{\to} \cY$ from the diagram of relevant groupoids. 

Now $X\times_\cY U \to X_U$ is an isomorphism, $U\to \cY$ is surjective and $X_U \to U$ is representable, hence $X \to \cY$ is representable, giving (1). Similarly $\cX_U \to X_U$ is proper and quasi-finite hence $\cX \to X$ is proper and quasi-finite, giving (3). Also (4) follows by flat base-change to $X_U$. Part (5) follows again since formation of moduli spaces commutes with flat base change: the coarse moduli space of $\cX_U \times_X X'$ is $X_U\times_X X'$ and similarly the coarse moduli space of $\cX_R \times_X X'$ is $X_R\times_X X'$.

Now consider the situation in (2). Since $X' \to \cY$ is representable, it is presented by the groupoid $X' \times_\cY R \double X' \times_\cY U$, which, by the universal property of coarse moduli spaces, is the target of a canonical morphism of groupoids from $X_R \double U_R$, giving the desired morphism $h\colon X \to X'$.
\end{proof}

Now we consider the tame case. 

\begin{defn}\label{D:tamedef}We say that  a morphism $f\colon \cX \to \cY$ with finite relative inertia is {\em tame}  if for some algebraic space $U$ and faithfully flat $U \to \cY$  we have that $\cX \times_\cY U$ is a tame algebraic stack. This notion is independent of the choice of $U$ by \cite[Theorem 3.2]{AOV}. Let $\cX \to X \to \cY$ be the relative coarse moduli space.
\end{defn}

\begin{proposition} Assume $f\colon \cX \to \cY$  is tame.
If $X'$   is an algebraic stack and $X'\to X$ is a representable morphism of stacks,   then $X'\to \cY$ is the relative coarse moduli space of $\cX \times_X X' \to \cY$.
\end{proposition}

This is proven as in part (5) of the theorem, using the fact that formation of coarse moduli spaces of tame stacks commutes with arbitrary base change.

Now we consider a special case which is relevant for our study of tame stacks.

Suppose that $V$ and $S$ are strictly henselian local schemes, with a local morphism $V \arr S$. Let $\Gamma$ be a finite \lr group scheme over $S$ acting on $V$, fixing a geometric point $s \to V$ mapping to the closed point of $V$. Let $\cX = [V/\Gamma]$, and consider a morphism $f\colon \cX \to \cY$, with $\cY$ also tame.

Let $K_{s}\subset \Gamma_{s}$ be the subgroup scheme of the pullback $\Gamma_{s}$ fixing the composite object $s \to V \to \cY$ (so $K$ is the kernel of the map $\Gamma_{s} \rightarrow \text{Aut}_{\cY}(f(s))$).

\begin{lemma}
There exists a unique subgroup scheme $K \subseteq \Gamma$, flat over $S$, whose fiber over $s$ coincides with $K_{s}$.
\end{lemma}

\begin{proof}
Let $H_{\Gamma/S}$ the Hilbert scheme of $\Gamma$ over $S$, and let $\widetilde{H}_{\Gamma/S}$ be the closed subscheme parametrizing flat subgroup schemes of $\Gamma$. The subgroup scheme $K_{s} \subseteq \Gamma_{s}$ corresponds to a morphism $s\arr \widetilde{H}_{\Gamma/S}$; we need to show that this extends uniquely to a section $S \arr \widetilde{H}_{\Gamma/S}$. For this it is enough to show that $\widetilde{H}_{\Gamma/S}$ is \'etale over $S$. By Grothendieck's infinitesimal criterion, it is enough to prove the following fact: if $A$ is an artinian ring with algebraically closed residue field, $B$ is a quotient of $A$, $\Gamma_{A}$ is a linearly reductive group scheme over $A$, and $K_{B}$ is a flat subgroup scheme of the restriction $\Gamma_{B}$ of $\Gamma_{A}$ to $\spec B$, then there exists a unique flat subgroup scheme $K_{A}$ of $\Gamma_{A}$ that pulls back to $K_{B}$.

Let $\Gamma_{A}^{0}$ be the connected component of the identity in $\Gamma_{A}$, and set $H_{A} \eqdef \Gamma_{A}/\Gamma_{A}^{0}$. Then, by the proof of \cite[Lemma~2.16]{AOV}, $H_{A}$ is a tame \'etale finite group scheme over $\spec A$, $\Gamma_{A}^{0}$ is diagonalizable, and $\Gamma_{A} = H_{A} \ltimes \Gamma_{A}^{0}$. If $K_{B}^{0}$ is the connected component of the identity in $K_{B}$, then obviously $K^{0}_{B} = K_{B} \cap \Gamma^{0}_{B}$. If $L_{B} \eqdef K_{B}/K^{0}_{B}$, then $L_{B}$ is a subgroup scheme of $\Gamma_{B}$, and $K_{B} = L_{B} \ltimes K^{0}_{B}$. Since $H_{A}$ is \'etale over $\spec A$, the subgroup scheme $L_{B}$ extends uniquely to an \'etale subgroup scheme $L_{A} \subseteq H_{A}$; similarly, by Cartier duality, the flat subgroup scheme $K^{0}_{B} \subseteq \Gamma^{0}_{B}$ extends uniquely to a flat subgroup scheme $K^{0}_{A} \subseteq\Gamma^{0}_{A}$. The action of $L_{B}$ on $K^{0}_{B}$ extends to an action of $K_{A}$ on $K^{0}_{A}$, which is obtained by restriction from the action of $H_{A}$ on $\Gamma^{0}_{A}$. The flat subgroup scheme $K_{A} \eqdef L_{A} \ltimes K^{0}_{A} \subseteq \Gamma_{A}$ is the desired subgroup; this is clearly unique.
\end{proof}

Consider the geometric quotient $U = V/K$, the quotient group scheme $Q = \Gamma/K$ and let $\ocX= [U/Q]$. There is a  natural projection $q\colon\cX \to \ocX$.

\begin{proposition}\label{Prop:B3} We use the notation above.
There exists a unique morphism $g\colon \ocX \to \cY$ such that $f = g\circ q$, and this factorization of $f$ identifies $\ocX$ as the relative coarse moduli space of $\cX \to \cY$. 
\end{proposition}

\begin{remark}\label{relative-example}
The result does not hold when $\cX$ is not tame. 
Let $R$ be a strictly henselian local ring over $\mathbb{F}_p$ ($p$ a prime), and let $x\in R$ be a nonzero element of the maximal ideal which is not a zero divisor. Let 
$$
\alpha _{p, R}:(R-\text{algebras})\rightarrow (\text{abelian groups})
$$
be the functor given by
$$
(f:R\rightarrow A)\mapsto \{s\in A \mid s^p = 0\},
$$
where the group law is given by addition.  Note that $\alpha _{p, R}$ is represented by the scheme
$$
\Sp (R[T]/(T^p)),
$$
which is finite and flat over $\Sp (R)$.  We have a morphism of group schemes
$$
\rho :\alpha _{p, R}\rightarrow \alpha _{p, R}
$$
whose value over $(f:R\rightarrow A)$ is the map
$$
\{s\in A\mid s^p = 0\}\rightarrow \{s\in A\mid s^p=0\}, \ \ s\mapsto f(x)\cdot s.
$$
Let $\mc X = \mc Y = B\alpha _{p, R}$ be the classifying stack of $\alpha _{p, R}$ over $\Sp (R)$, and let
$$
f:B\alpha _{p, R}\rightarrow B\alpha _{p, R}
$$
be the map induced by $\rho $.  Notice that the map $\rho $ restricts to the zero homomorphism over the closed point of $\Sp (R)$, so with the above notation we have $K = \Gamma = \alpha _{p, R}$, and $\overline {\mc X} = \Sp (R)$.  On the other hand, the map $f$ does not factor through $\Sp (R)$ since the map $\rho $ is an isomorphism over the nonempty (since $x$ is not a zero divisor) open subset $\Sp (R_x)\subset \Sp (R)$.
\end{remark}
\begin{proof}
Consider the commutative diagram 
$$\xymatrix{
\cX \ar@/^2pc/[rr]^f \ar^{f'}[r]\ar[rd] & \cY_\cX\ar[r]^h\ar[d] & \cY\ar[d]\\
& X \ar[r]^{\overline f} & Y
}$$
where $X, Y$ are the coarse moduli spaces of $\cX$ and $\cY$, respectively,  the square is cartesian, and $\overline f$ is induced from $f$.

Since $h$ is fully faithful we may replace $\cY$ by $\cY_\cX$. As the formation of coarse moduli space of tame stacks commutes with arbitrary base change, $X$ is the coarse moduli space of $\cY_\cX$.  We may therefore assume that $X = Y$.

It also suffices to prove the proposition after making a flat base change on $X$. 
We may therefore assume that the residue field of $V$ is algebraically closed. We can write $\cY = [P/G]$, where $P$ is a finite $X$-scheme and $G$ is a finite flat split \lr group scheme over $X$. We can assume $P$ is local strictly henselian as well, and that the action of $G$ on $P$ fixes the closed point $p\in P$. 

The map $f$ induces, by looking at residual gerbes, a map $\cB\rho_0\colon\cB\Gamma_0 \to \cB G_0$ where $\Gamma_0$ and $G_0$ denote the corresponding reductions to the closed point of $X$. This is induced by a group homomorphism $\rho_0\colon\Gamma_0 \to  G_0$,  uniquely defined up to conjugation. Since $\Gamma$ and $G$ are split, the homomorphism $\rho_0$ lifts to a homomorphism $\rho \colon\Gamma\to G$.

We obtain a 2-commutative diagram 
$$\xymatrix{
\cB\Gamma_0 \ar[r]^{\cB\rho_0} \ar[d]_{\bar s} & \cB G_0 \ar[d]^{\bar p} \\
[V/\Gamma] \ar[r]^f & [P/G].
}$$

Let $P'\to [V/\Gamma]$ be the pullback of the $G$-torsor $P \to [P/G]$

\begin{lemma}
There is an isomorphism of $G$-torsors over $[V/\Gamma]$
$$[V\times G/\Gamma] \to P'.$$
\end{lemma}

\begin{proof}
Note that $ V \times _{[V/\Gamma]} \cB\Gamma_0 = s$ and $P'\times _{[V/\Gamma]} \cB\Gamma_0 = \cB\Gamma_0 \times_{\cB G_0}p$.
Then the commutativity of the diagram 
$$\xymatrix{
s\ar[r]\ar[d] & p \ar[d] \\ \cB\Gamma_0 \ar[r]& \cB G_0
}$$
implies that there is an isomorphism of $G_0$-torsors over $\cB\Gamma _0$
$$\xymatrix {([V \times G/\Gamma] \times_{[V/\Gamma]} \cB \Gamma_0) \ar[rr]^{\sim} &&P'\times_{[V/\Gamma]} \cB \Gamma_0 .}$$
As in the proof of \cite[3.6]{AOV}, this lifts to an isomorphism as required.
\end{proof}
Now we have a morphism 

$$\xymatrix { V \ar[r]\ar@/^2pc/[rrr]^{\widetilde f} & [V \times G/\Gamma] \ar[r]& P'\ar[r] & P}$$
compatible with the actions of $\Gamma$ and $G$. By the universal property of the quotient, the map $\widetilde f$ factors through $U$. Passing to the quotient by the respective group actions, we get morphisms 
$$[V / \Gamma] \to [U/Q] \to [P/G].$$

 Let $\cX \to \cZ \to \cY$ be the relative coarse moduli space. Since $Q$ injects in $G$, the morphism $ [U/Q] \to [P/G]$ is representable, and therefore we have a morphism $\cZ \to  [U/Q]$ over $\cY$.

 To check this is an isomorphism it suffices to check after base change along the flat morphism $P \to \cY$. In other words, we need to check that  ${[V/\Gamma]} \times _{[P/G]} P\to [U/Q] \times _{ [P/G]} P$ is the coarse moduli space. This map can be identified with $[V\times G / \Gamma] \to U\times^Q G$. This follows by noting that 
 $$(\cO_V \otimes \cO_G)^\Gamma = ( (\cO_V)^K\otimes \cO_G)^Q = (\cO_U \otimes \cO_G)^Q.$$
 
\end{proof}

\section{Twisted stable maps}\label{Section:3}

This section relies on the results of Appendix \ref{AppendixA}, in which the stack of twisted curves is constructed.

\begin{definition} Let $\cM$ be a finitely presented 
%tame
 algebraic stack with finite inertia proper over a scheme $S$.  A \emph{twisted stable map} from an $n$-marked twisted curve $(f\colon\mc C\rightarrow S, \{\Sigma _i\subset \mc C\}_{i=1}^n)$ over $S$ to $\cM $ is a representable morphism of $S$-stacks
$$
g\colon\mc C\rightarrow \cM
$$
such that the induced morphism on coarse spaces
$$
\overline{g}\colon C\rightarrow M
$$
is a stable map (in the usual sense of Kontsevich) from the $n$-pointed curve $(C, \{p_1, \dots, p_n\}_{i=1}^n)$ to $M$ (where the points $p_i$ are the images of the $\Sigma _i$). 

If $M$ is projective over a field, the {\em target class} of $g$ is the class of $\overline{g}_*[C]$ in the chow group of 1-dimensional cycles up to numerical equivalence of $M$.
\end{definition}

\begin{proposition}\label{Prop:existence}
Let $\cM$ be a finitely presented algebraic stack proper over a
scheme $S$ with finite inertia. Then the stack $\cK_{g,n}(\cM)$ of twisted stable maps
from $n$-pointed genus $g$ curves into $\cM$  is an algebraic stack locally finitely presented
 over $S$. 
\end{proposition}

\begin{proof}
The statement is local in the Zariski topology of $S$, therefore we may assume $S$ affine. Since $\cM$ is of finite presentation, it is obtained by base change from the spectrum of a noetherian ring, and replacing $S$ by this spectrum we may assume $\cM$ is of finite type over a noetherian scheme $S$. 

As $\cM$ is now of finite type, there is an integer $N$ such that the degree of the automorphism group scheme of any geometric point of $\cM$ is $\leq N$. 

Consider the stack of twisted curves $\fM_{g,n}^{\tw}$ defined in \ref{3.19}. It contains an open substack $\fM_{g,n}^{\tw,\leq N}\subset\fM_{g,n}^{\tw}$ of twisted curves of index bounded by $N$ by \ref{1.7}. The natural morphism $\fM_{g,n}^{\tw,\leq N} \to \fM_{g,n}$ to the stack 
of prestable curves is  quasi-finite by \ref{1.7}.

Consider $\fB = \fM_{g,n}^{\tw,\leq N}\times_{\spec \ZZ} S$.  It is an algebraic stack of finite type over $B = \fM_{g,n} \times_{\spec \ZZ} S$. With the notation of Theorem~\ref{athm} in mind,
denote   $Y := M \times_{S}B$,  $Y_\fB := M \times_{S}\fB$ and similarly $\cY_\fB := \cM \times_{S}\fB$. These are all algebraic stacks.

Finally consider the universal curves $\fC^{\tw} \to \fM_{g,n}^{\tw}$ and $\fC\to \fM_{g,n}$, and their corresponding pullbacks $\cX_\fB:=\fC^{\tw}_{\fB}$, $X_\fB:=\fC_{\fB}$ and $X:=\fC_{B}$.

We consider the stack $$\hom_{\fB}^{\mathrm{rep}}(\cX_{\fB},\cY_{\fB}),$$ parametrizing representable morphisms of the universal twisted curve to $\cM$, as discussed in Appendix \ref{AppendixB}. Note that by definition of twisted curves, the morphism $\cX_{\fB} \to {\fB}$ is flat and proper, in particular it has finite diagonal. By  Theorem~\ref{athm}(i)   the stack is algebraic and locally of finite presentation over $\fB$. 

The Kontsevich moduli space $\overline{\cM}_{g,n}(M)$ is open in $\hom_{B}(X,Y)$. We  have 
   \[
   \cK_{g,n}(\cM) \ \ \simeq\ \   \hom_{\fB}^{\mathrm{rep}}(\cX_{\fB},\cY_{\fB}) \times_{\hom_{B}(X,Y)}\overline{\cM}_{g,n}(M),
   \]
hence $\cK_{g,n}(\cM) \to \overline{\cM}_{g,n}(M)$ is locally of finite presentation. The result follows.
\end{proof}

%As discussed in Appendix \ref{AppendixB},   there is an algebraic stack $$\hom_{\fB}^{\mathrm{rep}}(\fC^{\tw}_{\fB},\cM_{\fB}),$$ locally of finite type over $\fB$, parametrizing representable morphisms of the universal twisted curve to $\cM$. 

The main result of this paper is the following:

\begin{theorem}\label{Th:3.1}
Let $\cM$ be a finitely presented tame algebraic stack proper over a
scheme $S$ with finite inertia. Then the stack $\cK_{g,n}(\cM)$ of twisted stable maps
from $n$-pointed genus $g$ curves into $\cM$  is a locally finitely presented
algebraic stack over $S$, which is proper and quasi-finite over the stack of stable
maps into $M$. If $M$ is projective over a field, the substack $\cK_{g,n}(\cM,\beta)$   of twisted stable maps
from $n$-pointed genus $g$ curves into $\cM$ with target class $\beta$ is proper, and is open and closed in $\cK_{g,n}(\cM)$.\end{theorem}

\begin{proof}

Continuing with the notation of the proposition, we have that the stack
$\hom_{\fB}^{\mathrm{rep}}(\fC^{\tw}_{\fB},\cM_{\fB})$ is algebraic and locally of finite presentation over $S$.

We also have the analogous $\hom_{\fB}(\fC_{\fB},M_{\fB})$ and $\hom_{B}(\fC_{B},M_{B})$. We have the natural morphisms
\begin{equation}\label{compositemap}
   \xymatrix{
   \hom_{\fB}^{\mathrm{rep}}(\fC^{\tw}_{\fB},\cM_{\fB}) \ar[r]^a&
   \hom_{\fB}(\fC_{\fB},M_{\fB}) \ar[r]^b&
   \hom_{B}(\fC_{B},M_{B}).}
\end{equation}

 \begin{prop}\label{5.2} The morphism $a$ is quasi-finite.
 \end{prop}
 \begin{proof}
 It suffices to prove the following.  Let $k$ be an algebraically closed field with a morphism $\Sp (k)\rightarrow S$ and $\mc C/k$ a twisted curve with coarse moduli space $C$.  Let $f\colon C\rightarrow M$ be an $S$-morphism, and let $\mls G\rightarrow \mls C$ denote the pullback of $\mc M $ to $\mc C$ along the composite
 $$
 \xymatrix{\mc C\ar[r]& C\ar[r]^f& M.}
 $$
 We then need to show that the groupoid $\Sec (\mc G/\mc C)(k)$ (notation as in \ref{C.30}) has only finitely many isomorphism classes of objects.
 
 We first reduce to the case where $C$ is smooth. Let $C'$ denote the normalization of $C$ (so $C'$ is a disjoint union of smooth curves), and let $\mc C'\rightarrow \mc C$ be the maximal reduced substack of the fiber product $\mc C\times _CC'$.  For any node $x\in C$ let $\Sigma _x$ denote maximal reduced substack of $\mc C_x$ (so $\Sigma _x$ is isomorphic to $B\mu _r$ for some $r$), and let $\Sigma $ denote the disjoint union of the $\Sigma _x$ over the nodes $x$.  Then there are two natural inclusions $j_1, j_2\colon\Sigma \hookrightarrow \mc C'$, and by \cite[A.0.2]{AGV}  the functor
$$
\begin{CD}
\Hom(\mc C, \mc M)(k)\\
@VVV \\
 \Hom(\mc C', \mc M)(k)\times _{j_1^*\times j_2^*, \Hom(\Sigma , \mc M)(k)\times \Hom(\Sigma , \mc M)(k), \Delta }\Hom(\Sigma , \mc M)(k)
 \end{CD}
$$
is an equivalence of categories.  Since the stack $\Hom (\Sigma, \mc M)$ clearly has quasi-finite diagonal, it follows that the map 
$$
\Sec (\mls G/\mls C)\rightarrow \Sec (\mls G\times _{\mls C}\mls C'/\mls C')
$$
is quasi-finite.

 It follows that to prove \ref{5.2} it suffices to consider the case when $\mc C$ is smooth.
 
Let $\mc Y$ denote the normalization of $\mc G$.  Then any section $\mc C\rightarrow \mc G$ factors uniquely through $\mc Y$ so it suffices to show that the set
$$
\Sec (\mc Y/\mc C)(k)
$$
is finite.  To prove this we may without loss of generality assume that there exists a section $s\colon\mc C\rightarrow \mc Y$.

For any smooth morphism $V\rightarrow \mc C$ the coarse space of the fiber product $\mc G_V:=\mc G\times _{\mc C}V$ is equal to $V$ by \cite[Corollary 3.3 (a)]{AOV}. For a smooth morphism $V\rightarrow \mc C$, let 
$$
\xymatrix{\mc Y_V\ar[r]^c& \overline {\mc Y}_V\ar[r]^d& V}
$$
be the factorization of $\mc Y_V\rightarrow V$ provided by \cite[2.1]{bounded} (rigidification of the generic stabilizer).  The section $s$ induces a section $V\rightarrow \overline {\mc Y}_V$, and hence by \cite[2.9 (ii)]{homstack} the projection $d\colon\overline {\mc Y}\rightarrow V$ is in fact an isomorphism.  It follows that $\mc Y$ is a gerbe over $\mc C$.  If $\mathbb{G}\rightarrow \mc C$ denotes the automorphism group scheme of $s$, then in fact $\mc Y = \cB _{\mc C}\mathbb{G}$.

Let $\Delta \subset \mathbb{G}$ denote the connected component of $\mathbb{G}$, and let $\mathbb{H}$ denote the \'etale part of $\mathbb{G}$, so we have a short exact sequence 
 $$
 1\rightarrow \Delta \rightarrow \mathbb{G} \rightarrow \mathbb{H} \rightarrow 1
 $$
 of group schemes over $\mc C$.    
 
 Let $\mc C'\rightarrow \mc C$ be a finite \'etale cover such that $\mathbb{H}$ and the Cartier dual of $\Delta $ constant.  If $\mc C^{\prime \prime }$ denotes the product $\mc C'\times _{\mc C}\mc C'$ and $\mc Y'$ and $\mc Y^{\prime \prime }$ the pullbacks of $\mc Y$ to $\mc C'$ and $\mc C^{\prime \prime }$ respectively, then one sees using  a similar argument to the one used in \cite[\S 3]{homstack} that it suffices to show that
 $$
 \Sec (\mc Y'/\mc C')(k) \ \ \text{and} \ \ \Sec (\mc Y^{\prime \prime }/\mc C^{\prime \prime })(k)
 $$
 are finite. We may therefore assume that $\mathbb{H}$ and the Cartier dual of $\Delta $ are constant.
 
 The map $\mathbb{G}\rightarrow \mathbb{H}$ induces a projection over $\mc C$
 $$
 \cB _{\mc C}\mathbb{G}\rightarrow \cB _{\mc C}\mathbb{H}.
 $$
 Giving a section $z\colon\mc C\rightarrow \cB _{\mc C}\mathbb{H}$ is equivalent to giving an $\mathbb{H}$-torsor $P\rightarrow \mc C$.  Since $\mc C$ is normal and $\mathbb{H}$ is \'etale, such a torsor is determined by its restriction to any dense open subspace of $\mc C$.  Since $\mc C$ contains a dense open subscheme and the \'etale fundamental group of such a subscheme is finitely generated it follows that
$$
\Sec (\cB _{\mc C}\mathbb{H}/\mc C)(k)
$$
is finite.

Fix  a section $z\colon\mc C\rightarrow \cB _{\mc C}\mathbb{H}$, and let $\mc Q\rightarrow \mc C$ be the fiber product
$$
\mc Q:= \mc C\times _{z, \cB _{\mc C}\mathbb{H}}\cB _{\mc C}\mathbb{G}.
$$
Then $\mc Q$ is a gerbe over $\mc C$ banded by a twisted form of $\Delta $.  By replacing $\mc C$ by a finite \'etale covering as above, we may assume that in fact $\mc Q$ is a gerbe banded by $\Delta $ and hence isomorphic to $\mc C\times _{\Sp (k)}\cB \Delta _0$ for some finite diagonalizable group scheme $\Delta _0/k$.  Writing $\Delta _0$ as a product of $\mu _r$'s, we even reduce to the case when $\Delta _0 = \mu _r$.  In this case, as explained in \ref{B.18}, giving a section $w\colon\mc C\rightarrow \cB _{\mc C}\Delta $ is equivalent to giving a pair $(\mls L, \iota )$, where $\mls L$ is  a line bundle on $\mc C$ and $\iota \colon\mc O_{\mc C}\rightarrow \mls L^{\otimes r}$ is an isomorphism. We conclude that the set of sections $\mc C\rightarrow \mc Q$ is in bijection with the set
$$
\underline {\text{Pic}}_{\mc C/k}[r](k)
$$
of $r$-torsion points in the Picard space $\underline {\text{Pic}}_{\mc C/k}$. By \ref{3.7} this is a finite group, and hence this completes the proof of \ref{5.2}.
 \end{proof}
 
 Consider again the diagram \ref{compositemap}.
The morphism $a$ is quasi-finite by \ref{5.2} and of finite type by \ref{repfinite}. The second morphism is quasi-finite and of finite type as it is obtained by base change.   Therefore the composite \ref{compositemap} is also quasi-finite and of finite type.
It follows as in Proposition \ref{Prop:existence}
%, the Kontsevich moduli space $\overline{\cM}_{g,n}(M)$ is open in $\hom_{B}(\fC_{B},M_{B})$. We  have 
%   \[
%   \cK_{g,n}(\cM) \ \ \simeq\ \   \hom_{\fB}^{\mathrm{rep}}(\fC^{\tw}_{\fB},\cM_{\fB}) \times_{\hom_{B}(\fC_{B},M_{B})}\overline{\cM}_{g,n}(M),
%   \]
that $\cK_{g,n}(\cM) \to \overline{\cM}_{g,n}(M)$ is quasi-finite and of finite type.

Properness now follows from the valuative criterion, which is the content of the following Proposition \ref{Prop:valuative}.
\end{proof}

\begin{proposition}\label{Prop:valuative}
Let $R$ be  a discrete valuation ring, with spectrum $T$ having generic point $\eta$ and special point $t$. Fix a morphism $T \to S$ and let $\cC_{\eta} \to \cM$ be  an $n$-pointed  twisted stable map  over the generic point of $T$, and $C \to M$ a stable map over $T$ extending the coarse map $C_{\eta} \to M$. Then there is a discrete valuation ring $R_{1}$ containing $R$ as a local subring, and corresponding morphism of spectra $T_{1}\to T$, and  an $n$-pointed twisted stable map $\cC_{T_{1}} \to \cM$, such that the restriction $\cC_{T_{1}}\mid _{\eta_{1}} \to \cM$ is isomorphic to the base change $\cC_{\eta}\times_{T}T_{1} \to \cM$, and the coarse map coincides with $C\times_{T}T_{1} \to M$. 

Such an extension, when it exists, is unique up to a unique isomorphism.
\end{proposition}

Before proving this Proposition, we need to extend the following basic result known in case $\cM$ is a Deligne--Mumford stack.

\begin{lemma}[Purity Lemma] \label{Lem:purity}
Let $\cM$ be a separated tame stack with coarse moduli space $M$. Let $X$ be a locally noetherian separated scheme of dimension 2 satisfying Serre's condition S2. Let $P\subset X$ be a finite subset consisting of closed points, $U = X \setminus P$. Assume that the local fundamental groups of $U$ around the points of $P$ are trivial and that the local Picard groups of $U$ around points of $P$ are torsion free. 

Let $f\colon X \to M$ be a morphism. Suppose there is a lifting $\widetilde f_{U}\colon U \to \cM$: 
   \[
   \xymatrix{
   && \cM\ar[d]\\
   U\ar[r]\ar^{\widetilde f_{U}}[rru] & X \ar^{f}[r]& M 
   }
   \]
Then the lifting extends to $X$:
   \[
   \xymatrix{
   && \cM\ar[d]\\
   U\ar[r]\ar^{\widetilde f_{U}}[rru] & X \ar^{f}[r]\ar@{.>}|-{\widetilde f}[ru]& M 
   }
   \]
The lifting $\widetilde f$ is unique up to a unique isomorphism.
\end{lemma}
\begin{remark}
We will use this lemma principally in the following cases:
\begin{enumerate}
\item $X$ is regular,
\item $X$ is normal crossings, locally $\Sp (R[x,y]/xy)^{\sh}$, with $R$ regular, or
\item $X = X_0 \times G$, with $X_0$ one of the first two cases and $G$ a linearly reductive finite group scheme. 
\end{enumerate}
These schemes are S2 since they are l.c.i. \cite[Proposition 18.13]{Eisenbud}, and the local fundamental and Picard groups are easily seen to be trivial. 
\end{remark}
\begin{proof}
The question is local in the \'etale topology, so we may replace $X$ by its strict henselization over some point $p\in P$, and correspondingly we may replace $M$ and $\cM$ by the strict henselization at $f(p)$. 

By \cite[Theorem 3.2 (d)]{AOV} we can write $\cM =[V/G]$, where $V\to M$ is a finite morphism and $G\to M$ a linearly reductive group scheme acting on $V$. By  \cite[Lemma 2.20]{AOV} we have an exact sequence 
  $$
  1 \arr \Delta \arr G \arr H \arr 1
  $$
of group schemes over $M$, where $\Delta$ is diagonalizable and $H$ is tame and \'etale.  

The morphism $\widetilde f_{U}$ is equivalent to the data of
 a $G$-torsor $P_{U} \to U$ and a $G$-equivariant morphism $P_U \to V$. We first wish to extend $P_{U}$ over $X$.
 
Consider the $H$-torsor $Q_{U}\eqdef P_{U}/\Delta \to U$. Since the local fundamental groups of $U$ around $P$ are trivial, this $H$ torsor is trivial, and extends trivially to an $H$ torsor $Q \to X$. Note that $Q$ is S2 since it is \'etale over the S2 scheme $X$ (for example by \cite[IV.6.4.2]{EGA}).

Now $P_{U} \to Q_{U}$ is a $\Delta$-torsor. We claim that it extends uniquely to a $\Delta$-torsor $P \to Q$. Since $\Delta$ is diagonalizable, it suffices to treat the case $\Delta= \mmu_{r}$. In this case $P_{U} \to Q_{U}$ corresponds to an $r$-torsion line bundle with a chosen trivialization of its $r$-th power. The line bundle extends to a trivial line bundle on the strictly Henselian $X$ by the assumption on the local Picard groups, and the trivialization, which is a regular function on $X$, extends by the S2 assumption \cite[IV.5.10.5]{EGA}.  Note that the extension $P$ is S2, since it is flat over the S2 scheme $Q$ with S2 geometric fiber $\Delta$ \cite[IV.6.4.2]{EGA}.

Next we need to extend the principal action of $G$ on $P_{U}$ to $P$. The closure $\Gamma$ of the graph of the morphism $G \times_{M} P_{U} \to P_{U}$ inside the scheme $ G \times_{M} P \times_{X} P$ is finite over $G \times_{M} P$, and an isomorphism over $G \times_{M} P_{U}$. As $G \times_{M} P$ is finite and flat over the S2 scheme $P$ with S2 fibers, it is also S2 \cite[IV.6.4.2]{EGA}. It follows that the morphism $\Gamma \to G \times_{M} P$ is an isomorphism, hence the action extends to a morphism $G \times_{M} P \to P$. It is easy to check that this is a principal action, as required.

The same argument shows that the equivariant morphism $P_{U} \to V$ extends to an equivariant morphism $P \to V$, as required.
\end{proof}

\begin{proof}[Proof of Proposition~\ref{Prop:valuative}]
{\sc Step 1:} We extend $\cC_{\eta} \to \cM$ over the generic points of the special fiber $C_{s}$. 

The question is local in the \'etale topology of $C$, therefore we can replace $C$ by its strict henselization at one of the generic points of $C_{s}$. We may similarly replace $M$ by its strict henselization at the image point, and hence assume it is of the form $[V/G]$, where $V \to M$ is finite and $G$ a locally well split group scheme, see \cite[Theorem 3.2 (d) and Lemma 2.20]{AOV}. We again write $G$ as an extension  
  $$
  1 \arr \Delta \arr G \arr H \arr 1.
  $$
The morphism $\cC_{\eta}= C_{\eta} \to \cM$ is equivalent to the data of
 a $G$-torsor $P_{\eta} \to C_{\eta}$ and a $G$-equivariant morphism $P_{\eta} \to V$. By Abhyankar's lemma, the $H$-torsor $Q_{\eta}\eqdef P_{\eta}/\Delta \to C_{\eta}$ extends uniquely to $Q \to C$ after a base change on $R$. We want to extend the $\Delta$-torsor $P_{\eta}\to Q_{\eta}$ to a $\Delta$-torsor $P\to Q$. Factoring $\Delta$, we may assume $\Delta = \mmu_{r}$.  In this case $P_{\eta} \to Q_{\eta}$ corresponds to an $r$-torsion line bundle with a chosen trivialization of its $r$-th power. This line bundle is trivial since we have localized, hence it trivially extends to $Q$. The extension of the trivializing section may have a zero or pole along $C_{s}$, but after an $r$-th order base change on $R$ this zero or pole has multiplicity $r\cdot l$ for some integer $l$. Twisting the line bundle by $\cO_{C}(l C_{s})$, the trivialization extends uniquely over $C_{s}$ as required.
 
 We need to lift the principal action of $G$ on $P$. Write again $\Gamma$ for the closure of the graph of $G \times_{T} P_{\eta} \to P_{\eta}$ inside $G \times_{T} P_{\eta} \times_{C} P_{\eta}$. This is a $\Delta$-equivariant subscheme, with respect to the action $\delta\cdot (g,p_{1},p_{2}) = (\delta g,p_{1},\delta p_{2})$. The action is free, as it is already free on the last factor. Therefore the projection $\Gamma \to \Gamma/\Delta$ is a $\Delta$-torsor. The quotient $\Gamma/\Delta$ includes the subscheme $(\Gamma/\Delta)_{\eta}= Q_{\eta}$. We claim that this is scheme-theoretically dense in $\Gamma/\Delta$, because any nonzero sheaf of ideals with support over $s$ would pull back to a nonzero sheaf of ideals on $\Gamma$. But $\Gamma_{\eta}$ is by definition scheme theoretically dense in $\Gamma$. Since $Q$ is normal, it follows that $\Gamma/\Delta=Q$. Now as both $\Gamma$ and $P$ are $\Delta$-torsors over $Q$ and $\Gamma\to P$ is $\Delta$-equivariant, this morphism is an isomorphism. Hence we have a morphism $G \times P \to P$, extending the action. It is again easy to see this is a principal action with quotient $C$. It is also clearly unique.

Exactly the same argument shows that the morphism $P_{\eta} \to V$ extends uniquely to a $G$-equivariant $P \to V$, as required.

{\sc Step 2:} We  extend the twisted stable map $C \to  \cM$ over the general locus $C_{\mathrm{gen}}$, simply by applying the purity lemma (Lemma \ref{Lem:purity}). Uniqueness in the Purity Lemma implies that the extension is unique up to unique isomorphism.

{\sc Step 3:} We extend the twisted curve $\cC$ and the twisted stable map $\cC \to \cM$ along the smooth locus of $C$.

Consider the index $r_{i}$ of $\cC_{\eta}$ over a marking $\Sigma^{C}_{i}$ of $C$. Let $\cC_{\mathrm{sm}}$ be the stack obtained by taking $r_{i}$-th root of $\Sigma^{C}_{i}$ on $C$ for all $i$. Then $(\cC_{\mathrm{sm}})_{\eta}$ is uniquely isomorphic to $(\cC_{\eta})_{\mathrm{sm}}$.

We need to construct a representable map $\cC_{\mathrm{sm}} \to \cM$ lifting $C_{\mathrm{sm}} \to M$. The problem is local in the \'etale topology of $\cC_{\sm}$, so we may present $\cC_{\sm}$ at a point $p$ on $\Sigma^{C}_{i}\cap C_{s}$ as $[D/\mmu_{r_i}]$, where $D$ is smooth over $T$, with unique fixed point $q$ over $p$. Applying the purity lemma to the map $D\setminus \{q\} \to \cM$, we have a unique extension $D \to \cM$. The uniqueness applies also for $\mmu_{r_i}\times D \to \cM$, which  implies that the object $D \to \cM$ is $\mmu_{r_i}$-equivariant, giving a  unique morphism $[D/\mmu_{r_i}]\to \cM$, as needed.

{\sc Step 4:} We extend the twisted curve $\cC$ and the twisted stable map $\cC \to \cM$ over the closure of the singular locus of $C_{\eta}$.

This is similar to step 3.

{\sc Step 5:} Extension over isolated singular points.

Consider an isolated  node $p$ of $C$. Passing to an extension on $R$ we may assume it is rational over the residue field $k$. Its strict henselization is isomorphic to $(\spec R[u,v]/(uv-\pi^{l}))^{\sh}$, where $\pi$ is a uniformizer of $R$, and by Remark \ref{curveconstruct} there is a twisted curve $\cC_{l}$ which is regular over $p$ with index $l$. This twisted curve has a chart of type $[D/\mmu_{l}]$, where the strict henselization of $D$ looks like $\spec R[x,y]/(xy-\pi)$, where $\mmu_{l}$ acts via $(x,y) \mapsto (\zeta_{l}x,\zeta_{l}^{{-1}}y)$, and $u = x^{l}, v= y^{l}$. There is a unique point $q$ of $D$ over $p$.

To construct $\cC_{l}\to \cM$ it suffices to apply the purity lemma to the map $D \setminus \{q\} \to \cM$; the resulting extension $D \to \cM$ is  $\mmu_{l}$-equivariant by applying the Purity Lemma to $\bmu_l \times D \to \cM$.

We need to replace the morphism $\cC_{l}\to \cM$, as it is not necessarily representable, and the construction of $\cC_{l}$ does not commute with base change. 

Consider the morphism $\cC_{l} \to C \times \cM$. This morphism is proper and quasi-finite. Let  $\cC_{l} \to \cC \to C \times \cM$ be the relative moduli space. By Proposition \ref{Prop:B3}, $\cC$ is a twisted curve, hence $\cC \to \cM$ is a twisted stable map. Further, its formation commutes with further base change since the formation of moduli space does.
\end{proof}

\section{Reduction of Spaces of Galois admissible covers}\label{reductionsection}

Let $S$ be a scheme and $G/S$ a linearly reductive finite flat group scheme. Its classifying stack  $\cB G$ is tame by \cite[Theorem 3.2]{AOV}. We can then consider the proper $S$-stack $\cK_{g, n}(\cB G)$ of twisted stable maps from genus $g$ twisted curves with $n$ marked points to the classifying stack $\cB G$.  This is the stack which to any $S$-scheme $T$ associates the groupoid of data
$$
(\mc C, \{\Sigma _i\}_{i=1}^n, h\colon \mc C\rightarrow \cB G),
$$
where $(\mc C, \Sigma _i\}_{i=1}^n$ is an $n$-marked twisted curve such that the coarse space $(C, \{p_i\}_{i=1}^n)$ (with the marking induced by the $\Sigma _i$) is a stable $n$-pointed curve in the sense of Deligne-Mumford-Knudsen, and $h$ is a representable morphism of stacks over $S$.
 This is equivalent to the data
$$
(\mc C, \{\Sigma _i\}_{i=1}^n, P\to \mc C),
$$
where $P\to \cC$ is a principal $G$-bundle, with $P$ representable. When $G$ is tame and \'etale, $P\to C$ is an admissible $G$-cover, and the stack $\cK_{g, n}(\cB G)$ has been studied extensively. See \cite{ACV} for discussion and further references.

The main result of this section is the following:

\begin{thm}\label{flattheorem} The stack $\cK_{g, n}(\cB G)$ is flat over $S$, with \lci fibers.
\end{thm}
Here we say that an algebraic stack $\cX$ is \lci if it admits a smooth surjective morphism from a \lci algebraic space. This is equivalent to the condition that locally $\LL_\cX$ be a perfect complex supported in degrees $[-1,1]$.

In particular, it follows that $\cK_{g, n}(\cB G)$ is Cohen--Macaulay if $S$ is Cohen--Macaulay, and Gorenstein if $S$ is Gorenstein.

The smoothness in characteristic 0 does not extend: the example of $\cK_{1,1}(\cB \bmu_2^2)$ in Section~\ref{Sec:ex11} is evidently not smooth. This may not be too surprising since the group scheme $G/S$ itself is in general not smooth.

Given a scheme $X$ of finite type over a field $k$ and an extension $k'$ of $k$, it is well known that $X$ is a \lci if and only if $X_{k'}$ is; hence the statement of Theorem~\ref{flattheorem} is local in the fpqc topology on $S$.

To prove this theorem we in fact prove a stronger result.

Let $\mc M_{g, n}^{\text{tw}}(G)$ denote the stack over $S$ which associates to any $S$-scheme $T$ the groupoid of data 
$$
(\mc C, \{\Sigma _i\}_{i=1}^n, h\colon \mc C\rightarrow \cB G),
$$
where $(\mc C, \{\Sigma _i\}_{i=1}^n)$ is an $n$-marked twisted curve over $T$ and $h$ is a morphism of stacks over $S$ (so $(C, \{p_i\}_{i=1}^n)$ is not required to be stable and $h$ is not necessarily representable).

\begin{lem} The natural inclusion $\cK_{g, n}(\cB G)\subset \mc M_{g, n}^{\text{\rm tw}}(G)$ is representable by open immersions.
\end{lem}
\begin{proof} This is because the condition that an $n$-pointed nodal curve is stable is an open condition as is the condition that  a morphism of stacks $\mc C\rightarrow \cB G$ is representable, see Corollary \ref{repfinite} and the discussion preceding it.
\end{proof}

To prove \ref{flattheorem} it therefore suffices to show that $\mc M_{g, n}^{\text{tw}}(G)$ is flat over $S$ with local complete intersection fibers.  For this we in fact prove an even stronger result.  Let $\mc M_{g, n}^{\text{tw}}$ denote the stack defined in \ref{3.19}.  Forgetting the map to $\cB G$ defines a morphism stacks 
\begin{equation}\label{flatmap}
\mc M_{g, n}^{\text{tw}}(G)\rightarrow \mc M_{g, n}^{\text{tw}}\times _{\Sp (\Z)}S.
\end{equation}

Since $\mc M_{g, n}^{\text{tw}}$ is smooth over $\Z$ by \ref{3.19} the following theorem implies \ref{flattheorem}.
\begin{thm}\label{flattheorembis}
The morphism \ref{flatmap} is a flat morphism of algebraic stacks, with \lci fibers.

Equivalently: let $T$ be a scheme and $\cC \to T$ a twisted curve. Then $\underline{\text{\rm Hom}}_T(\cC, \cB G)$ is flat over $T$, with \lci fibers.
\end{thm}

The proof of \ref{flattheorembis} occupies the remainder of this section.

Note first of all that $\mc M_{g, n}^{\text{tw}}(G)$ is an algebraic stack locally of finite presentation over $\mc M_{g, n}^{\text{tw}}$.  Indeed the stack $\mc M_{g, n}^{\text{tw}}(G)$ is equal to the relative Hom-stack
$$
\underline {\text{Hom}}_{\mc M_{g, n}^{\text{tw}}}(\mc C^{\text{univ}}, \cB G\times \mc M_{g, n}^{\text{tw}}),
$$
where $\mc C^{\text{univ}}\rightarrow \mc M_{g, n}^{\text{tw}}$ denotes the universal twisted curve.  This stack is  locally of finite type over $\mc M_{g, n}^{\text{tw}}$ by \cite[1.1]{homstack} and Lemma \ref{Lem:relative-algebraic}. 
%It is also shown there that a $\text{Hom}$ stack is algebraic whenever the base is representable. This extends to the case where the base is an algebraic stack by the following well known result:

We prove \ref{flattheorembis} by first studying two special cases and then reducing the general case to these special cases:
\begin{enumerate}
\item In~\ref{Sec:G-tame-etale} we recall the easier case when $G$ is  tame and \'etale.
\item In~\ref{Sec:case-loc-diag}-\ref{Rem:loc-diag-gerbe} we prove the case where $G$ is locally diagonalizable. This involves torsion on the relative Picard stack of line bundles, and is precisely the place where \lci singularities rear their ugly head -- specifically Lemma~\ref{Lem:Pic[X]-lci}.  
\item In~\ref{Sec:fixed-points}-\ref{cohlem} we discuss fixed points of torsion of semiabelian schemes under a tame  action of a finite group. The main result here is \ref{5.9}, which shows that the fixed point set remains \lci. This is used later in  the general case in \ref{1.5.1}.  
\item The general case is set up in \ref{Sec:General-G-setup}. It is quite intricate; the strategy is outlined in \ref{Sec:General-G-strategy}.
\end{enumerate}

\subsection{The case when $G$ is a tame \'etale group scheme}\label{Sec:G-tame-etale}

In this case we claim that \ref{flatmap} in fact is \'etale.  Indeed to verify this it suffices to show that \ref{flatmap} is formally \'etale since it is a morphism locally of finite presentation.  If $T_0\hookrightarrow T$ is an infinitesimal thickening and $\mc C_T\rightarrow T$ is a twisted curve over $T$, then the reduction functor from $G$-torsors 
on $\mc C_T$ to $G$-torsors on $\mc C_{T_0}$ is an equivalence of categories since $G$ is \'etale.

\subsection{The case when $G$ is locally diagonalizable}\label{Sec:case-loc-diag}  Let $T$ be an $S$-scheme, and let $\mc C\rightarrow T$ be a twisted curve. The stack over $T$ of morphisms $\mc C\rightarrow \cB G$ is then equivalent to the stack $\text{\rm TORS}_{\mc C/T}(G)$  associating to any $T'\rightarrow T$ the groupoid of $G$-torsors on $\mc C\times _TT'$.  Let $X$ denote the Cartier dual $\text{Hom}(G, \mathbb{G}_m)$ so that $G = \Sp_S\cO_S[X]$.

We use the theory of Picard stacks -- namely ``commutative group-stacks" --  see \cite[XVIII.1.4]{SGA4}. Let $\scr Pic_{\mc C/T}$ denote the Picard stack of line bundles on $\mc C$, and let $\text{Pic}_{\mc C/T}$ denote the rigidification of $\scr Pic_{\mc C/T}$ with respect to $\mathbb{G}_m$, so that $\text{Pic}_{\mc C/T}$ is the relative Picard functor of $\mc C/T$.  By Lemma \ref{cokernel} we have that $\text{Pic}_{\mc C/T}$ is an extension of a semi-abelian group scheme by an \'etale group scheme.

Let $\text{Pic}_{\mc C/T}[X]$ denote the group scheme of homomorphisms
\begin{equation*}
X\rightarrow \text{Pic}_{\mc C/T}.
\end{equation*}

\begin{lemma}\label{Lem:Pic[X]-lci}
The scheme $\text{\rm Pic}_{\mc C/T}[X]$ is flat over $T$, with \lci fibers.
\end{lemma}
\begin{proof} 
The assertion is clearly local in the \'etale topology on $S$ so we may assume that $S$ is connected and $G$ diagonalizable.  We may write $G = \prod \bmu_{n_i}$ and $X = \prod \ZZ/(n_i)$. 
Then 
$\text{Pic}_{\mc C/T}[X] = \prod_i \text{Pic}_{\mc C/T}[\ZZ/(n_i)]$. 

It then suffices to consider the case  where $X = \ZZ/(n)$, in which case $\text{Pic}_{\mc C/T}[X]=\text{Pic}_{\mc C/T}[n]$. This is the fiber of the map $\text{Pic}_{\mc C/T}\stackrel{\times n}{\longrightarrow}\text{Pic}_{\mc C/T}$ over the identity section. This  is a map of smooth schemes of the same dimension. It has finite fibers since $\text{Pic}_{\mc C/T}$ is an extension of a semi-abelian group scheme by an \'etale group scheme. Therefore this map is flat with \lci fibers.  Hence $\text{Pic}_{\mc C/T}[n]\to T$ is flat with \lci fibers as well.
\end{proof}

Let $\scr Pic_{\mc C/T}[X]$ denote the Picard stack of morphisms of Picard stacks
\begin{equation*}
X\rightarrow \scr Pic_{\mc C/T},
\end{equation*}
where $X$ is viewed as a discrete stack.
Then $\scr Pic_{\mc C/T}[X]$ is a $G$--gerbe over $\text{Pic}_{\mc C/T}[X]$, hence flat over $T$.

The result in the locally diagonalizable case therefore follows from the following lemma:

\begin{lem}\label{Lem:tors-picx} There is an equivalence of categories
\begin{equation*}
\text{\rm TORS}_{\mc C/T}(G)\rightarrow \scr Pic_{\mc C/T}[X],
\end{equation*}
were $\text{\rm TORS}_{\mc C/T}(G)$ denotes the stack associating to any $T'\rightarrow T$ the groupoid of $G$-torsors on $\mc C\times _TT'$.\end{lem}
\begin{proof}
If  $p:P\rightarrow \mc C$ is a $G$-torsor, then $P$ is finite over $\mc C$ and therefore equal to the spectrum of the finite flat $\mc O_{\mc C}$-algebra
$$
\pi _*\mc O_{P}.
$$
The $G$-action on $P$ defines a $G$-action on $\pi _*\mc O_{\mc C}$, and therefore a decomposition
$$
\pi _*\mc O_{P} = \oplus _{x\in X}\mc L_x,
$$
where $\mc L_x\subset \pi _*\mc O_{P}$ is the subsheaf of element on which $G$ acts through the character $x$.  Each of the $\mc O_{\mc C}$-modules is locally free of rank $1$ as this can be verified locally when $P$ is the trivial torsor.  Moreover, the algebra structure on $\pi _*\mc O_P$ defines isomorphisms
$$
\mc L_x\otimes \mc L_{x'}\rightarrow \mc L_{x+x'}
$$
giving a morphism of Picard stacks
\begin{equation*}
F\colon X\rightarrow \scr Pic_{\mc C/T}.
\end{equation*}
Conversely given such a morphism $F$, let $\mc L_x$ denote $F(x)$.  The  isomorphisms
\begin{equation*}
F(x+x')\simeq F(x)+F(x')
\end{equation*}
define an algebra structure on 
\begin{equation*}
\oplus _{x\in X}\mc L_x,
\end{equation*}
and the $X$-grading defines a $G$-action on this sheaf of algebras such that
$$
\underline {\text{Spec}}_{\mc C}(\oplus _{x\in X}\mc L_x)\rightarrow \mc C
$$
is a $G$-torsor.

\end{proof}

\begin{remark} \label{Rem:loc-diag-gerbe}
With a bit more work, one can prove a more general result which may be of interest: given a twisted curve $\cC \to T$ and a $G$-gerbe $\cG \to \cC$, the stack $\Sec_T(\cG / \cC)$ is a $G$-gerbe over its rigidification $\rigSec_T(\cG / \cC)$, and the $T$-space $\rigSec_T(\cG / \cC)$ is a pseudo-torsor under the flat group-scheme  $\text{Pic}_{\mc C/T}[X]$. In particular  $\Sec_T(\cG / \cC)\to T$ is flat.
\end{remark}

\subsection{Observations on fixed points}\label{Sec:fixed-points}
Before we consider general $G$, 
we  make some observations about schemes of fixed points of group actions on semi-abelian group schemes. 

Let $S$ be a scheme, and let $A\rightarrow S$ be a smooth commutative group scheme.  Let $A^{\rmo}$ denote the connected component of the identity and assume we have an exact sequence of group schemes
\begin{equation}\label{Aseq}
0\arr A^{\rmo}\arr A\arr W\arr 0,
\end{equation}
with $W$ a finite \'etale group scheme over $S$.  In what follows we will assume that $W$ is a constant group scheme (this always hold after making an \'etale base change on $S$) and that the group scheme $A^{\rmo}$ is a semiabelian scheme.

Let $H$ be a finite group of order invertible on $S$ acting on $A$ (by homomorphisms of group schemes).  Let $N$ denote the order of $H$.    Since $A^{\rmo}$ is a semiabelian scheme, multiplication by $N$ is surjective and \'etale on $A^{\rmo}$.  It follows that multiplication by $N$ on $A$ is also \'etale and that there is an exact sequence
\begin{equation*}
\begin{CD}
A@>\times N>> A@>>> W/NW@>>> 0.
\end{CD}
\end{equation*}
In particular, the image of $\times N\colon A\rightarrow A$ is an open and closed subgroup scheme $A'\subset A$ preserved by the $H$-action.

\begin{lem}\label{smoothlem} The scheme of fixed points $A^H$ is a smooth group scheme over $S$.
\end{lem}

\begin{proof}
It is clear that $A^{H}$ is a group scheme; we need to show that it is smooth. This follows from the following Lemma (which is well known, but for which we don't have a good reference).
\end{proof}

\begin{lemma}
Let $X$ be a smooth finitely presented scheme over a scheme $S$, and let $G$ be a finite group whose order is invertible in $S$. Suppose that $G$ acts on $X$ via $S$-scheme automorphisms. Then the scheme of fixed points $X^{G}$ is smooth and finitely presented over $S$.\
\end{lemma}

\begin{proof}
The problem is local on $S$, so we may assume that $S = \spec R$ is affine. Then $X$, together with the action, will be obtained by base change from a smooth scheme over a finitely generated subring of $R$. So we may assume that $R$ is noetherian; and then the scheme $X^{G}$, which is by definition a closed subscheme in the $|G|$-th product $X \times\ldots \times X,$  is automatically finitely presented.

To check smoothness we are going to use Grothendieck's infinitesimal lifting criterion. Each fixed point of $X$ has an affine invariant neighborhood; we may replace $X$ by such neighborhood and assume that $X$ is affine; say $X = \spec O$.
Then $X^{G}$ is the spectrum of the largest invariant ring quotient of $O$; a morphism $\spec A \arr X^{G}$ correspond to an invariant homomorphism $O \arr A$.

Suppose that we are given a commutative diagram
   \[
   \xymatrix{
   R \ar[r]^{\rho}\ar[d] &A\ar[d]\\
   O \ar[r]^{\phi}\ar@{.>}[ur]^{\psi} &B
   }
   \]
(without the dotted arrow), where $A \arr B$ is a surjective homomorphism of artinian rings, whose kernel $I$ is such that $I^{2} = 0$, while $\rho$ and $\phi$ are invariant ring homomorphisms. We need to find an invariant ring homomorphism $\psi\colon O \arr A$ that fits above. Since $O$ is smooth over $R$, we can find a ring homomorphism $\psi$, which is not however necessarily invariant. For each $\sigma \in G$ the homomorphism $O \arr I$ defined by the rule $f \mapsto \psi(\sigma f) -- \psi(f)$ is an $R$-linear derivation; hence it defines a function from $G$ to the group $\der_{R}(O,I)$ of $R$-linear derivations from $O$ to $I$. There is a natural right action of $G$ on $\der_{R}(O,I)$, defined by the rule $D\sigma(f) = D(\sigma f)$; the function $G \arr \der_{R}(O,I)$ above is a $1$-cocycle with respect to such action. Since $G$ is tame we have $\H^{1}\bigl(G,\der_{R}(O,I)\bigr) = 0$, so there exists $D \in \der_{R}(O,I)$ such that
   \[
   \psi(\sigma f) -- \psi(f) = D(\sigma f) -- D(f)
   \]
for all $\sigma  \in G$ and all $f \in O$. The function $O \arr A$ defined by $f\mapsto \psi(f) -- D(f)$ is the desired invariant ring homomorphism.
\end{proof}

%\begin{proof}
%We verify the infinitesimal lifting property.  Let $T_0\hookrightarrow T$ be a closed immersion defined by a square-zero ideal $I\subset \mc O_T$, and let $p_0\in A^H(T_0)$ be an $H$-invariant point.  We show that $p_0$ can be lifted to an $H$-invariant point of $A(T)$.  Since $A$ is smooth we can after perhaps shrinking on $T$ lift $p_0$ to some point $p\in A(T)$.  Define
%\begin{equation*}
%\widetilde q:= \sum _{h\in H}p^h\in A(T),
%\end{equation*}
%where $p^h$ denotes the image of $p$ under $h:A\rightarrow A$.  Note that in fact $\widetilde q\in A^H(T)$ and that $\widetilde q$ reduces to $Np_0$ in $A(T_0)$ so that $\widetilde q\in A^{\prime H}(T)$.
%Since $N$ is invertible in $S$ multiplication by $N:A\rightarrow A'$ is \'etale.  Therefore after replacing $T$ by an \'etale extension there exists an element $q\in A(T)$ such that $Nq = \widetilde q$.  In particular, the reduction $q_0\in A(T_0)$  of $q$ is an element with $N(q_0-p_0) = 0$.  Since the group scheme $A[N]$ is \'etale over $T$ we can after changing $q$ by a point of $A[N]$ assume that $q$ reduces to $p_0$.    For any $h\in H$, the point $q^h-q\in A(T)$ is a point annihilated by $N$ since
%\begin{equation*}
%N(q^h-q) = (Nq)^h-(Nq) = \widetilde q^h-\widetilde q = 0.
%\end{equation*}
%Since the reduction map
%\begin{equation*}
%A[N](T)\rightarrow A[N](T_0)
%\end{equation*}
%is injective it follows that $q\in A^H(T)$.
%\end{proof}

Now let $X$ be a finitely generated $\Z/p^n$--module with $H$-action, where $p$ is prime to $N$.  Let $A[X]$ denote the scheme of homomorphisms $X\rightarrow A$.  The group $H$ acts on $A[X]$ as follows.  An element $h\in H$ sends a homomorphism $\rho \colon X\rightarrow A$ to the homomorphism
\begin{equation}\label{1.5.6}
\begin{CD}
X@>h^{-1}>> X@>\rho >> A@>h>> A.
\end{CD}
\end{equation}
By definition of fixed points, the group scheme $(A[X])^H$ is the fiber product of the diagram
$$
\begin{CD}
@.  A[X] \\
@. @VV\Delta V \\
A[X]@>{\prod _{h\in H}(\text{$h$-action})}>> \prod _{h\in H}A[X].
\end{CD}
$$
Note that the group scheme of fixed points $(A[X])^H$ is potentially quite unrelated to $A^H$.

\begin{prop}\label{5.9} The group scheme of fixed points $(A[X])^H$ is flat over $S$, with \lci fibers.
\end{prop}
\begin{proof}
Choose a presentation of $X$ as an $H$-representation
\begin{equation}\label{shortseq}
0\arr K \arr F \arr X \arr 0,
\end{equation}
where the underlying groups of $K$ and $F$ are free $\Z$-modules of finite rank.  Let $A[F]$ (resp. $A[K]$) denote the space of homomorphisms $F\rightarrow A$ (resp. $K\rightarrow A$).  Applying $\text{Hom}(F, -)$ (resp. $\text{Hom}(K, -)$) to \ref{Aseq} we see that $A[F]$ (resp. $A[K]$) sits in a short exact sequence
\begin{equation*}
0 \arr A^{\rmo}[F] \arr A[F] \arr \text{Hom}(F, W) \arr 0
\end{equation*}
\begin{equation*}
(\text{resp. } 0 \arr A^{\rmo}[K] \arr A[K] \arr \text{Hom}(K, W) \arr 0),
\end{equation*}
where $A^{\rmo}[F]$ (resp. $A^{\rmo}[K]$) is a semiabelian  scheme.  Note also that the relative dimension over $S$ of $A[F]$ (resp. $A[K]$) is equal to $\text{dim}(A)\cdot \text{rk}(F)$ (resp. $\text{dim}(A)\cdot \text{rk}(K)$).  Since $X$ is a torsion module we have $\text{rk}(F) = \text{rk}(K)$, so $A[F]$ and $A[K]$ have the same dimension.
 The inclusion $K\hookrightarrow F$ induces a homomorphism
\begin{equation}\label{map1}
A[F] \arr A[K]
\end{equation}
whose kernel is $A[X]$. 

\begin{lem} The induced map
\begin{equation*}
A^{\rmo}[F] \arr A^{\rmo}[K]
\end{equation*}
is surjective.
\end{lem}

\begin{proof}
The short exact sequence \ref{shortseq} induces an exact sequence of fppf-sheaves on $S$
\begin{equation*}
A^{\rmo}[F] \arr A^{\rmo}[K] \arr \scr Ext^1(X, A^{\rmo}).
\end{equation*}
Since $A^{\rmo}$ is divisible, it is an injective objects in the category of abelian groups, so $\scr Ext^1(X, A^{\rmo})  = 0$. The surjectivity of $A^{\rmo}[F] \arr A^{\rmo}[K]$ follows.
\end{proof}

  By \ref{smoothlem} the induced morphism 
\begin{equation}\label{kermap}
f\colon (A[F])^H \arr (A[K])^H
\end{equation}
is a morphism of smooth group schemes with kernel the group scheme $(A[X])^H$. This is finite since $A[X]$ is a finite flat group scheme.  Note that the connected component of the identity in $(A[F])^H$ (resp. $(A[K])^H$) is equal to the connected component of the identity in $A^{\rmo}[F]^H$ (resp. $A^{\rmo}[K]^H$).

\begin{lem} The morphism 
\begin{equation*}
A^{\rmo}[F]^H \arr A^{\rmo}[K]^H
\end{equation*}
induced by \ref{kermap} is surjective.
\end{lem}
\begin{proof}
Let $k$ be a field and suppose given an $H$-invariant homomorphism $f\colon K\rightarrow A^{\rmo}$ over $k$.    We wish to show that after possibly making a field extension of $k$ we can find an extension $\widetilde f\colon F\rightarrow A^{\rmo}$ of $f$ which is $H$-invariant. Replacing $S$ by $\Sp (k)$, for the rest of the proof we view everything as being over $k$.  Consider the short exact sequence of commutative group schemes with $H$-action
\begin{equation*}
0 \arr A^{\rmo}[X] \arr A^{\rmo}[F] \arr A^{\rmo}[K] \arr 0.
\end{equation*}
Viewing this sequence as an exact sequence of abelian sheaves on $\cB H_{\fppf}$ we obtain by pushing forward to $\Sp (k)_{\fppf}$ an exact sequence
\begin{equation*}
A^{\rmo}[F]^H(k) \arr A^{\rmo}[K]^H(k) \arr H^1(\cB H_{\fppf}, A^{\rmo}[X]).
\end{equation*}
The lemma therefore follows from  Lemma \ref{cohlem} below (which is stated in the generality needed later).
\end{proof}

\begin{lem}\label{cohlem} Let $\mc D$ be an abelian sheaf on $\cB H_{T, \fppf}$ such that every local section of $\mc D$ is torsion of order prime to $N$, the order of $H$.  Let 
\begin{equation*}
f\colon \cB H_{T, \fppf} \arr T_{\fppf}
\end{equation*}
be the topos morphism defined by the projection.  Then $R^if_*\mc D = 0$ for all $i>0$.
\end{lem}
\begin{proof} Let $T_\bullet \rightarrow \cB H_T $ be the simplicial scheme over $\cB H_T$
associated to the covering $T\rightarrow \cB H_T$ defined by the trivial torsor
(as in \cite[12.4]{L-MB}).  For $p\geq 0$ let $f_p\colon T_{p, \fppf}\rightarrow T_{\fppf}$ be the projection.  Note that $T_p$ is a finite disjoint union of copies of $T$, and therefore for any abelian sheaf $F$ on $T_{p, \fppf}$ we have $R^if_{p*}F = 0$ for $i>0$.  From this it follows that $Rf_*\mc D$ is quasi-isomorphic in the derived category of abelian sheaves on $T_{\fppf}$ to the complex 
\begin{equation*}
f_{0*}\mc D|_{T_0} \arr f_{1*}\mc D|_{T_1} \arr \cdots .
\end{equation*}
Therefore the sheaf $R^if_*\mc D$ is isomorphic to the sheaf associated to the presheaf which to any $T'\rightarrow T$ associates the group cohomology 
\begin{equation}\label{cohgroups}
H^i(H, \mc D(T')),
\end{equation}
where $\mc D(T')$ denotes the $H$-module obtained by evaluating $\mc D$ on the object $T'\rightarrow \cB H$ corresponding to the trivial torsor on $T'$.  Since $\mc D(T')$ is a direct limit of  groups of  order prime to the order of  $H$, it follows that the groups \ref{cohgroups} are zero for all $i>0$.
\end{proof}

Now we complete the proof of Proposition \ref{5.9}.
It follows that the map on connected components of the identity
\begin{equation*}
\left(A[F]^H\right)^{\rmo} \arr \left(A[K]^H\right)^{\rmo}
\end{equation*}
 is a surjective homomorphism of smooth group schemes of the same dimension, and hence a flat morphism with \lci fibers. 
Therefore the morphism
 \begin{equation*}
 A[F]^H \arr A[K]^H
 \end{equation*}
 is also flat and hence its fiber   over the identity section, namely $A[X]^H$, is flat over $S$, with \lci fibers.
\end{proof}

\subsection{General $G$: setup}\label{Sec:General-G-setup}
We return to the proof of \ref{flattheorembis}.

The assertion that \ref{flatmap} is flat is local in the fppf topology on $S$.  We may therefore assume that $G$ is well-split so that there is a split exact sequence
\begin{equation*}
1\arr D\arr G\arr H\arr 1
\end{equation*}
with $H$ \'etale and tame and $D$ diagonalizable. 
% Write $G = D\rtimes H$. 
  The map $G\rightarrow H$ induces a morphism
over $\mc M_{g, n}^{\text{tw}}\times S$
\begin{equation}\label{flatmap2}
\mc M_{g, n}^{\text{tw}}(G)\arr \mc M_{g, n}^{\text{tw}}(H).
\end{equation}
We may apply the case of \'etale group scheme (Section \ref{Sec:G-tame-etale}) 
to  the group scheme $H$. Therefore we know that $\mc M_{g, n}^{\text{tw}}(H)$ is \'etale over $\mc M_{g, n}^{\text{tw}}\times S$, so it suffices to show that the morphism \ref{flatmap2} is flat with \lci fibers.

For this in turn it suffices (for example by \cite[IV.11.5.1]{EGA}) to show that for any morphism $t\colon T\rightarrow \mc M_{g, n}^{\text{tw}}(H)$, with $T$ the spectrum of an artinian local ring, the induced map
\begin{equation}\label{relflat}
\mc M_{g, n}^{\text{tw}}(G)\times _{\mc M_{g, n}^{\text{tw}}(H)}T\arr T
\end{equation}
is flat with \lci fibers.

Let $(\mc C\rightarrow T, F\colon \mc C\rightarrow \cB H)$ be the $n$-pointed genus $g$ twisted curve together with the morphism $F\colon \mc C\rightarrow \cB H$  corresponding to $t$, and let $\mc C'\rightarrow \mc C$ denote the $H$-torsor obtained by pulling back the tautological $H$-torsor over $\cB H$.  So we have a diagram with cartesian  square
\begin{equation*}
\begin{CD}
\mc C'@>>> \mc C\\
@VVV @VVFV \\
T@>>> \cB H_T\\
@. @VVfV \\
 @. T.
\end{CD}
\end{equation*}

For the rest of the argument we may replace $S$ by $T$, pulling back all the objects over $S$ to $T$.
 
   The action of $H$ on $D$ defines a group scheme $\cD$ over $\cB H_T$ whose pullback along the projection $T\rightarrow \cB H_T$ is the group scheme $D$ with descent data defined by the action. Let $\mathbb{D}$ denote the pullback of this group scheme to $\mc C$.  The pullback of $\mathbb{D}$ to $\mc C'$ is equal to $D\times _T\mc C'$, but $\mathbb{D}$ is a twisted form of $D$ over $\mc C$. Similarly, the character group $X$ of $D$ has an action of the group $H$, giving a group scheme  $\scr X$ over $\cB H_T$ which is Cartier dual to $\cD$. This is summarized in the following Cartesian diagram:
   
   $$\xymatrix{
   D\times_T\cC' \simeq \cC'\times_\cC\mathbb{D}\ar[rr]\ar[d] & &  \mathbb{D}:=\cD \times_{\cB H_T}\cC\ar[d] \\
   D\simeq \Hom_T(X,\GG_m)\ar[rr]&& \cD=\Hom_{\cB H_T}(\scr X, \GG_m).
   }$$
   
   Let $\mc G\rightarrow \mc C$ denote the stack $\mc C\times _{\cB H}\cB G$, and let $\mc G'\rightarrow \mc C'$ denote the pullback to $\cC'$. The stack $\mc G$ is a gerbe over $\mc C$ banded by the sheaf of groups $\mathbb{D}$. Following a standard observation as in Appendix \ref{AppendixB}, we identify the stack on the left in  \ref{relflat} with the stack $\Sec_T (\mc G/\mc C)$, whose objects over a $T$-scheme $T'$ are sections $\cC_{T'} \to \cG_{T'}$ of   $\cG_{T'} \to \cC_{T'}$. Therefore 
to prove that the morphism \ref{relflat} is flat with \lci fibers we need to show that the stack $\Sec_T (\mc G/\mc C)$ is flat with \lci fibers over $T$. As $T$ is local artinian, we may assume that $\Sec_T (\mc G/\mc C) \to T$ is set theoretically surjective, and replacing $T$ by a finite flat cover we may also assume that $\Sec_T (\mc G/\mc C) \to T$ has a section over the reduction of $T$. 

Note that the gerbe $\mc G'\rightarrow \mc C'$ is trivial: there is an isomorphism $\mc G' \simeq \cC' \times BD$. This is because $T \to \cB H$ corresponds to the trivial torsor, and any trivialization gives an isomorphism $T\times_{\cB H} \cB G \simeq BD$. It follows that there is an isomorphism of $T$-schemes $\Sec_T (\mc G'/\mc C') \simeq \text{\rm TORS}_{\mc C'/T}(D)$. However this isomorphism is not canonical.

\subsection{Reduction to $\mc C'$ connected.}\label{SS:5.17} We wish to apply previous results such as \ref{Lem:tors-picx} to $\cC' \to T$ -- see \ref{Sec:structure-of-Sec} below -- 
but our discussion applied only when this is connected. We reduce to the connected case as follows: suppose the $H$-bundle $\cC' \to \cC$ has a disconnected fiber over the residue field of $T$, then since T is
artinian we can choose a connected component $\cC''\subset \cC'$. Let $H''\subset H$ be the subgroup sending $\cC''$ to itself. Then $\cC'' \arr \cC$ is an
 $H''$-bundle. If we denote  $G''=G \times_HH''$, there is an equivalence of categories between $G$-bundles $\cE\to \cC$ lifting the given $H$-bundle  $\cC'$ and $G''$-bundles $\cE'' \to \cC$
lifting the $H''$ bundle $\cC''$: one direction is by restricting $\cE$ to $\cC''\subset \cC'$, the other direction by inducing $\cE'' \to \cC$ from $G''$ to $G$. It therefore suffices to consider the case when $\cC'\to T$ has connected fibers.

\subsection{General $G$: strategy}\label{Sec:General-G-strategy} Our approach goes as follows:
\begin{enumerate}
\item By Lemma \ref{Lem:tors-picx} we have a precise structure, as a gerbe over a group scheme, of the stack $\Sec_T (\mc G'/\mc C')$.
\item 
%The action of $H$ on $\Sec_T (\mc G'/\mc C')$ is compatible with this structure, in a torsorial sense to be discussed below. 
This provides an analogous precise structure on $\Sec_{\cB H_T} (\mc G/\mc C)$ and its rigidification $\rigSec_{\cB H_T} (\mc G/\mc C)$: we have that $\Sec_{\cB H_T} (\mc G/\mc C)$ is a $\cD$-gerbe over $\rigSec_{\cB H_T} (\mc G/\mc C)$. However $ \rigSec_{\cB H_T} (\mc G/\mc C)$ is not a group-scheme but a torsor (see below).
\item The stack $$\Sec_T (\mc G/\mc C) = \Sec_T\left(\Sec_{\cB H_T} (\mc G/\mc C)\,/\, \cB H_T\right)$$ 
can be thought of as a stack theoretic version of push-forward in the fppf topology from $\cB H_T$ to $T$. Indeed, for a representable morphism $U \to \cB H_T$, the fppf sheaf associated to the space $\Sec_T(U/\cB H_T)$ coincides with $f_*(U)_{\fppf}$. We analyze the structure of the push-forward of its building blocks, namely the rigidification $\text{\rm Sec}_{\cB H_T}(\mc G/\mc C)$ and the group-schemes underlying the torsor and gerbe structures.
\end{enumerate}

\subsection{Structure of $\Sec_T (\mc G'/\mc C')$ and $\Sec_{\cB H_T} (\mc G/\mc C)$}\label{Sec:structure-of-Sec}
Let us explicitly state the structure in (1) above:

\begin{itemize}
\item $\Sec_T (\mc G'/\mc C') \to \rigSec_T (\mc G'/\mc C')$ is a gerbe banded by the group-scheme $D$, and
\item the rigidification $\rigSec_T (\mc G'/\mc C')$ of $\Sec_T (\mc G'/\mc C')$ is isomorphic to the $T$-group scheme $\text{Pic}_{\mc C'/T}[X]$ (this follows from Lemma \ref{Lem:tors-picx} and the paragraph preceding \ref{SS:5.17}).
\end{itemize}

We  proceed with the structure in (2).
\begin{lemma}\label{Lem:sec-BH-gerbe}
The automorphism group of any object of $\Sec_{\cB H_T} (\mc G/\mc C)$
  over $Z \to \cB H$ 
  is canonically isomorphic to $\cD(Z)$
\end{lemma}

\begin{proof}
This follows from the facts that $\mc G$ is a $\cD$-gerbe and $\cC \to
  \cB H$ connected.
  Indeed if $s\colon \mc C\rightarrow \mc G$ is a section, then an
  automorphism of $s$ is given by a map $\mc C\rightarrow \cD$. This
  necessarily factors through $\cB H_T$ since by connectedness $F_*\mc
  O_{\mc C} = \mc O_{\cB H_T}$ (where $F\colon \mc C\rightarrow \cB H_T$ is the
  structure morphism) and $\cD$ is affine over $\cB H_T$. 
\end{proof}

The stack $\Sec_{\cB H_T} (\mc G/\mc C)$ is algebraic by Lemma \ref{Lem:relative-algebraic}.
We can define $\rigSec_{\cB H_T} (\mc G/\mc C)$ to be the rigidification $\Sec_{\cB H_T} (\mc G/\mc C)\thickslash \cD$.
This is also the relative coarse moduli space since $\rigSec_{\cB H_T}
(\mc G/\mc C)\to \cB H$ is representable. 
We have  that $\Sec_{\cB H_T} (\mc G/\mc C)\to \rigSec_{\cB H_T} (\mc G/\mc
C)$ is a gerbe banded by $\cD$. 

We also have naturally that $\Sec_T (\mc G'/\mc C') =
T\times_{\cB H_T}\Sec_{\cB H_T} (\mc G/\mc C)$. The lemma above and the base
change property of rigidification gives
   \[
   \rigSec_T (\mc G'/\mc C') = T\times_{\cB H_T}\rigSec_{\cB H_T} (\mc G/\mc C).
   \]
Let us look at the structure underlying $\rigSec_{\cB H_T} (\mc G/\mc C)$.

 The group $H$ acts on $\text{Pic}_{\mc C'/T}$ and on $X$. We therefore also obtain  a left action of $H$ on
\begin{equation*}
\text{Pic}_{\mc C'/T}[X] = \text{Hom}(X, \text{Pic}_{\mc C'/T}),
\end{equation*}
where $h\in H$ sends a homomorphism $\rho \colon X\rightarrow \text{Pic}_{\mc C'/T}$ to the homomorphism
\begin{equation}
\begin{CD}
X@>h^{-1}>> X@>\rho >> \text{Pic}_{\mc C'/T}@>h>> \text{Pic}_{\mc C'/T}.
\end{CD}
\end{equation}

This defines an fppf sheaf on $\cB H_T$, which is easily seen to be represented by $\text{Pic}_{\mc C/\cB H_T}[\scr X]$. 

Note  that in general the identification  $\rigSec_T (\mc G'/\mc C') = T\times_{\cB H_T}\rigSec_{\cB H_T} (\mc G/\mc C)$ does not give an isomorphism of $\rigSec_{\cB H_T} (\mc G/\mc C)$ with $\text{Pic}_{\mc C/\cB H_T}[\scr X]$, because the descent data may  in general be different. However, we have the following:

\begin{lemma}
There is  a natural action of $\text{\rm Pic}_{\mc C/\cB H_T}[\scr X]$ on $\rigSec_{\cB H_T} (\mc G/\mc C)$, which pulls back to the natural action of $\text{\rm Pic}_{\mc C'/T}[X]$ on itself by translation. In particular  $\rigSec_{\cB H_T} (\mc G/\mc C)$ is a torsor under $\text{\rm Pic}_{\mc C/\cB H_T}[\scr X]$.
\end{lemma}

\begin{proof}
Let $\Phi\colon X\rightarrow \scr Pic_{\mc C/\cB H_T}$ be an object of $\scr Pic_{\mc C/\cB H_T}[\scr X]$ over some arrow $\xi\colon  Z \to \cB H_T$ and $s\colon \mc C\rightarrow \mc G$  a section, again over $Z$. Then we can define a new section $s^\Phi\colon \mc C\rightarrow \mc G$ as follows.  Let $P^\Phi$ denote the $\cD$-torsor corresponding to $\Phi$ (see Lemma \ref{Lem:tors-picx}).  The section $s$ defines an isomorphism (which we denote by the same letter)
\begin{equation*}
s\colon {\mc C}\times \cB\cD\rightarrow \mc G,
\end{equation*}
and we let $s^\Phi$ denote the composite
\begin{equation*}
\begin{CD}
\mc C@>1\times P^\Phi>> \mc C\times \cB\cD@>s>> \mc G.
\end{CD}
\end{equation*}
By the universal property of rigidification this defines an action of the group scheme $\text{Pic}_{\mc C/\cB H_T}[\scr X]$ on $\rigSec_{\cB H_T} (\mc G/\mc C)$. The fact that its pullback is identified with the action of $\text{Pic}_{\mc C'/T}[X]$ on itself by translation is routine.
\end{proof}

\subsection{The pushforward  of  $\text{Pic}_{\mc C/\cB H_T}[\scr X]$ and its torsor $\rigSec_{\cB H_T} (\mc G/\mc C)$}

We write $(\text{Pic}_{\mc C/\cB H_T}[\scr X])_{\fppf}$ for the sheaf in the fppf topology on $\cB H_T$ represented by $\text{Pic}_{\mc C/\cB H_T}[\scr X]$.
The  fppf pushforward $f_*((\text{\rm Pic}_{\mc C/\cB H_T}[\scr X])_{\fppf})$ under the morphism $f:\cB H_T \to T$ is  represented by the scheme of fixed points $\text{\rm Pic}_{\mc C'/T}[X]^H$.  Indeed if $T'$ is a $T$-scheme, then a section of $f_*((\text{\rm Pic}_{\mc C/\cB H_T}[\scr X])_{\fppf})$ over $T'$ is given by a morphism of sheaves on $BH_{T'}$
$$
\scr X\rightarrow \text{Pic}_{\mc C_{T'}/BH_{T'}}.
$$
Such a morphism  corresponds by descent theory to a morphism
$$
\rho :X\rightarrow \text{Pic}_{\mc C'_{T'}/T'}
$$
such that for any $T'$-scheme $T^{\prime \prime }$ and $h\in H_T(T^{\prime \prime })$ the diagram 
$$
\xymatrix{
X\ar[r]^-h\ar[d]^-\rho & X\ar[d]^-\rho \\
\text{Pic}_{\mc C'_{T^{\prime \prime }}/T^{\prime \prime }}\ar[r]^-h& \text{Pic}_{\mc C'_{T^{\prime \prime }}/T^{\prime \prime }}}
$$
commutes.  Equivalently, the point $\rho \in \text{Pic}_{\mc C'/T}[X](T')$ is $H$-invariant.

We have:

\begin{lem}\label{1.5.1} The scheme of fixed points $\text{\rm Pic}_{\mc C'/T}[X]^H$ is finite and  flat over $T$, with \lci fibers.
\end{lem}
\begin{proof}
Note that in the notation of Section \ref{3.7}, since $X$ is a torsion group we have
$$
\text{Pic}_{\mc C'/T}[X] = \text{Pic}^{\rmo} _{\mc C'/T}[X]
$$
where $\text{Pic}^{\rmo} _{\mc C'/T}$ is the group scheme of degree-0 line bundles defined in  \ref{Def:Pic^o}. Now  by \ref{3.7}, we have that $\text{Pic}^{\rmo}_{\mc C'/T}$ is an extension of a semiabelian  scheme by an \'etale group scheme. The latter \'etale group scheme is finite since $T$ is assumed Artinian local.  The lemma therefore follows from Proposition \ref{5.9}.
\end{proof}

Now consider the torsor  $\rigSec_{\cB H_T} (\mc G/\mc C)$ and the pushforward under $f$ of the corresponding sheaf $(\rigSec_{\cB H_T} (\mc G/\mc C))_{\fppf}$. We have the following.

\begin{lem}\label{torsorlem} Let $\mc E$ be an abelian sheaf on $\cB H_{T, \fppf}$ such that every local section is torsion of order prime to $N$, the order of $H$, and let $P\rightarrow \cB H_T$ be an $\mc E$--torsor.  Then the sheaf of sets $f_*P$ on $T_{\fppf}$ is a torsor under the sheaf $f_*\mc E$.
\end{lem}
\begin{proof} The fact that $P$ is a torsor under $\mc E$ immediately implies that $f_*P$ is a pseudo-torsor under $f_*\mc E$, namely 
$$\xymatrix{f_*\mc E \times_{\cB H_T} f_*P \ar[rr]^{(\text{action},\text{pr}_2)}&& f_*P \times_{\cB H_T} f_*P}$$
 is an isomorphism. So the only issue is to show that fppf locally on $T$ the sheaf $f_*P$ has a section.  Equivalently, we need to show that after making a flat surjective base change $T'\rightarrow T$ the torsor $P$ itself becomes trivial.  The class of the torsor is a class in $H^1(\cB H_T, \mc E)$. By \ref{cohlem} and the Leray spectral sequence this lies in $H^1(T, f_*\mc E)$. Any class in this group can be killed by making a flat surjective base change $T'\rightarrow T$.
\end{proof}

In our situation, with $\mc E=(\text{Pic}_{\mc C/\cB H_T}[\scr X])_{\fppf}$ and $P=\rigSec_{\cB H_T} (\mc G/\mc C)_{\fppf}$, we obtain that the sheaf $f_* (\rigSec_{\cB H_T} (\mc G/\mc C))_{\fppf}$ is represented by a torsor under $\text{\rm Pic}_{\mc C'/T}[X]^H$. In particular the space representing this sheaf is flat over $T$ with \lci fibers.

We denote this $T$-space by the shorthand notation $\rigSec^H$  -- a complete notation would look like $\Sec_T(\rigSec_{\cB H_T} (\mc G/\mc C) / \cB H_T)$.  

We turn our view to $\Sec_T(\mc G/\mc C)$.

\begin{lemma} 
The automorphism group-scheme of an object of  $\Sec_T(\mc G/\mc C)$ over a
$T$-scheme $B$ is  the group scheme $D^H$ representing $f_*\cD$.
\end{lemma}

\begin{proof}
Let $s\colon \mc C\rightarrow \mc G$ be a section over some $T$-scheme
$B$. Since  $\Sec_T(\mc G/\mc C) = \Sec_T( \Sec_{\cB H}(\mc G/\mc C)\,
/\, \cB H)$, an automorphism of $s$ is a section over $B \times \cB H$ of
the automorphism group-scheme of  $s$ viewed as an object of
$\Sec_{\cB H}(\mc G/\mc C)$ over $B\times \cB H$. By Lemma
\ref{Lem:sec-BH-gerbe} this is a section over $B\times \cB H$ of the
group scheme $\cD$, namely a section of $f_*\cD$, as required.    
\end{proof}

We can define $\rigSec_T(\mc G/\mc C)$ to be the rigidification. 
 We have that the morphism $\Sec_T(\mc G/\mc C)\to \rigSec_T(\mc G/\mc C)$ is a gerbe banded by the group scheme $D^H$ representing $f_*\cD$. This is automatically flat. It is smooth by the following:
\begin{lemma}
 Let $\Delta$ be a diagonalizable group scheme over an algebraically closed field. Then $\cB \Delta$ is smooth. 
\end{lemma}

\begin{proof}
Factoring $\Delta = \prod \bmu_{r_i}$ we have $\cB \Delta=\prod \cB \bmu_{r_i}$, and it suffices to consider each factor. But $\cB \bmu_{r_i} = [\GG_m/\GG_m]$, with the action through the $r_i$-th power map, in other words we have a smooth morphism $\GG_m \to \cB \bmu_{r_i}$ with smooth source, hence $\cB \bmu_{r_i}$ is smooth. 
\end{proof}

It therefore suffices to show that $\rigSec_T(\mc G/\mc C)$ is
flat with \lci fibers. There is a natural map 
$\Sec_T(\mc G/\mc C)\rightarrow \text{Sec}^H$ inducing  $\rigSec_T(\mc
G/\mc C)\rightarrow \text{Sec}^H$. 
The following clearly suffices:

\begin{proposition}\label{Prop:sec-T-gerbe}
The morphism $\Sec_T(\mc G/\mc C)\rightarrow \text{Sec}^H$ is a gerbe
banded by  the group scheme $D^H$ representing $f_*\cD$, and hence
$\rigSec_T(\mc G/\mc C)\rightarrow \text{Sec}^H$ is an isomorphism. 
\end{proposition}

Consider  a morphism $Z\rightarrow \text{Sec}^H$ and the fiber product
\begin{equation*}
\scr S\colon = Z\times _{\text{Sec}^H}\underline {\text{Sec}}(\mc G/\mc C).
\end{equation*}

\begin{lem} Locally in the fppf topology on $Z$ there exists a section $Z\rightarrow \scr S$.
\end{lem}
\begin{proof} Let $\bar s\in f_*(\rigSec_{\cB H_T} (\mc G/\mc C)\,(Z)$ denote the section defined by $Z\rightarrow \text{Sec}^H$. By definition it corresponds to a morphism $\cB H_Z \to \rigSec_{\cB H_T} (\mc G/\mc C)$ over the natural morphism $\cB H_Z \to \cB H_T$. Since $\Sec_{\cB H_T} (\mc G/\mc C)$ is a $\cD$-gerbe over $\rigSec_{\cB H_T} (\mc G/\mc C)$, the obstruction to finding a section of $\scr S$ over $\bar s$ is equal to the class of the $\cD$-gerbe of liftings of $\bar s:\cB H_Z \to \rigSec_{\cB H_T} (\mc G/\mc C)$ to a morphism $\cB H_Z \to \Sec_{\cB H_T} (\mc G/\mc C)$. This class lies in $H^2_{\fppf}(\cB H_Z,\cD)$.  Since by Lemma \ref{cohlem} we have $R^if_*\cD = 0$ for $i>0$, the spectral sequence puts this class in  $H^2_{\fppf}(Z,f_*\cD)$. This vanishes when pulled back to an fppf cover, and the lemma follows. 
\end{proof}

If $s, s'\in \scr S(Z)$ are two sections, then we claim that after making a flat base change on $Z$ the sections $s$ and $s'$ are isomorphic.  Indeed let $I$ be the sheaf over $\cB H_{Z, \fppf}$ of isomorphisms between the two morphisms $\cB H_Z \to \Sec_{\cB H_T} (\mc G/\mc C)$ corresponding to $s$ and $s'$.  Then $I$ is a $\mc D$-torsor over $\cB H_{Z, \fppf}$, and since $R^1f_*\mc D = 0$ it follows that fppf locally on $Z$ this torsor is trivial.  It follows that $\scr S$ is a gerbe banded by $f_*\mc D$. This completes the proof of \ref{Prop:sec-T-gerbe}, implying \ref{flattheorembis}. \qed

\section{Example: reduction of $\scr X(2)$ in characteristic 2}\label{examplesection}

\subsection{$\scr X(2)$ as a distinguished component in $\cK_{0,4}(\cB \bmu_2)$}

Let $\scr X(2)$ denote the stack over $\QQ$ associating to any scheme $T$ the groupoid of pairs $(E, \iota )$, where $E/T$ is a generalized elliptic curve with a full level $2$-structure $\iota :(\Z/2)^2\simeq E[2]$ in the sense of Deligne--Rapoport \cite{DR}. The theory of generalized elliptic curves is rather subtle in general, and we won't review it here, but let us at least recall the following basic facts about such pairs $(E, \iota )$:
\begin{enumerate}
\item $E\rightarrow T$ is a proper flat morphism such that for any geometric point $\bar t\rightarrow T$ the fiber $E_{\bar t}$ is either a connected smooth genus $1$ curve, or the 2-gon (a nodal curve obtained by gluing two copies of $\mathbb{P}^1$ along $0$ and $\infty $).
\item The smooth locus $E^\circ \subset E$ of $E\rightarrow T$ has the structure of a commutative group scheme over $T$ (so in particular there is a section).
\item The translation action of $E^\circ $ on itself extends to an action of $E^\circ $ on $E/T$.
\item The two-torsion subgroup scheme $E[2]\subset E^\circ $ is a finite flat group scheme over $T$ of rank $4$.
\item The involution $P\mapsto -P$ on $E^\circ $ extends to an involution of $E$.
\end{enumerate}

Order the elements of $(\Z/2)^2$ in some way.
For any  pair $(E, \iota )\in \scr X(2)(T)$, the involution $P\mapsto -P$ has fixed points the $2$-torsion points $E[2]\subset E$, and hence the stack-theoretic quotient $[E/\pm 1]$ of $E$ by this involution comes equipped with an ordered set of gerbes $\Sigma _i\subset [E/\pm 1]$, the images of the points $E[2]$,  in the smooth locus.  Furthermore, one checks by direct calculation on the geometric fibers of $E$ that the coarse space of $[E/\pm 1]$ with the resulting four marked points is a stable genus $0$ curve with four marked points.  In this way we obtain a morphism of stacks
$$
\scr X(2)\rightarrow \cK_{0,4}(\cB \bmu _2)_\QQ
$$
sending
$$
(E, \iota )\mapsto (E\rightarrow [E/\pm 1], \{\Sigma _i\subset [E/\pm 1]\}).
$$
One verifies immediately that this functor is fully faithful, and identifies $\scr X(2)$ with the closed substack  $\cK_\QQ\subset \cK_{0,4}(\cB \bmu _2)_\QQ$ classifying data $(B\rightarrow \scr P, \{\Sigma _i\})$ where $B\rightarrow \scr P$ is a $\bmu _2$-torsor such that the resulting map $B\rightarrow P$ to the coarse space of $\scr P$ is ramified over each of the marked points of $P$.

We have seen that  $ \cK_{0,4}(\cB \bmu_2)$ is flat over $\Sp (\ZZ)$. We want to have a nice description of the closure $\cK$ of $\cK_\QQ$, preferably as a naturally defined moduli stack flat over $\ZZ$. It turns out that this can be done in the best possible way: $\cK$ is an open and closed substack of $\cK_{0,4}(\cB \bmu_2)$, defined naturally in terms of the behavior of the $\bmu_2$-cover at the marked points. One can describe this directly, but we find that this gives us a good pretext for introducing the rigidified cyclotomic inertia stack and evaluation maps.

\subsection{Cyclotomic inertia}
Recall that in \cite[\S 3]{AGV} one defined the {\em cyclotomic inertia stacks}: for a fixed positive integer $r\geq 1$ we denote by $\cI_{\bmu_r}(\cX) \to \cX$  the stack whose objects are pairs $(\xi,\phi)$ where $\xi$ is an object of $\cX$ and $\phi: \bmu_r \to \Aut(\xi)$ is a monomorphism, and whose arrows are commutative diagrams as usual; and 
$$\cI_{\bmu}(\cX) = \coprod_{r=1}^\infty \cI_{\bmu_r}(\cX).$$

For a noetherian stack with finite diagonal, The argument of  \cite[Proposition 3.1.2]{AGV} shows that the morphism $\cI_{\bmu}(\cX) \to \cX$ is representable by quasi-projective schemes. Indeed, the fiber over an object $\xi\in \cX(T)$ is the scheme $\coprod_r Hom_{gr-sch/T}(\bmu_r,\Aut_T(\xi))$, which is quasi-projective as the relevant $r$ is bounded. 

\begin{proposition}
If $\cX$ is tame, the morphism $\cI_{\bmu}(\cX) \to \cX$ is finite and finitely presented. \end{proposition}
\begin{proof}
As the problem is local in the fppf topology of the coarse moduli space $X$ of $\cX$, we may assume $\cX = [U/G]$ with $G$ a \lr group-scheme \cite[Theorem 3.2 (c)]{AOV}. The pull-back $\cI_\cX \times_\cX U$ is a closed subgroup scheme of the finite flat \lr group scheme $G_U$. It therefore suffices to show the following Lemma.

\end{proof}
\begin{lemma}\label{grouphom} Let $G_1$ and $G_2$ be two finite flat finitely presented \lr group schemes over a scheme $U$. Then the scheme $\Hom_{\mathrm{gs}/U}(G_1,G_2)$ of homomorphisms of group schemes is finite, flat and finitely presented over $U$.
\end{lemma}

\begin{proof}
This is a local statement in the fppf topology; hence we may assume that $G_{1}$ and $G_{2}$ are of the form $H_{i} \ltimes \Delta_{i}$, where the $\Delta_{i}$ are diagonalizable group schemes whose orders are  powers of the same prime $p$, and $H_{i}$ are constant tame group stacks of orders not divisible by $p$. Let us also assume that $U$ is connected. If $V$ is a scheme over $U$, every homomorphism $G_{1, V} \arr G_{2,V}$ sends $\Delta_{1, V}$ into $\Delta_{2, V}$; hence we have an induced morphism of $U$-schemes
   \[
   \Hom_{\mathrm{gs}/U}(G_1,G_2) \arr 
   \Hom_{\mathrm{gs}/U}(\Delta_{1}, \Delta_{2}) \times_{U}
   \Hom_{\mathrm{gs}/U}(H_1,H_2).
   \]
Since $\Hom_{\mathrm{gs}/U}(\Delta_{1}, \Delta_{2})$ and $\Hom_{\mathrm{gs}/U}(H_1,H_2)$ are both unions of a finite number of copies of $U$, we may fix two homomorphisms of group schemes $\chi\colon \Delta_{1} \arr \Delta_{2}$ and $g\colon H_{1}\arr H_{2}$, and show that the open and closed subscheme $\Hom_{g,\chi}(G_{1}, G_{2})$ of $\Hom_{\mathrm{gs}/U}(G_1,G_2)$ inducing the homomorphisms $g$ and $\chi$ is finite, flat and finitely presented.

If the homomorphism $\chi\colon \Delta_{1} \arr \Delta_{2}$ is not $H_{1}$-equivariant, where $H_{1}$ acts on $\Delta_{2}$ through the homomorphism $f\colon H_{1} \arr H_{2}$, then $\Hom_{g,\chi}(G_{1}, G_{2})$ is empty; so we may assume that $\chi$ is $H_{1}$-equivariant. Then $\Hom_{g,\chi}(G_{1}, G_{2})$ is not empty: in fact, $g \times \chi \colon H_{1} \ltimes \Delta_{1} \arr H_{2} \ltimes \Delta_{2}$ is in $\Hom_{g,\chi}(G_{1}, G_{2})(U)$. If $V$ is a $U$-scheme, an element $f$ of $\Hom_{g,\chi}(G_{1}, G_{2})$ is determined by a morphism of schemes $\gamma\colon H_{1} \arr \Delta_{2}$ so that $f(ax) = g(a)\phi(a)\chi(x)$. By a straightforward calculation, the condition that $f$ be a homomorphism translates into the condition
   \[
   \phi(ab) = \phi(a)^{b}b,
   \]
where we have denoted by $(x, b) \mapsto x^{b}$ the action of $H_{1}$ on $\Delta_{2}(V)$. That is, $\phi$ should be a $1$-cocycle. Since the order of $H_{1}$ is prime to the exponent of $\Delta_{2}(V)$, every such cocyle is coboundary, that is, it is of the form $\phi(a) = d^{a}d^{-1}$ for some $d \in \Delta_{2}(V)$. The coboundary $\phi$ determines $d$ up to an element of $\Delta_{2}(V)^{H_{1}}$; hence we obtain a bijective correspondence between $\Delta_{2}(V)/\Delta_{2}(V)^{H_{1}}$ and $\Hom_{g,\chi}(G_{1}, G_{2})(V)$. Now, it is easy to see that the invariant subsheaf $\Delta_{2}^{H_{1}}$ is in fact a diagonalizable group scheme, and so is the quotient $\Delta_{2}/\Delta_{2}^{H_{1}}$, and we have just produced an isomorphism of functors of $\Delta_{2}/\Delta_{2}^{H_{1}}$ with $\Hom_{g,\chi}(G_{1}, G_{2})$. Since $\Delta_{2}/\Delta_{2}^{H_{1}}$ is finite, flat and finitely presented, this finishes the proof.
\end{proof}

\subsection{Rigidified cyclotomic inertia}

Each object of $\cI_{\bmu_r}(\cX) $ has $\bmu_r$ sitting in the center of its automorphisms. We can thus rigidify it and obtain a stack $ \overline{\cI_{\bmu_r}}(\cX)$ with a canonical morphism
$\cI_{\bmu_r}(\cX)\to \overline{\cI_{\bmu_r}}(\cX)$. Taking the disjoint union over $r\geq 1$ we obtain $\cI_{\bmu}(\cX)\to \overline{\cI_{\bmu}}(\cX)$. The stack $\overline{\cI_{\bmu}}(\cX)$ is called the {\em rigidified cyclotomic inertia stack} of $\cX$. 

The stacks $\overline{\cI_{\bmu_r}}(\cX)$ and $\overline{\cI_{\bmu}}(\cX)$ have another interpretation, discussed in \cite[Section 3.3]{AGV}. Consider the 2-category whose objects consist of morphisms of stacks $\alpha:\cG \to \cX$, where $\cG$ is a gerbe banded by $\bmu_r$ and  $\alpha$ is a {\em representable} morphism; morphisms and 2-arrows are defined in the usual way. Since $\alpha$ is representable,  this 2-category is equivalent to a category. This is isomorphic to the stack  $\overline{\cI_{\bmu_r}}(\cX)$. This gives the stack $\overline{\cI_{\bmu}}(\cX)$ the interpretation as the {stack of cyclotomic gerbes in $\cX$,}  and $\cI_{\bmu}(\cX)\to \overline{\cI_{\bmu}}(\cX)$ together with $\cI_{\bmu}(\cX)\to \cX$ is the universal gerbe in $\cX$.
The non-identity components of $\overline{\cI_{\bmu}}(\cX)$ are known as {\em twisted sectors}.

\subsection{Evaluation maps}

Consider now the stack of twisted stable maps $\cK_{g,n}(\cX,\beta)$. This comes with a diagram
$$
\xymatrix{
\Sigma_i \ar[r]\ar[dr] & \cC_{g,n}(\cX,\beta)\ar[r]\ar[d] & \cX \\
& \cK_{g,n}(\cX,\beta).
}
$$ 
where $ \cC_{g,n}(\cX,\beta)\to \cX$ is the universal twisted stable map and $\Sigma_i$ are the $n$ markings. For each $i$, the resulting diagram
$$
\xymatrix{
\Sigma_i \ar[d]\ar[r] & \cX \\
 \cK_{g,n}(\cX,\beta).
}
$$ 
is a cyclotomic gerbe in $\cX$ parametrized by $\cK_{g,n}(\cX,\beta)$. This gives rise to $n$ {\em evaluation maps}
$$e_i : \cK_{g,n}(\cX,\beta) \longrightarrow \overline{\cI_{\bmu}}(\cX).$$

Evaluation maps (and their sisters, twisted evaluation maps) are of central importance in Gromov--Witten theory of stacks: in the natural gluing maps
$$\xymatrix{\cK_{g_1,n_1+1}(\cX,\beta_1) \times_{\overline{\cI_{\bmu}}(\cX)}\cK_{g_2,1+n_2}(\cX,\beta_2) \ar[r]&\cK_{g_1+g_2,n_1+n_2}(\cX,\beta_2)}$$
the map $\cK_{g_2,1+n_2}(\cX,\beta_2) \to{\overline{\cI_{\bmu}}(\cX)}$ underlying the fiber product on the right is the evaluation map; the map $\cK_{g_1,n_1+1}(\cX,\beta_1) \to{\overline{\cI_{\bmu}}(\cX)}$ on the left is the {\em twisted evaluation map}, namely the evaluation map composed with inversion on the band.

\subsection{Back to $\cK \subset  \cK_{0,4}(\cB \bmu_2)$}
The rigidified cyclotomic inertia stack of $\cB \bmu_2$ has two components -- the identity component is $\cB \bmu_2$ itself, and the non-identity component -- the twisted sector --  is a copy of the base scheme $\spec \ZZ$. 

Since $\cK_\QQ$ corresponds to totally ramified maps,  it is precisely the locus in $\cK_{0,4}(\cB \bmu_2)$ where all four evaluation maps land in the twisted sector. Since ${\cI_{\bmu}}(\cB \bmu_2)$ is finite unramified, the inverse image of each component is open and closed. Therefore $\cK$, the closure of $\cK_\QQ$ in  $\cK_{0,4}(\cB \bmu_2)$ is open and closed. In particular it is flat over $\ZZ$. Also, since the generic fiber is irreducible, $\cK$ is irreducible.  It is simply the gerbe banded by $\bmu_2$ over the $\lambda$ line $\overline\cM_{0,4}\simeq \PP^1$, associated to the class in $H^2_{\fppf}( \PP^1, \bmu_2)$ associated to $c_1(\cO_{\PP^1}(1))$.

\subsection{What does $\cK_{\FF_2}$ parametrize?}

There is one little problem with the closure  $\cK$ of $X(2)_\QQ$: in characteristic 2 it has nothing to do with elliptic curves. 

Given a smooth rational curve with 4 marked points, say at $0,1,\infty$ and $\lambda$, there is a unique $\bmu_2$-bundle over the twisted curve which is degenerate over the coarse curve at all the markings, namely the scheme given in affine coordinates by the equation
$$y^2 = x(x-1)(x-\lambda).$$
In characteristic 2, this is a cuspidal curve of geometric genus 0, with its cusp at the point where $x^2 = \lambda$. This well known example  is possibly the simplest example of geometric interest  of a reduced principal bundle over a smooth scheme with singular total space.

It is not difficult to see that a similar situation holds at the node of a singular curve of genus 0: the singularity of the  $\bmu_2$-bundle is a tacnode.

\subsection{$\scr X(2)$ as a distinguished component in  $\cK_{1,1}(\cB \bmu_2^2)$.}\label{Sec:ex11}

Let $$\cK_{1, 1}^\circ (\cB \bmu _2^2)\subset  \cK_{1, 1} (\cB \bmu _2^2)$$ be the open substack classifying $\bmu _2^2$-torsors over smooth elliptic curves. For any object $(E/T, P\rightarrow E)\in \cK_{1, 1}^\circ (\cB \bmu _2^2)$ over some scheme $T$, there is a natural action of $\bmu _2^2$ on this object given by the action on $P$.  Let $\overline {\cK }_{1,1}^\circ (\cB \bmu _2^2)$ be the rigidification.

A simple, but not concrete,  description of $\overline {\cK }^\circ_{1,1} (\cB \bmu _2^2)$ was given in Lemma \ref{Lem:tors-picx}: let $E/\cM_{1,1}$ be the universal elliptic curve. Then as discussed before, $\overline {\cK }^\circ_{1,1} (\cB \bmu _2^2) = \pic_{E/\cM_{1,1}}[X]$ where $X = (\ZZ/(2))^2$. It follows that $$\overline {\cK }^\circ_{1,1} (\cB \bmu _2^2) = \pic_{E/\cM_{1,1}}[2]\times_{\cM_{1,1}}\pic_{E/\cM_{1,1}}[2].$$ The compactification over the boundary of $\ocM_{1,1}$ is not hard to describe as well. But we wish to describe this stack in terms of more classical objects.

The stack $\overline {\cK }_{1,1}^\circ (\cB \bmu _2^2)$ can be reinterpreted as follows.

Given an object $(E/T, \pi :P\rightarrow E)\in \cK_{1, 1}^\circ (\cB \bmu _2^2)$ over a scheme $T$, the sheaf $\pi _*\mls O_P$ on $E$ has a $\bmu _2^2$-action given by the action on $P$, and therefore decomposes as a direct sum
$$
\pi _*\mls O_P = \oplus _{\chi \in (\Z/2)^2}L_\chi .
$$
Furthermore, since $P$ is a torsor over $E$ each $L_\chi $ is a locally free sheaf of rank $1$ on $E$.   The sheaves $L_\chi $ therefore define a homomorphism
$$\begin{array}{rcl}
\phi _P:(\Z/2)^2 &\rightarrow &\underline {\text{Pic}}^0(E) \simeq E, \\ 
\chi & \mapsto & [L_\chi ].\end{array}
$$

Let $\scr H(2)$ denote the stack over $\Z$ associating to any scheme $T$ the groupoid of pairs $(E, \phi )$, where $E/T$ is an elliptic curve and $\phi :(\Z/2)^2\rightarrow E[2]$ is a homomorphism of group schemes over $T$.  The above construction gives a functor
$$
\cK_{1, 1}^\circ (\cB \bmu _2^2)\rightarrow \scr H(2)
$$
and a straightforward verification shows that this functor induces an isomorphism
$$
\overline {\cK}_{1, 1}^\circ (\cB \bmu _2^2)\simeq \scr H(2).
$$

For any quotient of finite abelian groups $\pi :(\Z/2)^2\rightarrow A$, let $\scr H^A(2)$ denote the stack associating to any scheme $T$ the groupoid of pairs $(E, \phi _A)$, where $E/T$ is an elliptic curve and $\phi _A:A\rightarrow E[2]$ is an $A$-structure in the sense of Katz--Mazur \cite[1.5.1]{KM} (recall that a homomorphism $\phi _A:A\rightarrow E$ is called an $A$-structure if the effective Cartier divisor $D = \sum _{a\in A}\phi (a)$ is a subgroup scheme of $E$).  Sending such a pair $(E, \phi _A)$ to 
$$
(E, \xymatrix{(\Z/2)^2\ar[r]^-\pi & A\ar[r]^{\phi _A}& E[2]})
$$
defines a morphism of stacks
\begin{equation}\label{KMinclusion}
\scr H^A(2)\rightarrow \scr H(2).
\end{equation}
Using \cite[1.6.2]{KM} one sees that this in fact is a closed immersion.

We can describe the above in more classical notation as follows.  Define stacks $\scr Y(2)$, $\scr Y_1(2)$, and $\scr Y(1)$ by associating to any scheme $T$ the following groupoids:
\begin{enumerate}
\item [$\scr Y(2):$]  The groupoid of pairs $(E, \phi )$, where $E/T$ is an elliptic curve and $\phi :(\Z/2)^2\rightarrow E[2]$ is a $(\Z/2)^2$-generator;
\item [$\scr Y_1(2):$] The groupoid of pairs $(E, \phi )$, where $E/T$ is an elliptic curve and $\phi :(\Z/2)\rightarrow E[2]$ is a $(\Z/2)$-generator;
\item [$\scr Y(1):$] The groupoid of elliptic curves over $T$.
\end{enumerate}
Then \ref{KMinclusion} gives a closed immersion
$$
i^{(2)}:\scr Y(2)\hookrightarrow \scr H(2)
$$
corresponding to the identity map $(\Z/2)^2\rightarrow (\Z/2)^2$, three closed immersions
$$
i_{K}^{(1)}:\scr Y_1(2)\hookrightarrow \scr H(2), \ \ j=1,2,3,
$$
corresponding to the three surjections $(\Z/2)^2\rightarrow \Z/2$ (indexed by index $2$ subgroups $K\subset (\Z/2)^2$), and a closed immersion
$$
i^{(0)}:\scr Y(1)\hookrightarrow \scr H(2)
$$
corresponding to the unique surjection $(\Z/2)^2\rightarrow 0$.  The resulting map
\begin{equation}\label{coprod}
\scr Y(2)\,\sqcup\, \bigg(\scr Y_1(2)\,\sqcup\, \scr Y_1(2)\,\sqcup\, \scr Y_1(2)\bigg)\,\sqcup\, \scr Y(1)\ \ \longrightarrow \ \ \scr H(2)
\end{equation}
is then a proper surjection, which over $\Z[1/2]$ is an isomorphism.  Note in particular that the forgetful map $\scr H(2)\rightarrow \scr Y(1)$ sending $(E, \phi )$ to $E$ has degree $16$, as $\scr H(2)\rightarrow \scr Y(1)$ is flat by \ref{flattheorembis} and over $\Z[1/2]$ the map clearly has degree $16$.

Over $\mathbb{F}_2$, however, the map \ref{coprod} joins together the various components in an interesting way.  If $E/k$ is an ordinary elliptic curve over a field $k$ of characteristic $2$, then any homomorphism
\begin{equation}\label{fibermap}
\phi :(\Z/2)^2\rightarrow E[2]\simeq \Z/(2)\times \mu _2
\end{equation}
has a nontrivial kernel, and the irreducible components of $\scr Y(2)_{\mathbb{F}_2}$ are indexed by these kernels. Since a $Z/(2)^2$-structure is surjective on geometric points, the kernel cannot be the whole group.  Therefore $\scr Y(2)_{\mathbb{F}_2}$ has $3$ irreducible components: for a subgroup $K\subset (\Z/2)^2$ of index $2$ we get an irreducible component $\scr Y(2)^K_{\mathbb{F}_2}$ corresponding to pairs $(E, \phi )$, where the map \ref{fibermap} has kernel $K$. 
% There is also an irreducible component $\scr Y(2)^\emptyset _{\mathbb{F}_2}$ classifying pairs $(E, \phi )$, where \ref{fibermap} is the zero map.  
Note also that for any ordinary elliptic curve $E/k$ over a field $k$, the fiber product of the diagram
$$
\begin{CD}
@. \scr Y(2)^K_{\mathbb{F}_2}\\
@. @VVV \\
\Sp (k)@>E>> \scr Y(1)_{\mathbb{F}_2}
\end{CD}
$$
has length $2$ over $k$.

Similarly the reduction $\scr Y_1(2)_{\mathbb{F}_2}$ has two irreducible components 
$$
\scr Y_1(2)_{\mathbb{F}_2} = \scr Y_1(2)'\cup \scr Y_1(2)^\emptyset ,
$$
where $\scr Y_1(2)'$ (resp. $\scr Y_1(2)^\emptyset $) classifies pairs $(E, \phi :\Z/2\rightarrow E[2])$ where the map $\phi $ is injective (resp. the zero map). This time only the fiber of  $\scr Y_1(2)'$ over a $k$-point of $\scr Y(1)_{\mathbb{F}_2}$ has length 2.

The closed fiber $\scr H(2)_{\mathbb{F}_2}$ then has four irreducible components.  For a subgroup $K\subset (\Z/2)^2$ of index $2$, we have an irreducible component $Z_K$ which is set-theoretically identified with $\scr Y(2)^K_{\mathbb{F}_2}$ as well as the component $\scr Y_1(2)'$ given by  the inclusion $i_K^{(1)}$.  The fourth component $Z_\emptyset $ is  set-theoretically identified with  $\scr Y(1) _{\mathbb{F}_2}$ via $i^{(0)}$, and with all the components $\scr Y_1(2)^\emptyset $ via the inclusions $i_K^{(1)}$. 

Over the point of $\scr Y(1)(\overline {\mathbb{F}}_2)$ given by the supersingular elliptic curve $E$, there is only one homomorphism (the zero map)
$$
\phi :(\Z/2)^2\rightarrow E[2],
$$
and hence the four irreducible components of $\scr H(2)_{\mathbb{F}_2}$ all intersect at this point (and nowhere else).  Note also that over the ordinary locus of $\scr Y(1)$, each of the components $Z_K$ has length $4$  over $\scr Y(1)_{\mathbb{F}_2}$ as does the component $Z_\emptyset $ (so $\scr H(2)_{\mathbb{F}_2}$ has length $16$ over $\scr Y(1)_{\mathbb{F}_2}$, as we already knew).

%=-=-=-
The reduced fiber over $\FF_2$ looks roughly like this:

%\begin{center}
\includegraphics[scale=.5]{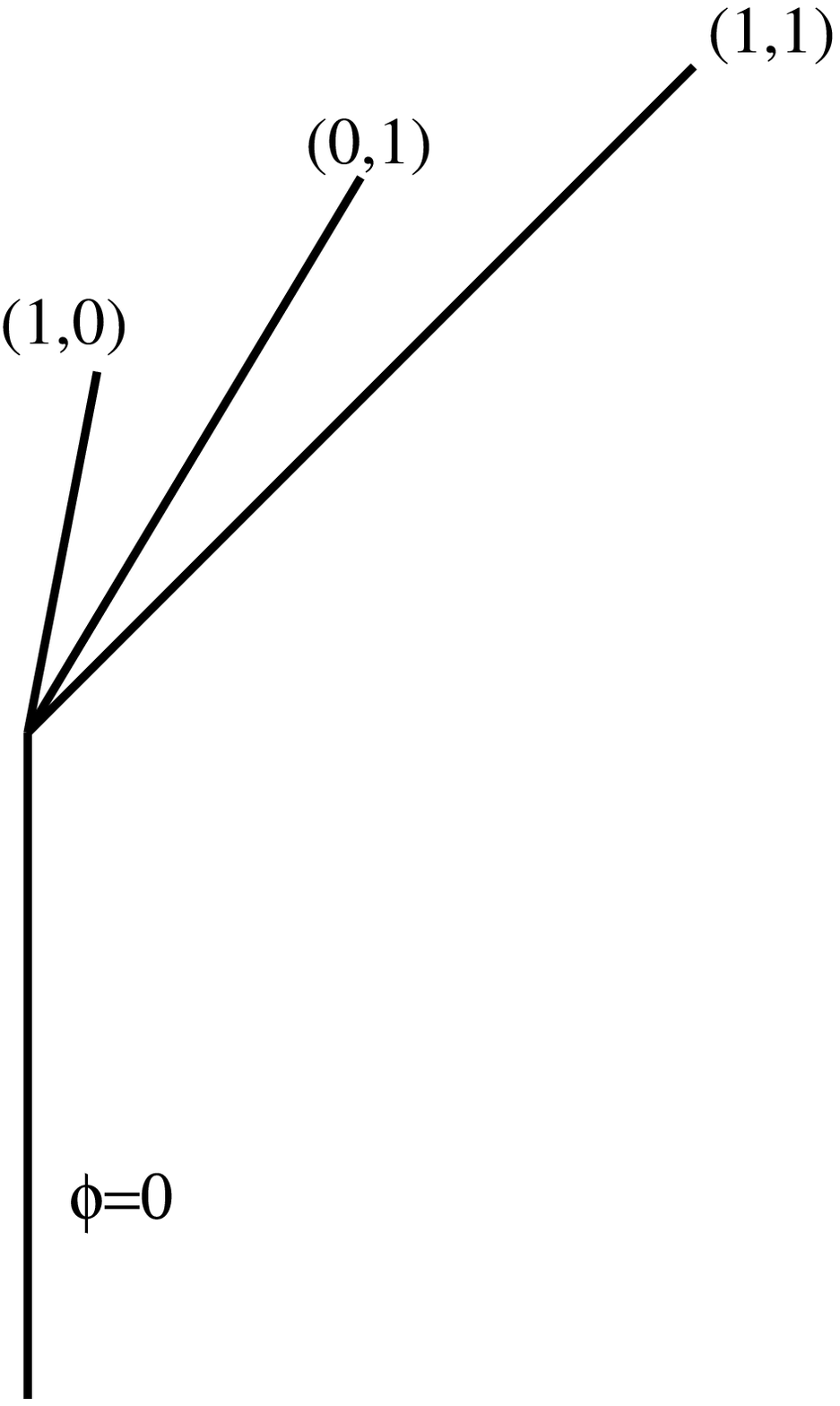}
%\end{center}

The following figure schematically describes the main component $\scr Y(2)$:
\begin{center}
\includegraphics[scale=.5]{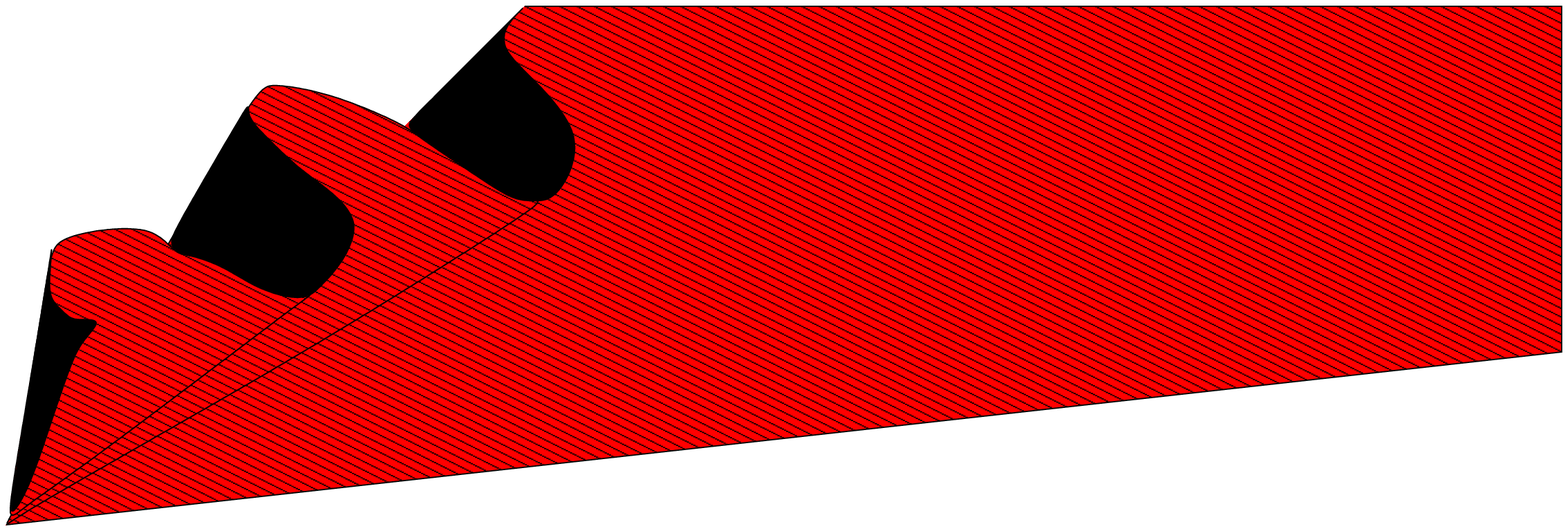}
\end{center}
Here is the component $\scr Y_1(2)^K$ for $K$ generated by $(1,1)$:
\begin{center}
\includegraphics[scale=.5]{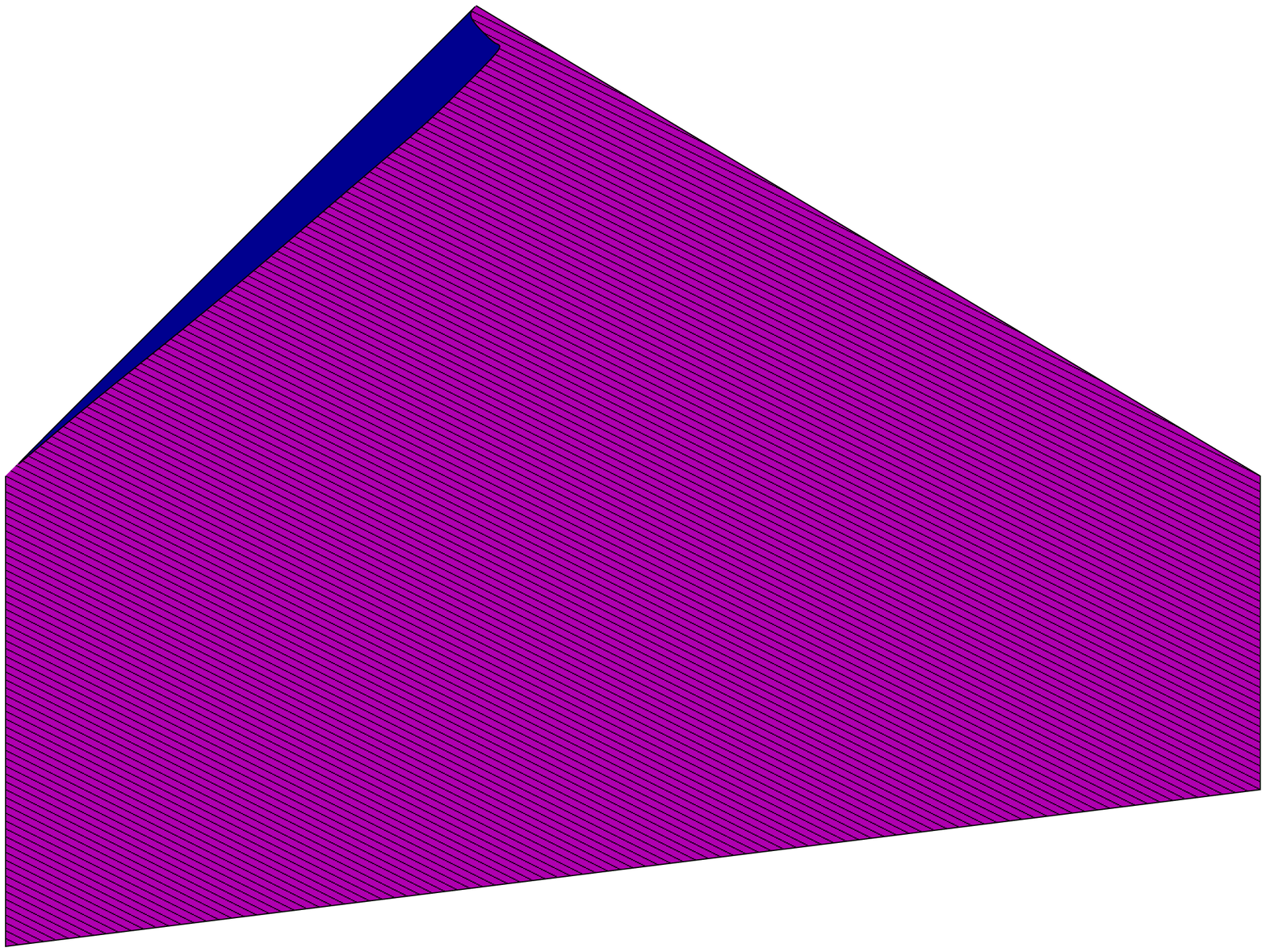}
\end{center}
There are of course two more components  corresponding to the two other choices for a subgroup $K$ of order 2.
In addition  of course there is $\scr Y(1)$:

%\begin{center}
\includegraphics[scale=.5]{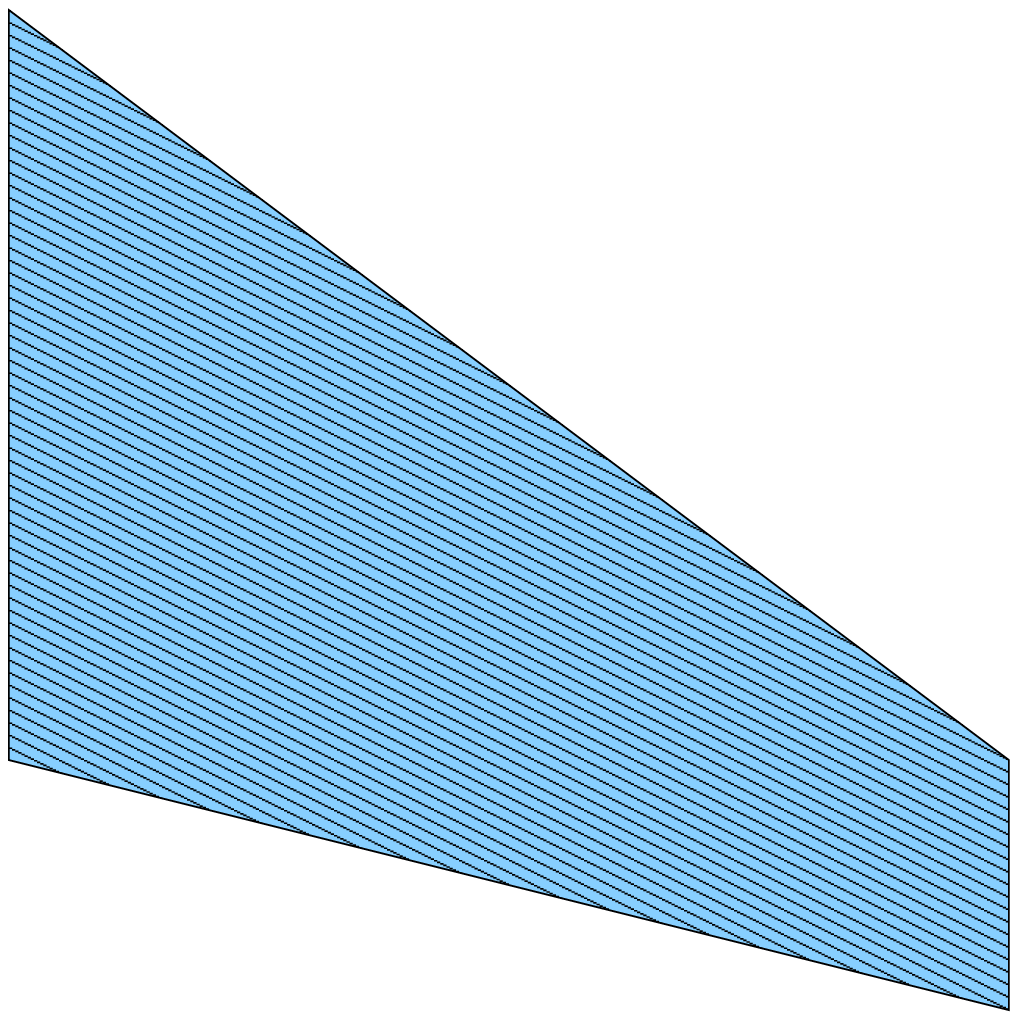}
%\end{center}

The way the main component $\scr Y(2)$ meets 
$\scr Y_1(2)^K$ is described 
as follows:   
\begin{center}
\includegraphics[scale=.5]{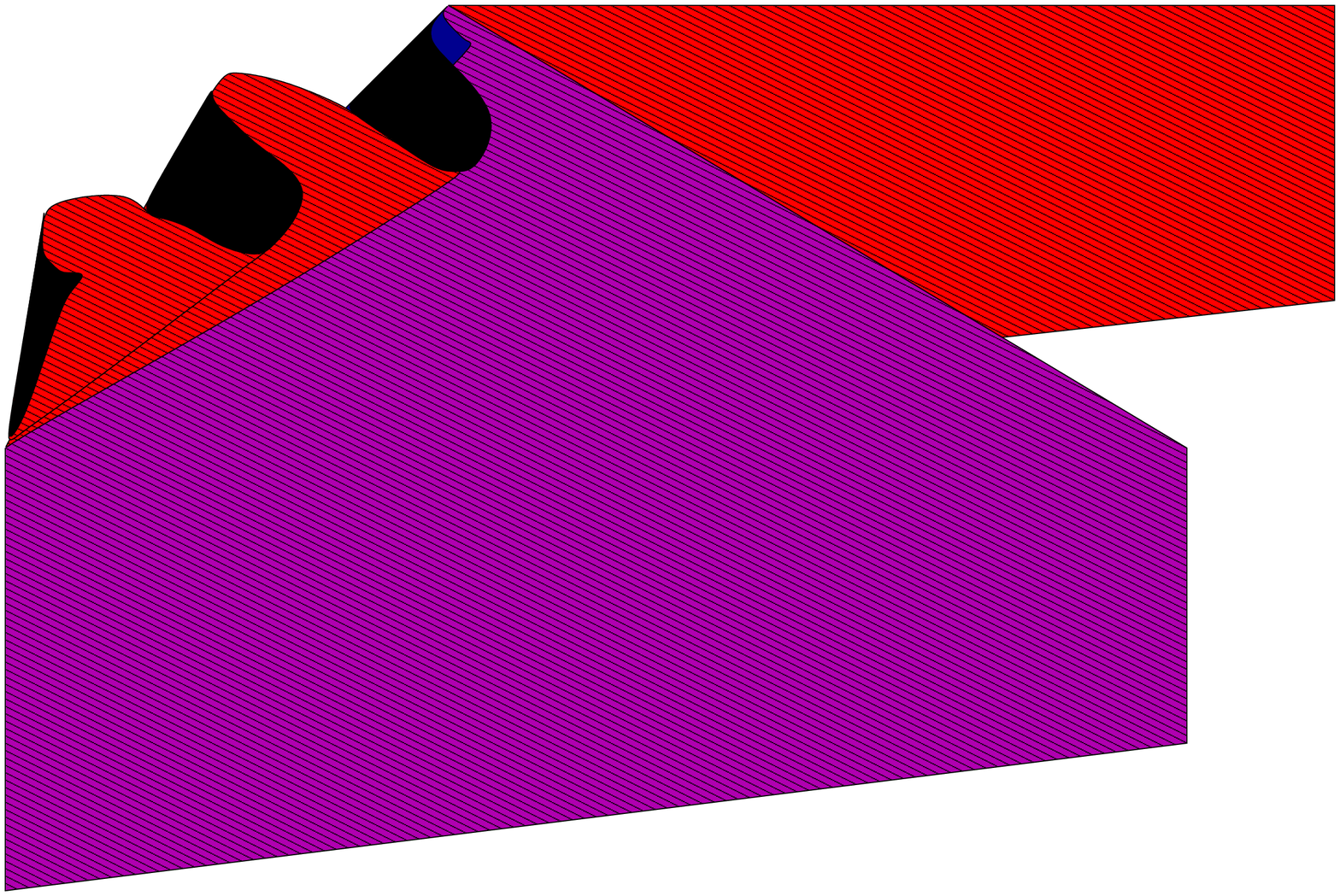}
\end{center}
and the other two are similar, of course along the corresponding
component $\scr Y(2)^K_{\FF_2}$. 

 The main component $\scr Y(2)$ meets $\scr Y(1)$ only at the supersingular
 point, but $\scr Y_1(2)^K$ meets $\scr Y(1)$ along the component corresponding
 to $\phi=0$:
\begin{center}
\includegraphics[scale=.5]{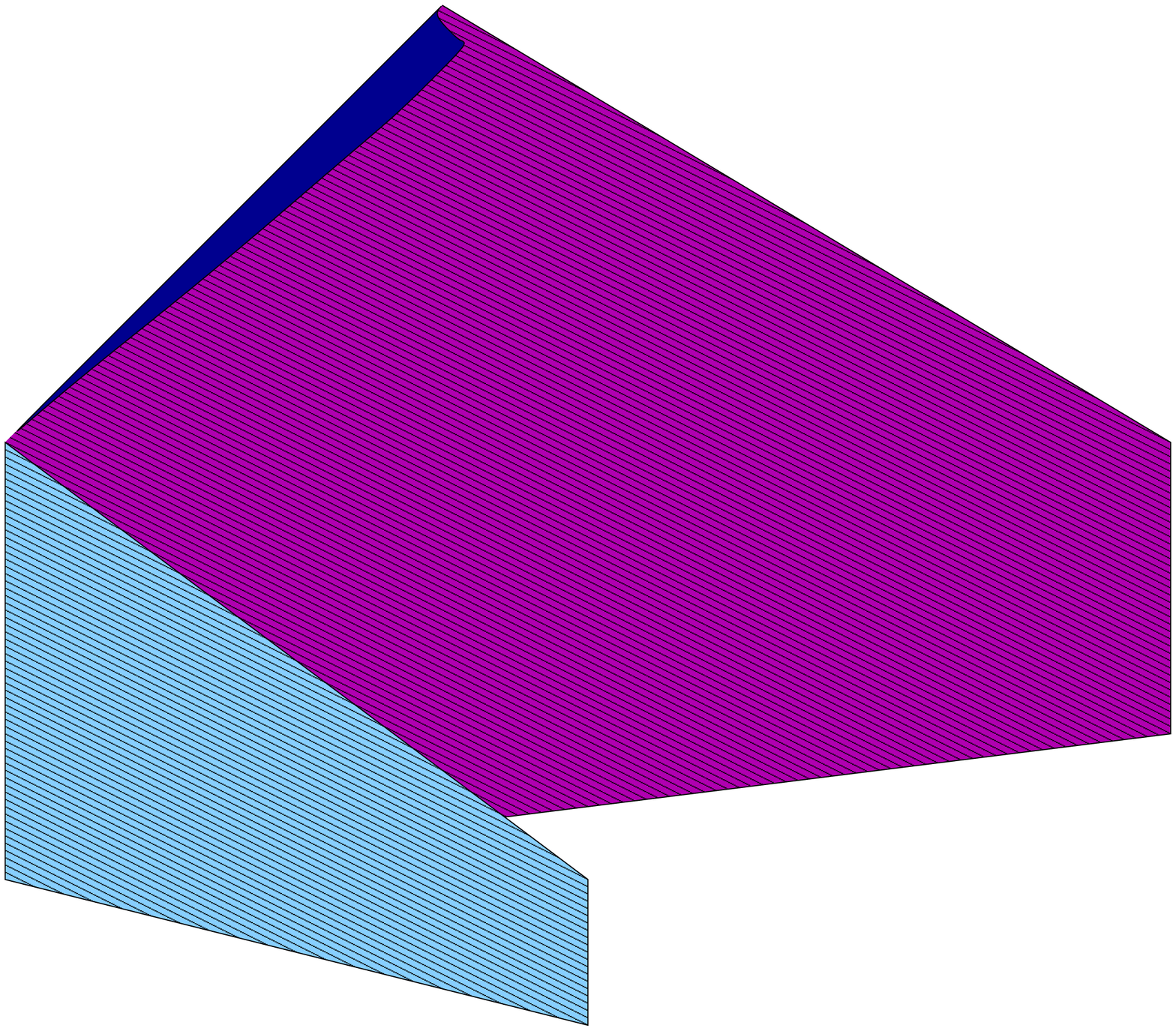}
\end{center}
The complete picture is something like this:
\begin{center}
\includegraphics[scale=.5]{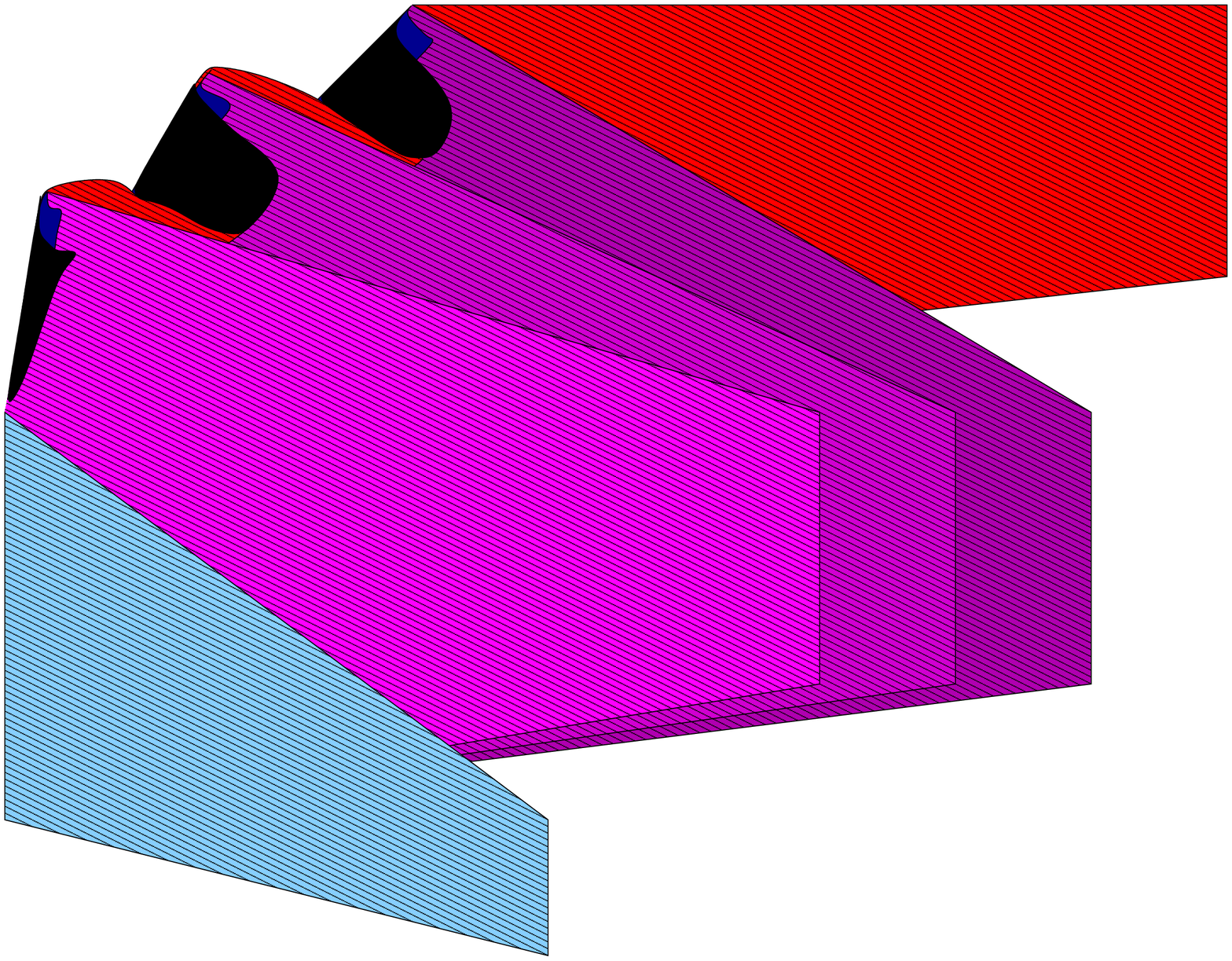}
\end{center}

%=-=-=-

An important question is that of describing explicitly the ``most important" part of $\overline {\cK }_{g,n} (\cB \bmu _m^{2g})$, in this  example ${\scr X}(2)$. In genus 1, Katz and Mazur gave an ideal solution: the moduli stack of generalized elliptic curves with $(\ZZ/m)^2$ structure, which is a closed substack of our  $\overline {\cK }_{1,1} (\cB \bmu _m^{2})$, is regular and flat over $\ocM_{1,1}$. But for higher genus a satisfactory conclusion is missing:  in \cite[Section 1.9]{KM} Katz and Mazur proposed a solution using norms, but it  was shown to be ill behaved in  \cite[Appendix A]{Chai-Norman}. It still may be of interest to study  the closure of the $\QQ$-stack $\ocM_{g}^{(m)}$  of level $m$ curves in  $\overline {\cK }_{g} (\cB \bmu _m^{2g})$.

\begin{appendix}

\section{Twisted  curves and log curves\\ by Martin Olsson}\label{AppendixA}

In \cite{logtwisted} we introduced a notion of log twisted curve, and proved that the category of log twisted curves (over some scheme $S$) is naturally equivalent to the category of twisted Deligne-Mumford curves over $S$.  The definition of a log twisted curve in loc. cit. included various assumptions about certain integers being invertible on the base scheme.  In this appendix we show that if we remove these assumptions in the definition of log twisted curve, we obtain a category equivalent to the category of twisted  curves.

\begin{defn}[{\cite[3.1]{M-O}}]\label{locallyfree}
Let $X$ be an Artin stack.

(i) A fine log structure $\mc M$ on $X$ is called \emph{locally free}
if for every geometric point $\bar x\rightarrow X$ the monoid
$\overline {\mc M}_{\bar x}:= \mc M_{\bar x}/\mc O_{X, \bar x}^*$
is isomorphic to $\mathbb{N}^r$ for some $r$.

(ii) A morphism $\mc M\rightarrow \mc N$ of locally free log
structures on $X$ is called \emph{simple} if for every geometric
point $\bar x\rightarrow X$ the monoids $\overline {\mc M}_{\bar
x}$ and $\overline {\mc N}_{\bar x}$ have the same rank,  the
morphism $\varphi :\overline {\mc M}_x\rightarrow \overline {\mc
N}_{\bar x}$ is injective, and for every irreducible element $f\in
\overline {\mc N}_{\bar x}$ there exists an irreducible element
$g\in \overline {\mc M}_{\bar x}$ and a positive integer $n$ such
that $\varphi (g) = nf$.
\end{defn}

\begin{rem} This differs from \cite[1.5]{logtwisted}  where the integer $n$ in (\ref{locallyfree} (ii)) was assumed invertible in $k(\bar x)$.
\end{rem}

Let $f:C\rightarrow S$ be a nodal curve over a scheme $S$.  As discussed in \cite[\S 3]{logtwisted}, there exist canonical log structures $\mc M_C$ and $\mc M_S$ on $C$ and $S$ respectively, and an extension of $f$ to a log smooth morphism
\begin{equation*}
(C,\widetilde  {\mc M}_C)\rightarrow (S, \mc M_S).
\end{equation*}

\begin{defn}\label{1.5} A \emph{$n$-pointed log twisted curve}  over a
scheme $S$ is a collection of data
\begin{equation*}
(C/S,  \{\sigma _i, a_i\}_{i=1}^n, \ell :\mc M_S\hookrightarrow
\mc M_S'),
\end{equation*}
where $C/S$ is a nodal curve, $\sigma _i:S\rightarrow C$ are
sections, the $a_i$ are positive  integer-valued locally constant functions on $S$,  and
$\ell :\mc M_S\hookrightarrow \mc M_S'$ is a simple morphism of
log structures on $S$, where $\mc M_S$ denotes the canonical log
structure on $S$ mentioned above.
\end{defn}

Let  $(C/S, \{\sigma
_i, a_i\}, \ell :\mc M_S\hookrightarrow \mc M_S')$ be a log twisted  curve over $S$, and let
\begin{equation}\label{logmap}
(C, \mc M_C)\longrightarrow (S, \mc M_S)
\end{equation}
be the morphism of log schemes obtained from the pointed curve
$$(C/S, \{\sigma _i\})$$ as in \cite[3.10]{logtwisted} (note that $\mc M_C$ is not equal to $\widetilde {\mc M}_C$ mentioned above, as $\mc M_C$ also takes into account the marked points). We construct a
twisted $n$--pointed curve $(\mc C, \{\Sigma _i\})/S$ from the log
twisted curve $(C/S, \{\sigma _i, a_i\}, \ell :\mc
M_S\hookrightarrow \mc M_S')$ as in \cite[\S 4]{logtwisted}.

Define $\mc C$ to be the fibered category over $S$ which to any
$h:T\rightarrow S$ associated the groupoid of data consisting of a
morphism $s:T\rightarrow C$ over $h$ together with a commutative
diagram of locally free log structures on $T$
\begin{equation}\label{3.15.2}
\begin{CD}
h^*\mc M_S@>\ell >> h^*\mc M_S'\\
@VVV @VV\tau V \\
s^*\mc M_C@>k>> \mc M_C',
\end{CD}
\end{equation}
where:

\noindent (\arabic{section}.\arabic{subsection} (i)) The map $k$
is a simple, and for every geometric point $\bar t\rightarrow T$,
the map $\overline {\mc M}'_{S, \bar t}\rightarrow \overline {\mc
M}_{C, \bar t}'$ is either an isomorphism, or of the form
$\mathbb{N}^r\rightarrow \mathbb{N}^{r+1}$ sending $e_i$ to $e_i$
for $i<r$ and $e_r$ to either $e_r$ or $e_{r}+e_{r+1}$.

\noindent (\arabic{section}.\arabic{subsection} (ii)) For every
$i$ and geometric point $\bar t\rightarrow T$ with image under $s$
in $\sigma _i(S)\subset C$, the group
\begin{equation}\label{3.14.3}
\text{Coker}(\overline {\mc M}_{S, \bar t}^{\prime
\text{gp}}\oplus \overline {\mc M}_{C, \bar t
}^{\text{gp}}\longrightarrow \overline {\mc M}_{C, \bar t}^{\prime
\text{gp}})
\end{equation}
is a cyclic group of order $a_i(\bar t)$.

For every $i$, define $\Sigma _i\subset \mc C$ to be the substack
classifying morphisms $s:T\rightarrow C$ which factor through
$\sigma _i(S)\subset C$ and diagrams (\ref{3.15.2}) such that for
every geometric point $\bar t\rightarrow T$ the image of
\begin{equation*}
(\mc M_{C, \bar t}' -- \tau (h^*\mc M_{S, \bar t}'))\rightarrow \mc
O_{T, \bar t}
\end{equation*}
is zero.

\begin{prop} The data $(\mc C, \{\Sigma _i\}_{i=1}^n)$ is a twisted $n$-pointed  curve.
\end{prop}
\begin{proof} This follows from the argument given in \cite[\S 4]{logtwisted}.
\end{proof}

The main result of this appendix is the following:

\begin{thm}\label{Bthm} Let $S$ be a scheme.  The functor
\begin{equation}\label{thefunctor}
(\text{\rm $n$-pointed log twisted curves})\rightarrow (\text{\rm twisted $n$-pointed  curves})
\end{equation}
sending  $(C/S, \{\sigma
_i, a_i\}, \ell :\mc M_S\hookrightarrow \mc M_S')$ to the twisted $n$-pointed  curve $(\mc C, \{\Sigma _i\}_{i=1}^n)$ constructed above is an equivalence of categories.  Moreover, this equivalence of categories is compatible with base change $S'\rightarrow S$.
\end{thm}

We prove this theorem  below.  Before giving the proof, however, let us record the main consequences of this result that we will use.

 Fix integers
$g$ and $n$, and  let $\scr S_{g, n}$ denote the fibered category
over $\Z$ which to any scheme $S$ associates the groupoid of all
(not necessarily stable) $n$--pointed genus $g$ nodal curves
$C/S$. The stack $\scr S_{g, n}$ is algebraic, and as explained in
section \cite[\S 5]{logtwisted} the substack $\scr S_{g, n}^0\subset \scr
S_{g, n}$ classifying smooth curves defines a log structure $\mc
M_{\scr S_{g, n}}$ on $\scr S_{g, n}$.  The same argument used in \cite[\S 5]{logtwisted} yields the following theorem:

\begin{thm}\label{3.19} Let $\fM _{g, n}^{\text{\rm tw}}$ denote the fibered
category over $\Z$ which to any scheme $T$ associates the groupoid
of $n$--marked genus $g$ twisted curves $(\mc C, \{\Sigma _i\})$
over $T$. Then $\fM _{g, n}^{\text{\rm tw}}$ is a smooth Artin
stack, and the natural map
\begin{equation}\label{3.19.1}
\pi :\fM _{g, n}^{\text{\rm tw}}\longrightarrow \scr S_{g, n}
\end{equation}
sending $(\mc C, \{\Sigma _i\})$ to its coarse moduli space with
the marked points induced by the $\Sigma _i$ is representable by
tame  stacks. Moreover, there is a natural locally free
log structure $\mc M _{\fM _{g, n}^{\text{\rm tw}}}$ on $\fM _{g,
n}^{\text{\rm tw}}$.
\end{thm}

\begin{rem}
Consider a field $k$ and an object $(\mc C,  \{\Sigma _i\})\in \mc
M_{g, n}^{\text{\rm tw}}(k)$.  Let $(C, \{\sigma _i\})$ be the
coarse moduli space, and let $R$ be a versal deformation space for
the object $(C, \{\sigma _i\})\in \scr S_{g, n}(k)$.  Let $q_1,
\dots , q_m\in C$ be the nodes and let $r_i$ be the order of the
stabilizer group of a point of $\mc C$ lying above $q_i$.  As in
\cite[1.5]{D-M}, there is a smooth divisor $D_i\subset \Sp (R)$
classifying deformations where $q_i$ remains a node.  In other
words, if $t_i\in R$ is an element defining $D_i$ then in an
\'etale neighborhood of $q_i$ the versal deformation $\widetilde
C\rightarrow \Sp (R)$ of $(C, \{\sigma _i\})$ is isomorphic to
\begin{equation*}
\Sp (R[x,y]/(xy-t_i)).
\end{equation*}
It follows from the argument in \cite[\S 5]{logtwisted} that 
the fiber product
\begin{equation*}
\fM _{g, n}^{\text{tw}}\times _{\scr S_{g, n}}\Sp (R)
\end{equation*}
is isomorphic to the stack-theoretic quotient of 
\begin{equation*}
\Sp (R[z_1, \dots , z_m]/(z_1^{r_1}-t_1, \dots , z_m^{r_m}-t_m))
\end{equation*}
by the action of $\mu _{r_1}\times \cdots \mu _{r_m}$ for which $(\zeta _1, \dots , \zeta _r)\in \mu _{r_1}\times \cdots \mu _{r_m}$ sends $z_i$ to $\zeta _iz_i$. In particular $\fM _{g, n}^{\text{tw}}$ is flat over $\scr S_{g, n}$, and hence also flat over $\Z$.
\end{rem}

We also get a generalization of \cite[1.11]{logtwisted}:
\begin{cor}\label{1.7} For any integer $N>0$, let $\fM _{g, n}^{\text{\rm tw},
\leq N}$ denote the substack of $\fM _{g, n}^{\text{\rm tw}}$
classifying $n$--pointed genus $g$ twisted  curves such that the
order of the stabilizer group at every point is less than or equal
to $N$. Then $\fM _{g, n}^{\text{\rm tw}, \leq N}$ is an open
substack of $\fM _{g, n}^{\text{\rm tw}}$ and the map $\fM _{g,
n}^{\text{\rm tw}, \leq N}\rightarrow \scr S_{g, n}$ is of finite
type and quasi-finite.
\end{cor}

\begin{rem}\label{curveconstruct} Let $R$ be a discrete valuation ring with uniformizer $\pi $ and separably closed residue field, and let $C/R$ be a nodal curve.  Let $p_1, \dots, p_r$ be the nodes of $C$ in the closed fiber.  Then we have a log smooth morphism
$$
(C, \mc M_C)\rightarrow (\Sp (R), \mc M_R),
$$
where the log structure $\mc M_R\rightarrow R$ admits a chart $\mathbb{N}^r\rightarrow R$ such that the image of the $i$-th standard generator is equal to $\pi ^{l_i}$, where $l_i\in \mathbb{N}\cup \{\infty \}$ is an element such that in an \'etale neighborhood of $p_i$ the curve $C_i$ is isomorphic to 
$$
\Sp (R[x,y]/(xy-\pi ^{l_i}),
$$
where by convention if $l_i=\infty $ we set $\pi ^{l_i} = 0$.

Now assume some $l_i$ is finite.  Then we obtain a twisted curve by taking the stack corresponding to the morphism of log structures $\mc M_R\rightarrow \mc M_{R}'$, where $\mc M_R'$ is the log structure associated to the map $\mathbb{N}^r\rightarrow R$ sending $e_j$ to $\pi ^{l_j}$ for $j\neq i$, and $e_i$ to $\pi $.
\end{rem}

\subsection*{Proof of \ref{Bthm}}\label{proofofBthm}

The proof of \ref{Bthm} follows the same outline as the proof of \cite[1.9]{logtwisted}.  We review the argument here indicating the necessary changes for this more general setting.

\begin{defn}[{\cite[3.3]{logtwisted}}]\label{esssems}
A log smooth morphism of fine log  schemes $f:(X, \mc M_X)\rightarrow (S, \mc M_S)$ is
\emph{essentially semi-stable} if for each geometric point $\bar
x\rightarrow X$ the monoids $(f^{-1}\overline {\lsm {}}_S)_{\bar
x}$ and $\overline {\lsm {}}_{X, \bar x}$ are free monoids, and if
for suitable isomorphisms $(f^{-1}\overline {\lsm {}}_S)_{\bar
x}\simeq \mathbb{N}^r$ and $\overline {\lsm {}}_{X, \bar x}\simeq
\mathbb{N}^{r+s}$ the map
\begin{equation*}
(f^{-1}\overline {\lsm {}}_S)_{\bar x}\rightarrow \overline {\lsm
{}}_{X, \bar x}
\end{equation*}
is of the form
\begin{equation}\label{chartformula}
e_i\mapsto \left \{\begin{array}{cl}
    e_i & \mbox{if $i\neq r$}\\
    e_r+e_{r+1}+\cdots +e_{r+s} & \mbox{if $i=r,$}
    \end{array}\right.
\end{equation}
where $e_i$ denotes the $i$-th standard generator of
$\mathbb{N}^r$.
\end{defn}

With notation as in \ref{esssems}, if $U\rightarrow X$ is a smooth surjection of schemes and $\mc M_U$ denotes the pullback of $\mc M_X$, then it follows immediately from the definition that $(X, \mc M_X)\rightarrow (S, \mc M_S)$ is essentially semistable if and only if the morphism $(U, \mc M_U)\rightarrow (S, \mc M_S)$ is essentially semistable (the property of being semistable is local in the smooth topology on $X$).  It follows that the notion of a morphism being essentially semistable extends to Artin stacks: If  $(\scr X, \mc M_{\scr X})$ is an Artin stack with a fine log structure, and $f:(\scr X, \mc M_{\scr X})\rightarrow (S, \mc M_S)$ is a morphism to a fine log scheme, then $f$ is \emph{essentially semistable} if for some smooth surjection $U\rightarrow \scr X$ with $U$ a scheme the induced morphism of log schemes
\begin{equation}\label{testmap}
(U, \mc M_{\scr X}|_U)\rightarrow (S, \mc M_S)
\end{equation}
is essentially semistable.  As usual if for some smooth surjection $U\rightarrow \scr X$ the morphism \ref{testmap} is essentially semistable then for any smooth surjection $V\rightarrow \scr X$ the induced morphism $(V, \mc M_{\scr X}|_V)\rightarrow (S, \mc M_S)$ is essentially semistable.

As explained in \cite[\S 3]{logtwisted}, if $(X, \mc M_X)\rightarrow (S, \mc M_S)$ is an essentially semistable morphism of log schemes, then for any geometric point $\bar s\rightarrow S$ there is a canonical induced map
\begin{equation*}
s_{X_{\bar s}}\colon \{\text{singular points of $X_{\bar s}$}\}\rightarrow \text{Irr}(\overline {\mc M}_{S, \bar s}),
\end{equation*}
where $\text{Irr}(\overline {\mc M}_{S, \bar s})$ denotes the set of irreducible elements in the monoid $\overline {\mc M}_{S, \bar s}$.

If $\scr Y$ is an Artin stack over a field $k$, let $\pi _0(\scr Y_{\text{sing}})$ denote the set of connected components of the complement of the maximal open substack $\scr U\subset \scr Y$ which is smooth over $k$.  Note that if $\scr Y'\rightarrow \scr Y$ is a smooth morphism of $k$-stacks, then there is an induced map
$$
\pi _0(\scr Y'_{\text{sing}})\rightarrow \pi _0(\scr Y_{\text{sing}}).
$$
If $\scr Y'\rightarrow \scr Y$ is also surjective then this map is also surjective.

\begin{defn} Let $f:(\scr X, \mc M_{\scr X})\rightarrow (S, \mc M_S)$ be an essentially semistable morphism from a log Artin stack to a fine log scheme, and let $\bar s\rightarrow S$ be a geometric point.  We say that $f$ is \emph{special at $\bar s$} if for any smooth surjection $U\rightarrow \scr X$ with $U$ a scheme the map
$$
s_{U_{\bar s}}\colon \{\text{singular points of $U_{\bar s}$}\}\rightarrow \text{Irr}(\overline {\mc M}_{S, \bar s})
$$
factors through the composite (which is surjective)
$$
\begin{CD}
\{\text{singular points of $U_{\bar s}$}\}@>>> \pi _0(U_{\bar s, \text{sing}})@>>> \pi _0(\scr X_{\bar s, \text{sing}})
\end{CD}
$$
to give an isomorphism
$$
 \pi _0(\scr X_{\bar s, \text{sing}})\rightarrow \text{Irr}(\overline {\mc M}_{S, \bar s}).
 $$
\end{defn}

\begin{thm}[{Generalization of \cite[3.6]{logtwisted}}]\label{canlogthm} Let $(f:\mc C\rightarrow S, \{\Sigma _i\})$ be an $n$-marked  twisted  curve.  Then there exist log
structures $\widetilde {\mc M}_{\mc C}$ and $\mc M_S'$ on $\mc C$
and $S$ respectively, and a special morphism
\begin{equation*}
(f, f^b):(\mc C, \widetilde {\mc M}_{\mc C})\longrightarrow (S,
\mc M_S').
\end{equation*}
Moreover, the datum $(\widetilde {\mc M}_{\mc C}, \mc M_S', f^b)$
is unique up to unique isomorphism.
\end{thm}
\begin{proof}
The proof will be in several steps (ending in paragraph following \ref{A18}).

The uniqueness statement follows as in \cite[3.6]{logtwisted} from the argument proving the uniqueness in \cite[2.7]{Olsson-tohoku}.

Given the uniqueness, to prove existence we may work \'etale locally on $S$, and by a limit argument as in the proof of \cite[3.6]{logtwisted} may even assume that $S$ is the spectrum of a strictly henselian local ring.  Let $\bar s\in S$ be the closed point.

Let $p_1, \dots, p_n\in C$ be the nodes of the closed fiber, and choose for each $i=1, \dots , n$ an affine open set $U_i\subset C$ containing $p_i$ and no other nodes.  Let $\mc U_i\subset \mc C$ denote the inverse image of $U_i$. 

\begin{lem}\label{A14} For any quasi-coherent sheaf $\mc F$ on $\mc U_i$, we have $$H^j(\mc U_i, \mc F) = 0$$ for $j>0$.
\end{lem}
\begin{proof}
This follows from the same argument proving \ref{cohlemb}.
\end{proof}

Let $t_i\in \mc O_S$ be an element such that the fiber product
$$
\mc C\times _C\Sp (\mc O_{C, \bar p_i})
$$
is isomorphic to
$$
[\Sp (\widetilde {\mc O_S[z, w]/(zw-t_i)})/\mu _n]
$$
as in (\ref{A1} (v)).  Let $\mc M_S^i$ be the log structure on $S$ associated to the morphism $\mathbb{N}\rightarrow \mc O_S$ sending $1$ to $t_i$

\begin{defn} Let $f:X\rightarrow S$ be a morphism of schemes.  A \emph{$t_i$-semistable log structure on $X$} is a pair $(\mc M_X, f^b)$, where $\mc M_X$ is a locally free log structure on $X$ and $f^b:f^*\mc M_S^i\rightarrow \mc X$ is a morphism of log structures, such that the following hold:
\begin{enumeratei}
\item The morphism of log schemes
$$
(f, f^b):(X, \mc M_X)\rightarrow (S, \mc M_S^i)
$$
is log smooth;
\item For every geometric point $\bar x\rightarrow X$ the induced map of free monoids
$$
\mathbb{N}\rightarrow \overline {\mc M}_{S, f(\bar x)}^i\rightarrow \overline {\mc M}_{X, \bar x}
$$
is the diagonal map.
\end{enumeratei}
\end{defn}

\begin{rem}\label{smoothdes} By \cite[3.4]{logtwisted}, if $X\rightarrow S$ admits a $t_i$-semistable log structure, then \'etale locally on $X$ there exists a smooth morphism
$$
X\rightarrow \Sp (\mc O_S[z,w]/(zw-t_i)).
$$
\end{rem}

To prove \ref{canlogthm}, it suffices by the same argument used in \cite[proof of 3.6]{logtwisted} to show that there exists a $t_i$-semistable log structure on each $\mc U_i$.

\begin{lem}\label{localexist} Let $X\rightarrow \mc U_i$ be a smooth morphism with $X$ a scheme.  Then \'etale locally on $X$ there exists a $t_i$-semistable log structure.
\end{lem}
\begin{proof}
It suffices to prove the existence in some \'etale neighborhood of a geometric point $\bar x\rightarrow X$ mapping to the node $p_i$ of $C$.  Making an \'etale base change on $C$, it therefore suffices to show that if
$$
g:X\rightarrow [\Sp (\mc O_S[z,w]/(zw-t_i)/\mu _n]
$$
is a smooth morphism then there exists a $t_i$-semistable log structure on $X$.  For this note that there is a log structure $\mc M$ on $[\Sp (\mc O_S[z,w]/(zw-t_i)/\mu _n]$ and a morphism of log stacks
\begin{equation}\label{stackmap}
([\Sp (\mc O_S[z,w]/(zw-t_i)/\mu _n], \mc M)\rightarrow (S, \mc M_S)
\end{equation}
induced by the commutative diagram
$$
\begin{CD}
\mathbb{N}^2@>\beta >> \mc O_S[z,w]/(zw-t_i)\\
@A\Delta AA @AAA \\
\mathbb{N}@>1\mapsto t_i>> \mc O_S
\end{CD}
$$
where $\beta $ sends $(1, 0)$ to $z$ and $(0,1)$ to $w$.  By \cite[5.23]{Olsson-ENS}, the morphism \ref{stackmap} is log \'etale.  It follows that the pullback $g^*\mc M$ on $X$ is a $t_i$-semistable log structure on $X$. 
\end{proof}

Let $SS_{t_i}$ denote the presheaf on the lisse-\'etale site $\text{Lis-Et}(\mc U_i)$ of $\mc U_i$ which to any smooth morphism $X\rightarrow \mc U_i$ associates the set of isomorphism classes of $t_i$-semistable log structures on $X$.  As  in the proof of \cite[3.18]{Olsson-tohoku} a $t_i$-semistable log structure admits no nontrivial automorphisms, and hence $SS_{t_i}$ is in fact a sheaf.  In fact, $SS_{t_i}$ is a torsor under a certain sheaf of abelian groups which we now describe.

For any smooth morphism $X\rightarrow \mc U_i$, there exists by \ref{localexist} and \ref{smoothdes} \'etale locally on $X$ a smooth morphism
$$
\rho :X\rightarrow \Sp (\mc O_S[z,w]/(zw-t_i)).
$$
As explained in \cite[3.12]{Olsson-tohoku} the ideal $J:= (z,w)\mc O_X\subset \mc O_X$ is independent of the choice of the smooth morphism $\rho $. It follows that these locally defined sheaves of ideals descend to a sheaf of ideals $\mc J\subset \mc O_{\mc U_i}$.  Let $\mc D\subset \mc U_i$ be the closed substack defined by this sheaf of ideals. The local description (\ref{A1} (v)) of the stack $\mc U_i$ implies that there is an isomorphism
$$
\mc D\simeq \cB\mu _n\times _{\Sp (\Z)}\Sp (\mc O_S/(t_i)).
$$  

Let $K_{t_i}\subset \mc O_S$ be the kernel of multiplication by $t_i$ on $\mc O_S$, and let $K_{t_i}^{\mc U_i}$ denote the kernel of multiplication by $t_i$ on $\mc O_{\mc U_i}$.  Note that since $\mc U_i$ is flat over $S$ the sheaf $K_{t_i}^{\mc U_i}$ is equal to the pullback of $K_{t_i}$.  Let $\mc Z\subset \mc U$ be the closed substack defined by $K_{t_i}^{\mc U_i}\cdot \mc J$. Define
$$
G:= \text{Ker}(\mc O_{\mc U_i}^*\rightarrow \mc O_{\mc Z}^*),
$$
and let $G_2\subset \mc O_{\mc U_i}^*$ denote the subsheaf of units $u$ such that $ut_i=t_i$.  There is a natural inclusion $G\subset G_2$, and as explained in the proof of \cite[3.18]{Olsson-tohoku} the sheaf $SS_{t_i}$ is naturally a torsor under $G_2/G$.    To prove the existence of a $t_i$-semistable log structure on $\mc U_i$ we therefore must show that the class of this torsor
$$
o\in H^1(\mc U_i, G_2/G)
$$
is zero.

\begin{lem}\label{A18} The map $H^1(\mc U_i, G_2/G)\rightarrow H^1(\mc D, \mc O_{\mc D}^*)$ induced by the composite
$$
G_2\subset \mc O_{\mc U_i}^*\rightarrow \mc O_{\mc D}^*
$$
is injective.
\end{lem}
\begin{proof}
This follows from the same argument proving \cite[3.7]{logtwisted} (and using \ref{A14}).
\end{proof}
The image of the class $o$ in $H^1(\mc D, \mc O_{\mc D}^*)$ corresponds to an invertible sheaf $\mc L$ on $\mc D$.  As in \cite[proof of 3.16]{Olsson-tohoku} this invertible sheaf $\mc L$ can be described as follows.  Consider the inclusion
$$
\cB\mu _n\hookrightarrow [\Sp (\widetilde {\mc O_S[z,w]/(zw-t_i)})/\mu _n]\simeq \mc U_i\times _{U_i}\Sp (\mc O_{C, \bar p_i}).
$$
Then $\mc L$ corresponds to the representation of $\mu _N$ with basis $z\dot w$.  By the assumptions on the $\mu _n$-action in (\ref{A1} (v)) it follows that $\mc L$ is trivial, and hence $o=0$.
\end{proof}

As in \cite[3.10]{logtwisted}, the log structure $\widetilde {\mc M}_{\mc C}$ is not the ``right'' log structure on $\mc C$ as it does not take into account the markings.   Exactly as in loc. cit., for each $i=1, \dots , n$ the ideal sheaf $\mc I_i\subset \mc O_{\mc C}$ defining $\Sigma _i$ defines a log structure $\mc N_i$ on $\mc C$.  Set
$$
\mc M_{\mc C}:= \widetilde {\mc M}_{\mc C}\oplus _{\mc O_{\mc C}^*}(\oplus _{i=1, \mc O_{\mc C}^*}^n\mc N_i).
$$
The map $f^*\mc M_S'\rightarrow \widetilde {\mc M}_{\mc C}$ then induces a log smooth morphism
$$
(\mc C, \mc M_{\mc C})\rightarrow (S, \mc M_S').
$$

If $g:(C, \sigma_1, \dots, \sigma _n)\rightarrow S$ is the coarse moduli space of $\mc C$ with its $n$ sections defined by the $p_i$, then the above construction applied to $(C, \sigma _1, \dots, \sigma_n)$ yields log structures $\mc M_C$ and $\mc M_S$ on $C$ and $S$ respectively and a special morphism
$$
(C, \mc M_C)\rightarrow (S, \mc M_S).
$$

\begin{prop}[{Generalization of \cite[4.7]{logtwisted}}] Let $\pi :\mc C\rightarrow C$ be the projection. There exists canonical morphisms of log structures $\pi ^b:\pi ^*\mc M_C\rightarrow \mc M_{\mc C}$ and $\ell :\mc M_S\hookrightarrow \mc M_S'$ such that the diagram of log stacks
$$
\begin{CD}
(\mc C, \mc M_{\mc C})@>(\pi , \pi ^b)>> (C, \mc M_C)\\
@VfVV @VVgV\\
(S, \mc M_S')@>(\text{id}, \ell)>> (S, \mc M_S)
\end{CD}
$$
commutes.  Moreover, the map $\ell:\mc M_S\hookrightarrow \mc M_S'$ is a simple extension.
\end{prop}
\begin{proof} This follows from the same argument proving \cite[4.7]{logtwisted}.  The key point is the local description (\ref{A1} (iv) and (v)) of a twisted  curve which enables one to rewrite \cite[4.6]{logtwisted} verbatim in the present situation.
\end{proof}
Finally note that for any geometric point $\bar s\rightarrow S$ the gerbe $\Sigma _{i, \bar s}$ is necessarily trivial and hence isomorphic to $\cB\mu _{a_i(\bar s)}$ for some integer $a_i(\bar s)$. One verifies immediately that the $a_i$ are positive integer valued locally constant functions on $S$.

The association
$$
(\mc C, \{\Sigma _i\})\mapsto (C/S, \{\sigma _i, a_i\}, \ell:\mc M_S\hookrightarrow \mc M_S')
$$
therefore defines a functor
\begin{equation}\label{inversefunctor}
 (\text{\rm twisted $n$-pointed  curves})
\rightarrow (\text{\rm $n$-pointed log twisted curves})
\end{equation}
As in \cite[4.8]{logtwisted} one sees that \ref{thefunctor} and \ref{inversefunctor} are inverse functors thereby proving \ref{Bthm}. \qed

\section{Remarks on $\text{Ext}$-groups  and base change\\ By Martin Olsson}\label{AppendixC}

In this appendix we gather together some fairly standard results about $\text{Ext}$-groups and base change, which will be used in the following appendix \ref{AppendixB}.  We include the results here as there does not seem to be an adequate treatment in the literature that covers the case when one does not have boundedness of cohomology (such as for Artin stacks).

\subsection{Base change for cohomology}

\begin{pg}\label{B:setup} Let $S$ be an integral noetherian scheme, and let $f:\mc X\rightarrow S$ be a proper algebraic stack over $S$.  Let $\mc F$ be a coherent sheaf on $\mc X$. For any morphism
$$
g:S'\rightarrow S
$$
we can then form the cartesian square
$$
\xymatrix{
\mc X'\ar[d]_-{f'}\ar[r]^-h& \mc X\ar[d]^-f\\
S'\ar[r]^-g& S,}
$$
and for every $i\in \mathbb{Z}$ we have a base change morphism
\begin{equation}\label{B:basechange}
g^*R^if_*\mc F\rightarrow R^if'_*h^*\mc F.
\end{equation}
\end{pg}

\begin{thm}\label{T:basechange}
With notation as in \ref{B:setup}, for every integer $i\in \mathbb{Z}$ there exists a dense open subset $U\subset S$ (which depends on $i$) such that for any morphism $g:S'\rightarrow S$ which factors through $U$, the base change morphism (\ref{B:basechange}) is an isomorphism.
\end{thm}
\begin{proof} First of all, by shrinking on $S$ we may assume that $S$ is affine, say $S = \Sp (R)$, and that $\mc F$ is flat over $S$.   Furthermore, it suffices to show that there exists a dense open subset $U\subset \Sp (R)$ such that for any morphism of affine schemes
$$
g:S' = \Sp (R')\rightarrow \Sp (R)
$$
which factors through $U$, the natural map
$$
H^i(\mc X, \mc F)\otimes _RR'\rightarrow H^i(\mc X', h^*\mc F)
$$
is an isomorphism.

Let $X_\cdot \rightarrow \mc X$ be a smooth hypercover of $\mc X$, with each $X_n$ a finite disjoint union of affine schemes.  For $n\geq 0$ set
$$
F_n:=\Gamma (X_n, \mc F|_{X_n}),
$$
so we obtain a complex (using the standard {\v C}ech differentials)
$$
C^\cdot : F^0\rightarrow F^1\rightarrow F^2\rightarrow \cdots 
$$
which computes the cohomology $R\Gamma (\mc X, F)$.

For any morphism of affine schemes
$$
g:S' = \Sp (R')\rightarrow \Sp (R)
$$
the base change of $X_\cdot $ to $S'$ is a smooth hypercover $X'_\cdot \rightarrow \mc X'$ which similarly can be used to compute the cohomology of $h^*\mc F$.  From this it follows that it suffices to show that for any $i\in \mathbb{Z}$ there exists a dense open subset $U\subset S$ such that if $g$ factors through $U$, then the natural map
$$
H^i(C^\cdot )\otimes _RR'\rightarrow H^i(C^\cdot \otimes _RR')
$$
is an isomorphism (note that here $C^\cdot \otimes _RR'$ is \emph{not} the derived tensor product).  This follows from the following lemma.
\end{proof}

\begin{lem} Let $R$ be a noetherian integral domain, and let $C^\cdot $ be a complex of flat $R$-modules such that for every $j\in \mathbb{Z}$ the $R$-module $H^j(C^\cdot )$ is finitely generated.  Then for every $i\in \mathbb{Z}$ there exists an open subset $U\subset \Sp (R)$ such that for every morphism $\Sp (R')\rightarrow U$ the natural map
$$
H^i(C^\cdot )\otimes _RR'\rightarrow H^i(C^\cdot \otimes _RR')
$$
is an isomorphism.
\end{lem}
\begin{proof}
Let $\delta _{\leq i+1}C^\cdot $ denote the complex whose $j$-th term is $C^j$ if $j\leq i+1$ and $0$ if $j>i+1$:
$$
\delta _{\leq i+1}C^\cdot : \cdots \rightarrow C^i\rightarrow C^{i+1}\rightarrow 0\rightarrow \cdots .
$$
There is a natural surjection
$$
C^\cdot \rightarrow \delta _{\leq i+1}C^\cdot
$$
which for every ring homomorphism $R\rightarrow R'$ (in particular the identity $R\rightarrow R$) induces an isomorphism
$$
H^j(C^\cdot \otimes _RR')\rightarrow H^j(\delta _{\leq i+1}C^\cdot \otimes _RR'), \ \ \text{for $j\leq i$.}
$$
Replacing $C$ by $\delta _{\leq i+1}C^\cdot $ we are then reduced to the case when $C^\cdot $ is bounded above, where the result is standard.
\end{proof}

\subsection{Base change for Ext}

Let $f:\mc X\rightarrow S$ be a finite type morphism between noetherian algebraic stacks, and assume that $S$ is an integral scheme.  Fix  $L^\cdot \in D^-_{\text{coh}}(\mc X)$ (the derived category of bounded above complexes of $\mc O_{\mc X}$-modules with coherent cohomology sheaves) and a coherent sheaf $\mls J\in \text{Coh}(\mc X)$.

\begin{lem}\label{L:B6} For every integer $n$, there exists a dense open subset (which depends on $n$) $U\subset S$ such that for any cartesian diagram
$$
\begin{CD}
\mc X'@>h>> \mc X\\
@Vf'VV @VVfV \\
S'@>g>> S,
\end{CD}
$$
where $g$ factors through $U$, the natural map
$$
h^*\mls Ext^p(L^\cdot, \mls J)\rightarrow \mls Ext^p(\LL h^*L^\cdot, \LL h^*\mls J)
$$
is an isomorphism for $p\leq n$.
\end{lem}
\begin{proof}
The assertion is local in the fppf topology on $\mc X$, $S$, and $S'$.  We may therefore assume that $\mc X = \Sp (R)$ and $S = \Sp (A)$ are affine schemes, and that $L^\cdot $ can be represented by a bounded above complex of projective $R$-modules of finite type, which we again denote by $L^\cdot $.  Let $J$ denote the $R$-module corresponding to the sheaf $\mls J$, and let $F^\cdot $ denote the bounded below complex of finite type $R$-modules
$$
F^\cdot := \text{Hom}^\cdot (L^\cdot, J).
$$
After shrinking on $S$ we may assume that $J$ is flat over $A$, in which case each $F^j$ is flat over $A$.
We need to show that after possibly replacing $S$ by a dense affine open subset, the natural map
$$
H^p(F^\cdot )\otimes _AA'\rightarrow H^p(F^\cdot \otimes _AA')
$$
is an isomorphism for all ring homomorphisms $A\rightarrow A'$ and all $p\leq n$.  This is a standard argument and we leave it to the reader (see for example \cite[IV.9.4.3]{EGA} where a similar argument is made).
\end{proof}

\begin{thm}\label{cohboundedness} Let $f:\mc X\rightarrow S$ be a  finite type proper morphism of noetherian algebraic stacks with $S = \Sp (A)$ an integral affine scheme. Let $L^\cdot \in D^-_{\text{\rm coh}}(\mc X)$ and $\mls J\in \text{\rm Coh}(\mc X)$. Then 
for every integer $n$, there exists a dense open subset (which depends on $n$) $U\subset S$ such that for any cartesian diagram
$$
\begin{CD}
\mc X'@>h>> \mc X\\
@Vf'VV @VVfV \\
S'=\Sp (A')@>g>> S,
\end{CD}
$$
where $g$ factors through $U$, the natural map
\begin{equation}\label{targetmorphism}
\text{\rm Ext}^p(L^\cdot, \mls J)\otimes _AA'\rightarrow \text{\rm Ext}^p(\LL h^*L^\cdot, \LL h^*\mls J)
\end{equation}
is an isomorphism for $p\leq n$.
\end{thm}
\begin{proof} 
First after shrinking on $S$ we may assume that $\mls J$ is flat over $S$ and  that for any morphism $g:\Sp (A')\rightarrow \Sp (A)$ the pullback map
$$
h^*\mls Ext^p(L^\cdot , \mls  J)\rightarrow \mls Ext^p(\LL h^*L^\cdot, \LL h^*\mls J)
$$
is an isomorphism for all $p\leq n$ (by \ref{L:B6}).  By \ref{T:basechange} we may,  after shrinking some more on $S$, assume that each of the terms $E_r^{pq}$ for $p+q\leq n$ in the spectral sequence
$$
E_2^{pq} = H^p(\mc X, \mls Ext^q(L^\cdot, \mls J))\implies \text{Ext}^{p+q}(L^\cdot, \mls J)
$$
are flat over $S$, and that their formation commute with arbitrary base change $\Sp (A')\rightarrow \Sp (A)$.  That \ref{targetmorphism} is an isomorphism then follows from consideration of the morphism of spectral sequences
$$
\begin{CD}
E_2^{pq} = H^p(\mc X, \mls Ext^q(L^\cdot, \mls J))\implies \text{Ext}^{p+q}(L^\cdot, \mls J)\\
@VVV \\
E_2^{pq} = H^p(\mc X', \mls Ext^q(\LL h^*L^\cdot, \LL h^*\mls J))\implies \text{Ext}^{p+q}(\LL h^*L^\cdot, \LL h^*\mls J).
\end{CD}
$$
\end{proof}

\section{Another boundedness theorem for Hom-stacks\\ By Martin Olsson}\label{AppendixB}

\subsection{Statement of Theorems}

\begin{thm}\label{athm} Let $\mc B$ be  an Artin stack, and let $\mc X$ and $\mc Y$ be Artin stacks of finite presentation over $\mc B$ with finite diagonals.  Assume that $\mc X$ is flat and proper over $\mc B$.  Then:

(i) Let  $\underline {\text{\rm Hom}}(\mc X, \mc Y)$ denote the stack over $\mc B$ whose fiber over a $\mc B$-scheme $T$ is the groupoid of $T$-morphisms of stacks
$$
\mc X\times _{\mc B}T\rightarrow \mc Y\times _{\mc B}T.
$$
Then $\underline {\text{\rm Hom}}(\mc X, \mc Y)$ is an algebraic stack locally of finite presentation over $\mc B$, with quasi-compact and separated diagonal over $\mc B$.

(ii) Assume further that $\mc Y\rightarrow \mc B$ is a tame morphism in the sense of \ref{D:tamedef}.  If $Y$ denotes the relative coarse moduli space of $\mc Y/\mc B$ (see \ref{Th:rel-cms}), then the natural map
$$
\underline {\text{\rm Hom}}(\mc X, \mc Y)\rightarrow \underline {\text{\rm Hom}}(\mc X, Y)
$$
is of finite type.
\end{thm}

\begin{rem} If $\mc X$ is a tame stack and $\mc B = B$ is an algebraic space, then the condition that $\mc X/B$ is flat implies that the coarse moduli space $X/B$ of $\mc X$ is also flat by \cite[3.3 (b)]{AOV} and its formation commutes with arbitrary base change by \cite[3.3]{AOV}. It follows that $\Hom (X, Y)$ is an algebraic space locally of finite presentation over $B$, and by the universal property of coarse moduli space there is a canonical isomorphism
$$
\Hom (\mc X, Y)\simeq \Hom (X, Y).
$$
\end{rem}

In the context of \ref{athm}  one can also consider the substack 
$$
\Hom ^{\text{rep}}(\mls X, \mls Y)\subset \Hom (\mls X, \mls Y)
$$
classifying representable morphisms $\mls X\rightarrow \mls Y$.  As in \cite[1.6]{homstack} this substack is an open substack, and therefore the following corollary follows from (\ref{athm} (ii)):

\begin{cor}\label{repfinite} With assumptions as in \ref{athm}, the natural map 
$$
\Hom ^{\text{\rm rep}}(\mls X, \mls Y)\rightarrow \Hom (\mls X, Y)
$$
is of finite type. 
\end{cor}

The rest  of this appendix is devoted to the proof of \ref{athm}.  

First of all observe that it suffices to consider the case when $\mc B$ is an algebraic space (which we will henceforth denote by $B$), by the following well-known result:

\begin{lemma}\label{Lem:relative-algebraic}
Let $\cX \to \cY$ be a morphism of stacks, and assume the following.
\begin{enumerate}
\item The stack $\cY$ is algebraic.
\item There is an algebraic space $Z$ and a smooth surjective morphism $Z \to \cY$ such that $\cX_Z:=Z\times_\cY \cX$ is algebraic.
\end{enumerate}
Then $\cX$ is algebraic.
\end{lemma}

\begin{proof}
First we note that the relative diagonal $\cX \to \cX \times_\cY \cX$ is representable. Indeed, let $U$ be an algebraic space and $U \to \cX \times_\cY \cX$ a morphism. Form the pullback $U_Z = U \times_\cY Z$. We have a morphism $U_Z \to \cX_Z \times_Z \cX_Z$, and by assumption $\cX_Z\times_{\cX_Z \times_Z \cX_Z} U_Z$ is representable. Applying the same to $Z\times_{\cY} Z$ we obtain flat descent data for $\cX \times_{\cX \times_\cY \cX} U$, showing that the latter is representable. Therefore the 
relative diagonal $\cX \to \cX \times_\cY \cX$ is representable.

Now $\cX \times_\cY \cX \to \cX \times \cX$ is the pullback of the diagonal
$\cY \to \cY \times\cY$ via $ \cX \times \cX \to \cY \times\cY$. Since the diagonal
$\cY \to \cY \times\cY$ is representable we have that $\cX \times_\cY \cX \to \cX \times \cX$ is representable. 

Composing we get that the diagonal $\cX  \to \cX \times \cX$ is representable. It remains to produce a smooth presentation for $\cX$.
 
There exists a scheme $V$ with a smooth surjective  $V\to \cX_Z$. Since $\cX_Z \to \cX$ is smooth surjective, we have $V \to \cX$ smooth surjective.
\end{proof}

Furthermore, by a standard limit argument it suffices to consider the case when $B$ is of finite type over an excellent Dedekind ring, which we assume henceforth.

\begin{pg}\emph{Algebraicity of $\Hom (\mc X, \mc Y)$.}\label{C.30}
That $\Hom (\mc X, \mc Y)$ is an algebraic stack locally of finite presentation over $B$ is shown in \cite[1.1]{Aoki} (see also the erratum \cite{Aokierratum}).

For a morphism of stacks $\mc G\rightarrow \mc X$ over $B$, let 
$$
\Sec (\mc G/\mc X)
$$
denote the stack over $B$ whose fiber over a $B$-scheme $T$ is the groupoid of sections of the base change $\mc G_T\rightarrow \mc X_T$.  The algebraicity of the $\Hom $-stacks implies that if $\mc X/B$ is a proper flat algebraic stack with finite diagonal, and if $\mc G\rightarrow \mc X$ is a morphism of algebraic stacks with $\mc G$ a finitely presented Artin stack over $B$ with finite diagonal, then the stack $\Sec (\mc G/\mc X)$ is an algebraic stack locally of finite presentation over $B$, as it
is equal to the fiber product of the diagram
$$
\xymatrix{
& \Hom (\mc X, \mc G)\ar[d]\\
S\ar[r]^-{\text{id}_{\mc X}}& \Hom (\mc X, \mc X).}
$$
\end{pg}

\begin{pg}\label{P:C6}\emph{The diagonal of $\Hom (\mc X,  Y)$ is quasi-compact and separated.}
Under  the assumptions of (\ref{athm} (i)), we now show that the diagonal of $\Hom (\mc X, \mc Y)$ is quasi-compact and separated in the case when $\mc Y = Y$ is an algebraic space (the general case will be treated in \ref{P:C13} below).

To see that the diagonal of $\Hom (\mc X, Y)$ is separated, let 
$$
f_1, f_2:\mc X\rightarrow Y
$$
be two morphisms, and let  $I$ denote the fiber product of the diagram
$$
\xymatrix{
& B\ar[d]^-{f_1\times f_2}\\
\Hom (\mc X, Y)\ar[r]^-\Delta & \Hom (\mc X, Y)\times \Hom (\mc X, Y).}
$$
The space $I$ 
represents the functor on $B$-schemes which to any $T/B$ associates the unital set if the base changes of the maps $f_1$ and $f_2$ to $T$ are equal, and the empty set otherwise.  In particular, the morphism $I\rightarrow B$ is a monomorphism, which implies that the diagonal
$$
I\rightarrow I\times _BI
$$
is an isomorphism.  Therefore $I$ is separated.  

To see that the diagonal of $\Hom (\mc X, Y)$ is quasi-compact, we may without loss of generality assume that $B$ is reduced, and in this case it suffices by noetherian induction to exhibit a dense open subset of $B$ over which $I$ is quasi-compact.  By shrinking on $B$ we may therefore assume that the coarse moduli space $X$ of $\mc X$ is flat over $B$, and that its formation commutes with arbitrary base change on $B$.  In this case we have an isomorphism
$$
\Hom (\mc X, Y)\simeq \Hom (X, Y),
$$
and the quasi-compactness of the diagonal follows from \cite[5.1]{bounded}.
\end{pg}

Next we use the results of appendix \ref{AppendixC} to generalize three basic results about modifications in \cite[\S 4]{bounded}, which will be necessary to proceed further in the proof of \ref{athm}.

\

\noindent {\bf Passage to the maximal reduced substack.}

\begin{pg} Let $S$ be a noetherian scheme, and let $\mc G\rightarrow \mc X$ be a finite type morphism between Artin stacks of finite type over $S$ with finite diagonals.  Assume further that $\mc X$ is flat and proper over $S$.  
\end{pg}

\begin{prop}\label{finitelem} Let $\mc X_0\hookrightarrow \mc X$ be a closed immersion
defined by a nilpotent ideal $\mc J\subset \mls O_{\mc X}$, and assume $\mc X_0$ is flat over $S$. Let $\mls G_0$ denote the base change $\mls G\times _{\mc X}\mc X_0$. Then the natural map
$$
\underline {\text{\rm Sec}}(\mls G/\mc X)\rightarrow \underline
{\text{\rm Sec}}(\mls G_0/\mc X_0)
$$
is of finite type with quasi-compact diagonal.
\end{prop}
\begin{proof}
This is essentially the same as in \cite[5.11]{homstack}.

By noetherian induction, it suffices to show that the morphism is of finite type with quasi-compact diagonal over a dense open subset of $S$.  We may further assume that $S$ is reduced (since if $U$ and $V$ are $S$-schemes locally of finite type and $g:U\rightarrow V$ is a morphism, then $g$ is of finite type if and only if the base change of $g$ to $S_{\text{red}}$ is of finite type), and using the same argument given in \cite[paragraph following proof of 5.11]{homstack} that $\mc J^2 = 0$.

It suffices to show that if $T\rightarrow \Sec (\mls G_0/\mc X_0)$ is a
morphism corresponding to a section $s:\mc X_{0, T}\rightarrow \mls G_{0, T}$,
then the fiber product
$$
\mc P:= \Sec (\mls G/\mc X)\times _{\Sec (\mls G_0/\mc X_0)}T
$$
is of finite type over $T$ with quasi-compact diagonal.  Furthermore, we may assume that $T$
is an integral noetherian affine scheme.  The stack $\mc P$ associates to any
$w:W\rightarrow T$ the groupoid of liftings $\widetilde s:\mc X_{W
}\rightarrow \mls G$ over $\mc X$ of the composite
$$
\mc X_{0, W}\rightarrow \mls G_0\rightarrow \mls G,
$$
where the first map is the one induced by $s$.  

To prove that $\mc P$ is quasi--compact with quasi-compact diagonal, it suffices
by Noetherian induction to exhibit a dense open set $U\subset T$
such that $\mc P_U$ is quasi--compact with quasi-compact diagonal.

Let $L_{\mls G_0/\mc X_0}$ be the cotangent complex of $\mls G_0/\mc X_0$. By \ref{cohboundedness},  after replacing $T$ by a dense open subscheme we may assume that the groups
$$
\text{Ext}^i(s^*L_{\mls G_0/\mc X_0}, \mc J), \ \ i=-1, 0, 1
$$
are projective modules on $T$ of finite type, and that the formation of these modules commutes with arbitrary base change on $T$. 

By \cite[1.5]{deform} there is a
canonical obstruction
$$
o\in \text{Ext}^1(s^*L_{\mc G_0/\mc X_0}, \mc J)
$$
whose vanishing is necessary and sufficient for the existence of a
lifting $\widetilde s$ as above, and the formation of this obstruction is functorial in $T$.   After replacing $T$ by the closed subscheme defined by the condition that $o$ vanishes, we may assume that $o=0$.  In this case the set of isomorphism classes of liftings $\widetilde s$ form a torsor under
$$
\text{Ext}^0(s^*L_{\mc G_0/\mc X_0}, \mc J)
$$
and the group of infinitesimal automorphisms of $\widetilde s$ is canonically isomorphic to
$$
\text{Ext}^{-1}(s^*L_{\mc G_0/\mc X_0}, \mc J).
$$
It follows that $\mc P$ is an $\text{Ext}^{-1}(s^*L_{\mc G_0/\mc X_0}, \mc J)$--gerbe over $$\text{Ext}^0(s^*L_{\mc G_0/\mc X_0}, \mc J),$$ and in particular is quasi--compact with quasi-compact diagonal.
\end{proof}

\noindent {\bf The case of a finite morphism.}

\begin{pg} Let $S$ be a noetherian scheme, and $\mc X/S$ a proper flat Artin stack with finite diagonal.  Let $\mc G\rightarrow \mc X$ be a finite morphism.
\end{pg}

\begin{prop}\label{finitecase} The stack $\Sec (\mc G/\mc X)$ is an algebraic space separated and of finite type over $S$.
\end{prop}
\begin{proof}
Let us first verify that the diagonal of $\Sec (\mc G/\mc X)$ is separated and quasi-compact.  Equivalently, let
$$
s_1, s_2:\mc X\rightarrow \mc G
$$
be two sections and let $I$ denote the fiber product of the diagram
$$
\xymatrix{
& S\ar[d]^-{s_1\times s_2}\\
\Sec (\mc G/\mc X)\ar[r]^-\Delta & \Sec (\mc G/\mc X)\times \Sec (\mc G/\mc X).}
$$
We need to show that $I$ is quasi-compact and separated.

That $I$ is separated is immediate, for as in \ref{P:C6} the fact that the map $\mc G\rightarrow \mc X$ is representable implies that $I$ is a subfunctor of $B$, and so the diagonal of $I/B$ is an isomorphism.

To see the quasi-compactness of $I$, let $\mc Z\hookrightarrow \mc X$ be the closed substack  defined as the fiber product of the diagram (note that since $\mc G/\mc X$ is separated the diagonal of $\mc G/\mc X$ is a closed immersion)
$$
\xymatrix{
& \mc G\ar[d]^-\Delta \\
\mc X\ar[r]^-{s_1\times s_2}& \mc G\times _{\mc X}\mc G.}
$$
Then $I$  represents the functor which to any $B$-scheme $T$ associates the unital set if $\mc Z_T= \mc X_T$, and the empty set otherwise.  

Now to prove that $I$ is quasi-compact, we can without loss of generality replace $S$ by $S_{\text{red}}$ and hence may assume that $S$ is reduced, and by noetherian induction it suffices to show that the restriction of $I$ to some dense open subset of $S$ is quasi-compact.  We may therefore assume that $S$ is integral, and that the stack $\mc Z$ is flat over $S$.  In this case $I$ is isomorphic to an open subset of $S$.  Indeed let $K$ denote the kernel of the surjection
$$
\mc O_{\mc X}\rightarrow \mc O_{\mc Z}.
$$
Since $\mc Z$ is flat over $S$, for any $S$-scheme $T$ the sequence
$$
0\rightarrow K\rightarrow \mc O_X\rightarrow \mc O_{\mc Z}\rightarrow 0
$$
remains exact when pulled back to $\mc X_T$.  It follows that if $\mc R\subset \mc X$ denotes the support of the coherent sheaf $K$, then $I$ is represented by the complement of the closed image (since $\mc X/S$ is proper) of $\mc R$ in $S$.  This completes the verification that the diagonal of $\Sec (\mc G/\mc X)$ is quasi-compact and separated.

To verify that $\Sec (\mc G/\mc X)$ is of finite type over $S$, it suffices to show that its pullback to $S_{\text{red}}$ is of finite type.  We may therefore assume that $S$ is reduced.  Furthermore, we may by \ref{finitelem} (applied with $\mc X_0 = \mc X_{\text{red}}$) assume that $\mc X$ is also reduced.  Furthermore, by noetherian induction it suffices to exhibit a dense open subset of $S$ over which $\Sec (\mc G/\mc X)$ is of finite type.

By Chow's lemma for Artin stacks \cite[1.1]{chow}, there exists a proper surjection $h:X'\rightarrow \mc X$ with $X'$ a projective $S$-scheme.  After shrinking on $S$ we may also assume that $X'$ is flat over $S$ (since $S$ is reduced).  Since $\mc X$ is reduced the map $\mc O_{\mc X}\rightarrow h_*\mc O_{X'}$ is injective.  After shrinking on $S$, we may assume that the formation of $h_*\mc O_{X'}$ commutes with arbitrary base change $S'\rightarrow S$, and that for any such base change the map $\mc O_{\mc X_{S'}}\rightarrow h_*\mc O_{X_{S'}'}$ is also injective.  Let $G'\rightarrow X'$ denote the pullback $\mc G\times _{\mc X}X'$.

By \cite[5.10]{homstack} the stack (in fact an algebraic space) $\Sec ( G'/X')$ is of finite type over $S$.  It therefore suffices to show that the pullback map
\begin{equation}\label{5.9.1}
\Sec (\mc G/\mc X)\rightarrow \Sec ( G'/X')
\end{equation}
is of finite type.

Let $\mc A$ be the coherent sheaf of algebras defined by $\mc G = \underline {\text{\Sp }}_{\mc X}(\mc A)$, and let $\mc B$ denote the sheaf of $\mc O_{\mc X}$-algebras $h_*\mc O_{X'}$.  For a morphism $T\rightarrow S$, let $\mc A_T$ (resp. $\mc B_T$) denote the pullback of $\mc A$ (resp. $\mc B$) to $\mc X_T:= \mc X\times _ST$.  Then $\Sec (G'/X')$ is equal to the functor which to any scheme $T\rightarrow S$ associates the set of morphisms of $\mc O_{\mc X_T}$-algebras $\varphi :\mc A_T\rightarrow \mc B_T$ over $\mc X_T$.  The stack $\Sec (\mc G/\mc X)$ is the subfunctor of morphisms $\varphi $ which factor through $\mc O_{\mc X_T}\subset \mc B_T$.

Fix a morphism $\varphi :\mc A_T\rightarrow \mc B_T$ over a noetherian scheme $T$ defining a $T$-valued point of $\Sec (G'/X')$, and set
$$
P:= \Sec (\mc G/\mc X)\times _{\Sec (G'/X'), \varphi }T.
$$
To prove that $P$ is of finite type, we may again assume that $T$ is reduced, and by noetherian induction it suffices to show that $P\rightarrow T$ is of finite type over a dense open subset of $T$.  Let $M$ denote the cokernel of the map $\mc O_{\mc X_T}\hookrightarrow \mc B_T$, and let $\bar \varphi :\mc A_T\rightarrow M$ be the map induced by $\varphi $.    Let $Q$ be the cokernel of $\bar \varphi $, and let $K$ be the kernel of the map $M\rightarrow Q$ so that there is an exact sequence
\begin{equation}\label{Qseq}
0\rightarrow K\rightarrow M\rightarrow Q\rightarrow 0.
\end{equation}
After shrinking on $T$ we may assume that $Q$ is flat over $T$ in which case \ref{Qseq} remains exact after arbitrary base change $T'\rightarrow T$.  In this case, let $\mc Z\subset \mc X_T$ be the support of the coherent sheaf $K$, and let $W\subset T$ be the image of $\mc Z$ (which is closed since $\mc X_T\rightarrow T$ is proper).  Then $P$ is represented by the complement of $W$ in $T$.

To verify that $\Sec (\mc G/\mc X)$ is separated, note that since we already know that the diagonal is quasi-compact and separated, it suffices to verify the valuative criterion for properness.  This amounts to the following.
  Assume that $S$ is the spectrum of a valuation ring and let $\eta \in S$ be the generic point.  Assume given two sections $s_1, s_2:\mc X\rightarrow \mc G$ whose restrictions to $\eta $ are equal.  Then we need to show that $s_1 = s_2$.  For this note that since $\mc G$ is finite over $\mc X$, we have $\mc G = \underline {\text{Spec}}_{\mc X}(\mc A)$ for some coherent sheaf of $\mc O_{\mc X}$-algebras $\mc A$ on $\mc X$, and the sections $s_1$ and $s_2$ are specified by two morphisms of $\mc O_{\mc X}$-algebras $\rho _1, \rho _2:\mc A\rightarrow \mc O_{\mc X}$.  Let $j:\mc X_{\eta }\hookrightarrow \mc X$ be the inclusion of the generic fiber.  Since $\mc X$ is flat over $S$, the natural map $\mc O_{\mc X}\rightarrow j_*\mc O_{\mc X_{\eta }}$ is an inclusion.  Therefore it suffices to show that the composite maps
$$
\xymatrix{
\mc A\ar[r]^-{\rho _i}& \mc O_{\mc X}\ar[r]& j_*\mc O_{\mc X_\eta }}
$$
are equal.  Equivalently that the maps on the generic fiber
$$
\rho _{i, \eta }:\mc A_\eta \rightarrow \mc O_{\mc X_\eta }
$$
are equal, which holds by assumption.
\end{proof}

\noindent {\bf Behavior with respect to proper modifications of $\mc X$.}

\begin{pg}\label{P:C.14} Let $S$ be a noetherian scheme, and let $\mc X$ and $\mc Y$ be Artin stacks of finite type over $S$ with finite diagonals.   Assume that the following conditions hold:
\begin{enumerate}
\item [(i)] The formation of the coarse spaces $\pi _{\mc X}:\mc X\rightarrow X$ and $\pi _{\mc Y}:\mc Y\rightarrow Y$ commutes with arbitrary base change $S'\rightarrow S$.  
\item [(ii)] $\mc X$ is proper and flat over $S$.

\item [(iii)] The coarse space $X$ of $\mc X$ is flat over $S$ (note that $X$ is automatically proper over $S$ since $\mc X$ is proper over $S$).
\end{enumerate}
\end{pg}

\begin{prop}\label{modprop} Let $m:\mc X'\rightarrow \mc X$ be a proper representable surjection with $\mc X'$ a proper and flat algebraic stack over $S$ with finite diagonal.  Let $f:\mc X\rightarrow Y$ be a morphism, and let $\mc G$ (resp. $\mc G'$) denote the pullback along $f$ (resp. $f\circ m$) of $\mc Y$.  Assume that $\Sec (\mc G/\mc X)$ and $\Sec (\mc G'/\mc X')$ are algebraic stacks locally of finite type over $S$ with quasi-compact diagonals.  Then the pullback map
\begin{equation}\label{pullmap}
\Sec (\mc G/\mc X)\rightarrow \Sec (\mc G'/\mc X')
\end{equation}
is of finite type.
\end{prop}
\begin{proof}  
This follows from the same argument proving \cite[4.13]{bounded}. 
\end{proof}

\subsection{Proof of (\ref{athm} (i))}\label{P:C13}

\

Let $\mc X$ and $\mc Y$ be as in (\ref{athm} (i)).
To complete the proof of (\ref{athm} (i)), it remains to show that the diagonal of $\Hom (\mc X, \mc Y)$ is separated and of finite type.

  Since
$$
\Hom (\mc X, Y)
$$
is locally of finite type with quasi-compact and separated diagonal (by \ref{P:C6}), to prove that $\Hom (\mc X, \mc Y)$ has separated and quasi-compact diagonal over $B$, it suffices to show that $\Hom (\mc X, \mc Y)$ has separated and quasi-compact diagonal over $\Hom (\mc X, Y)$.  Equivalently, let $f:\mc X\rightarrow Y$ be a fixed morphism, and set
$$
\mc G:= \mc X\times _{f, Y}\mc Y.
$$
Then it suffices to show that $\Sec (\mc G/\mc X)$ has separated and quasi-compact diagonal.

Let  $s_1, s_2:\mls X\rightarrow \mls G$ be two sections, and let $I$ denote the fiber product of the diagram (which is a finite stack over $\mls X$ since $\mls G$ has finite diagonal)
$$
\begin{CD}
@. \mls X\\
@. @VVs_1\times s_2V \\
\mls G@>\Delta >> \mls G\times _{\mls X}\mls G.
\end{CD}
$$
Then the fiber product of the diagram
$$
\xymatrix{
& B\ar[d]^-{s_1\times s_2}\\
\Sec (\mc G/\mc X)\ar[r]^-{\Delta }& \Sec (\mc G/\mc X)\times _B\Sec (\mc G/\mc X)}
$$
is isomorphic to
 $\Sec (I/\mls X)$ which is quasi-compact and separated by \ref{finitecase}.
 
 This completes the proof of (\ref{athm} (i)). \qed 

\subsection{Reduction of (\ref{athm} (ii)) to a special case}

\begin{pg} Turning to the proof of (\ref{athm} (ii)), note first that it suffices to show that 
if  $f:\mc X\rightarrow Y$ is a morphism and $\mc G$ denotes the fiber product
$$
\mc G:= \mc X\times _{f, Y}\mc Y,
$$
then the stack $\Sec (\mc G/\mc X)$ is of finite type over $B$.  This we do by a series of reductions.
\end{pg}
\begin{pg}
To prove that  $\Sec (\mls G/\mls X)$  is quasi-compact, it suffices by noetherian induction to exhibit a dense open subset of $B$ where this is so.  Furthermore, we may without loss of generality assume that $B$ is integral.  By Chow's lemma for Artin stacks \cite[1.1]{chow} there exists a proper surjection $X'\rightarrow \mls X$ with $X'$ an algebraic space.  After shrinking on $B$ we may assume that $X'$ is flat over $B$.
Applying \ref{modprop} to the stack $\Sec (\mls G/\mls X)$, we are then reduced to proving that $\Sec (\mls G'/X')$ is of finite type over $B$.  This reduces to the case when $\mc X$ is an algebraic space, which we assume for the rest of the proof (and for consistency we write $X$ instead of $\mc X$).
\end{pg}

\begin{pg}
Furthermore, we can without loss of generality assume that $B$ is integral.   By noetherian induction, it suffices to exhibit a dominant morphism $B'\rightarrow B$ such that the restriction of $\Sec (\mls G/X)$ to $B'$ is of finite type.
\end{pg}

\begin{pg} By the same argument used in \cite[6.4]{bounded} we then reduce the proof of (\ref{athm} (ii)) to the case when $X$ is integral  and $\mls G\rightarrow X$ has a section $s:X\rightarrow \mls G$.  Let $\underline {\text{Aut}}(s)$ denote the group scheme of automorphisms of $s$ (a scheme over $X$).  Since $X$ is reduced there exists a dense open subset $U\subset X$ such that the restriction of $\underline {\text{Aut}}(s)$ to $U$ is flat over $U$.  Let $H\subset \underline {\text{Aut}}(s)$ denote the scheme--theoretic closure of $\underline {\text{Aut}}(s)|_U$.  By \cite[Premi\'ere partie 5.2.2]{GR} there exists a blow-up $X'\rightarrow X$ with center a proper closed subspace of $X$ such that the strict transform of $H$ in $\underline {\text{Aut}}(s)$ is flat over $X'$. After further shrinking on $S$ we may assume that $X'$ is flat over $S$.  Since the map
$$
\Sec (\mls G/X)\rightarrow \Sec (\mls G\times _XX'/X')
$$
is of finite type by \ref{modprop}, we may therefore replace $X$ by $X'$ and hence can assume that $H$ is a finite flat subgroup scheme of $\underline {\text{Aut}}(s)$.  Since all the geometric fibers of $H$ are linearly reductive being closed subgroup schemes of linearly reductive groups (see \cite[2.7]{AOV}), the $X$-group scheme $H$ is in fact linearly reductive. If $\widetilde X$ is the normalization of $X$ the map
$$
\Sec (\mls G/X)\rightarrow \Sec (\mls G\times _X\widetilde X/\widetilde X)
$$
is also of finite type by \ref{modprop}.  This enables us to reduce to the following situation: $X$ is normal, $s:X\rightarrow \mls G$ is a section, and $H\subset \underline {\text{Aut}}(s)$ is a finite closed  subgroup scheme which is flat and linearly reductive over $X$ and over a dense open subset of $X$ the group scheme $H$ is equal to $\underline {\text{Aut}}(s)$.
\end{pg}

\begin{lem} The morphism $\pi :\cB H\rightarrow \mls G$ induced by $s$ identifies $\cB H$ with the normalization of $\mls G_{\text{\rm red}}$.
\end{lem}
\begin{proof}
Étale locally on $X$ we can write $\mls G = [W/GL_n]$ for some $n$, where $W$ is an $X$--scheme.  Indeed by \cite[3.2]{AOV} we can \'etale locally on $X$ write $\mls G = [Z/G]$ with $G$ a linearly reductive group scheme.  Choosing any embedding $G\hookrightarrow GL_n$ we take $W$ to be the quotient of $Z\times GL_n$ by the diagonal action of $G$.

Let $P\rightarrow X$ be the pullback $X\times _{s, \mls G}W$.  Then $P$ is a $GL_n$--torsor over $X$ with  a $GL_n$--equivariant morphism $f:P\rightarrow W$.  Since $P$ is smooth over $X$ the space $P$ is in particular normal, and since $s$ is proper and quasi--finite the morphism $f$ is finite. After replacing $X$ by an \'etale cover we may assume given a section $p\in P$.  Then $\underline {\text{Aut}}(s)$ is the closed subgroup scheme of $GL_{n, X}$ given by the fiber product of the diagram
$$
\begin{CD}
W\times GL_n@. \\
@V\rho VV @. \\
W\times W@<\Delta \circ f\circ p<< X,
\end{CD}
$$
where $\Delta :W\rightarrow W\times W$ is the diagonal and $\rho $ is the map sending $(v, g)$ to $(v, g(v))$.  In particular, we obtain an embedding $H\subset GL_{n, X}$.  Let $\overline P$ denote the quotient $P/H$.  The space $\overline P$ is normal.  By the construction the map $f$ induces a morphism $\bar f:\overline P\rightarrow W^*$, where $W
^*$ denotes the normalization of $W_{\text{red}}$.  This map $\bar f$ is finite, surjective, and birational and hence by Zariski's Main Theorem an isomorphism.  On the other hand, we have $\overline P\simeq W\times _{\mls G}\cB H$, and hence if $\mls G^*$ denotes the normalization of $\mls G_{\text{red}}$ we find that the base change of $\cB H\rightarrow \mls G^*$ to $W\times _{\mls G}\mls G^*$ is an isomorphism.  It follows that $\cB H\rightarrow \mls G^*$ is an isomorphism. 
\end{proof}

\begin{pg}
Let $\mls G^*$ denote the normalization of $\mls G$, so that $\mls G^*$ is a gerbe over $X$.
As in \cite[3.14]{bounded}, after shrinking on $B$ we may assume that for every field valued point $b\in B(k)$ the fiber $\mls G_b^*\rightarrow \mls G_b$ is the normalization of $\mls G_{b, \text{red}}$ and that $X_b$ is normal.  It then follows that the map
$$
\Sec (\mls G^*/X)\rightarrow \Sec (\mls G/X)
$$
is surjective on field valued points and hence surjective.  We may therefore replace $\mls G$ by $\mls G^*$ and therefore may assume that $\mls G$ is a gerbe over $X$.

Next we consider the case of a gerbe.  If the generic point of $B$ has characteristic $0$, then after shrinking on $B$ we can assume $\mls G$ is Deligne--Mumford in which case the result follows from \cite[1.1]{bounded}.  We may therefore assume that the generic point of $B$ has characteristic $p>0$, and hence after shrinking on $B$ may assume that $B$ is an $\mathbb{F}_p$--scheme.  In this case, for any $B$--scheme $T$ and $t\in \mls G(T)$, the automorphism group scheme $G_t$ of $t$ is canonically an extension
$$
1\rightarrow \Delta _t\rightarrow G_t\rightarrow H_t\rightarrow 1,
$$
where $\Delta _t$ is locally diagonalizable and $H_t$ is \'etale over $T$.  Indeed, the group scheme $G_t$ is tame so we can take $\Delta _t$ to be the subfunctor of $G_t$ classifying elements of $G_t$ killed by some power of $p$.  This normal subgroup $\Delta _t$ is functorial in the pair $(T, t)$.

Define $\mls H$ to be the stack (with respect to the fppf topology) associated to the prestack whose objects are the same as the objects of $\mls G$ but for which a morphism between two $t, t'\in \mls G(T)$ is defined to be the set
$$
\Delta _t\backslash \underline {\text{Hom}} (t, t') = \Delta _t\backslash \underline {\text{Hom}} (t, t')/\Delta _{t'} = \underline {\text{Hom}} (t, t')/\Delta _{t'}.
$$
Then $\mls H$ is a Deligne--Mumford stack and there is a natural map $\mls G\rightarrow \mls H$ over $X$.  By \cite[1.1]{bounded} the stack $\Sec (\mls H/X)$ is of finite type.  On the other hand, for any section $s:X\rightarrow \mls H$, the fiber product
$$
\Sec (\mls G/X)\times _{\Sec (\mls H/X), s}B
$$
is isomorphic to $\Sec (\mls G\times _{\mls H, s}X/B)$.  Now the stack $\mls G\times _{\mls H, s}X$ is a gerbe over $X$ whose stabilizer groups are all diagonalizable. This further reduces the proof to the case when the $H_t$'s in the above discussion are all trivial.
\end{pg}

\begin{pg}\label{B.33} In this case, there exists a canonical locally diagonalizable group scheme $\Delta $ on $X$ such that $\mls G$ is bound by $\Delta $. Indeed fppf locally on $X$ there exists a section $s:X\rightarrow \mls G$ giving rise to a group scheme $\Delta _s$.  If $s':X\rightarrow \mls G$ is a second section then locally $s$ and $s'$ are isomorphic so we obtain an isomorphism $\Delta _s\simeq \Delta _{s'}$.  This isomorphism is independent of the choice of the isomorphism between $s$ and $s'$ since $\Delta _s$ and $\Delta _{s'}$ are abelian.  It follows that the $\Delta _s$'s descend to a group scheme $\Delta $ on $X$.

Since the Cartier dual of $\Delta $ is a locally constant sheaf of finite abelian groups on $X_{\text{et}}$, there exists  a finite \'etale covering $X'\rightarrow X$ such that the pullback of $\Delta $ to $X'$ is a diagonalizable group scheme.  Since 
$$
\Sec (\mls G/X)\rightarrow \Sec (\mls G\times _XX'/X')
$$
is of finite type by \ref{modprop}, this reduces the proof to the case when $\mls G$ is a gerbe over $X$ bound by a diagonalizable group scheme $\Delta $.  In other words, when $\mls G$ corresponds to a cohomology class in $H^2(X, \Delta )$.  Write $\Delta = \mu _{n_1}\times \cdots \times \mu _{n_r}$ so that
$$
H^2(X, \Delta ) = \prod _iH^2(X, \mu _{n_i}).
$$
Using the resulting decomposition of the cohomology class of $\mls G$, we see that $\mls G$ is isomorphic to a product $\mls G_1\times _X\cdots \times _X\mls G_r$, where $\mls G_i$ is a $\mu _{n_i}$--gerbe over $X$.  Then 
$$
\Sec (\mls G/X)\simeq \prod _i\Sec (\mls G_i/X).
$$
This therefore finally reduces the proof to the case of a $\mu _n$--gerbe over $X$ which is the case treated in the following subsection.
\end{pg}

\subsection{$\mls G$--twisted line bundles}\label{Gtwisted}

\begin{pg}\label{case1} Fix an integer $n$.  Assume that $\mls X = X$ and that $\mls G$ is a $\mu _n$--gerbe over $X$.  Assume furthermore that if $f:X\rightarrow B$ denotes the structural morphism then the map $\mls O_B\rightarrow f_*\mls O_X$ is an isomorphism, and the same remains true after arbitrary base change $B'\rightarrow B$.  This implies that for any scheme $T$ and object $s\in \mls G(T)$ we have an isomorphism $\mu _n\simeq \underline {\text{Aut}}(s)$ and these isomorphisms are functorial in the pair $(T, s)$.  If $\mls L$ is a line bundle on $\mls G$ then for any such pair $(T, s)$ we obtain a line bundle $s^*\mls L$ on $T$ which also comes equipped with an action of $\mu  _n=\underline {\text{Aut}}(s)$.  The line bundle $\mls L$ is called \emph{a $\mls G$--twisted invertible sheaf on $X$} if this action of $\mu _n$ coincides with the standard action induced by the embedding $\mu _n\hookrightarrow \mathbb{G}_m$.

Note that if $\chi $ denotes the standard character $\mu _n\hookrightarrow \mathbb{G}_m$ then for any line bundle $\mls L$ on $\mls G$ the action of $\mu _n = \underline {\text{Aut}}(s)$ on $s^*\mls L$ is given by $\chi ^i$ for some integer $i$ (locally constant on $T$).  It follows that to check that a line bundle $\mls L$ is $\mls G$--twisted it suffices to verify that the two actions coincide for pairs $(T, s)$ with $T$ the spectrum of an algebraically closed field.

If $\mls L$ is a $\mls G$--twisted sheaf on $X$, then $\mls L^{\otimes n}$ descends canonically to an invertible sheaf on $X$ since the stabilizer actions of $\mu _n$ are trivial. We usually write just $\mls L^{\otimes n}$ for the sheaf on $X$ obtained from $\mls L^{\otimes n}$.
\end{pg}

\begin{prop}\label{B.18} There is a natural equivalence between $\Sec (\mls G/X)$ and the stack which to any $B$--scheme $T$ associates the groupoid of pairs $(\mls L, \iota )$, where $\mls L$ is a $\mls G$--twisted sheaf on $X$ and $\iota :\mls L^{\otimes n}\simeq \mls O_X$ is an isomorphism of invertible sheaves on $X$.
\end{prop}
\begin{proof}
If $s:X\rightarrow \mls G$ is a section, then there is a canonical action of $\mu _n=\underline {\text{Aut}}(s)$ on $\mc F:=s_*\mls O_X$ and hence $\mc F$ decomposes canonically as $\oplus _{\chi \in \Z/(n)}\mc F_\chi $.

\begin{lem} Each $\mc F_\chi $ is locally free of rank $1$ on $\mls G$ and for any two $\chi , \epsilon \in \Z/(n)$ the natural map $\mc F_\chi \otimes \mc F_\epsilon \rightarrow \mc F_{\chi +\epsilon }$ is an isomorphism. If $\chi _1:\mu _n\rightarrow \mathbb{G}_m$ denotes the standard inclusion, then $\mc F_{\chi _1}$ is a $\mc G$--twisted sheaf on $X$.
\end{lem}
\begin{proof}
The choice of the section $s$ identifies $\mls G$ with $B\underline {\text{Aut}}(s)\simeq \cB\mu _n\times X$.  The fiber product $X\times _{s, \cB\mu _n, s}X$ is equal to $X\times \cB\mu _n$ from which we see that $s^*\mc F$ is equal to $\mls O_X\otimes _{\Z}\Z[X]/(X^n-1)$ with the natural action of $\mu _n$.  From this the lemma follows.
\end{proof}

In particular, by sending a section $s:X\rightarrow \mc G$ to $\mc F_{\chi _1}$ with the isomorphism $\mc F_{\chi _1}^{\otimes n}\simeq \mls O_{\mls G}$ we obtain a functor
\begin{equation}\label{sectwist}
\Sec (\mls G/X)\rightarrow (\text{stack of pairs $(\mls L, \iota )$ as in \ref{B.18}}).
\end{equation}
Conversely, given a $\mls G$--twisted sheaf $\mls L$ with an isomorphism $\iota :\mls L^{\otimes n}\simeq \mls O_X$ we can form the cyclic algebra $\mls A = \oplus _{i=0}^{n-1} \mls L^{\otimes i}$ with multiplication induced by the natural maps $\mls L^{\otimes i}\otimes \mls L^{\otimes j}\rightarrow \mls L^{\otimes i+j}$ and the map $\iota $.  We can then consider 
$$
\Sp _{\mls G}(\mls A)\rightarrow \mls G.
$$
If $\mls G = \cB\mu _n\times X$ and $X\rightarrow \mls G$ is the map defined by the trivial torsor, then the restriction of $(\mls L, \iota )$ to $X$ is  an invertible sheaf $L$ on $X$ with an isomorphism $L^{\otimes n}\simeq \mls O_X$ and the   
restriction of $\mls A$ to $X$ is just the cyclic algebra $\oplus _{i=0}^{n-1}L^{\otimes i}$ with $\mu _n$ acting $L^{\otimes i}$ through the character $u\mapsto u^i$.  It follows that the projection $\Sp _{\mls G}(\mls A)\rightarrow X$ is an isomorphism and hence defines a point of $\Sec (\mls G/X)$.  We leave to the reader that this defines an inverse to \ref{sectwist}. 
\end{proof}

Let $\scr {P}ic_{X/B}^{\mls G}$ denote the stack over $B$ which to any $B$--scheme $T$ associates the groupoid of $\mls G_T$--twisted invertible sheaves on $X_{T}$ (where $\mls G_T$ denotes $\mls G\times _BT$ etc.).

\begin{prop} {\rm (i)} The stack $\scr Pic^{\mls G}_{X/B}$ is an algebraic stack locally of finite presentation over $B$.

{\rm (ii)} The sheaf (with respect to the fppf topology) associated to the presheaf
$$
T\mapsto \{\text{isomorphism classes in $\scr Pic^{\mls G}_{X/B}$}\}
$$
is representable by a separated  algebraic space locally of finite presentation over $B$.
\end{prop}
\begin{proof} This follows for example by the same argument used in \cite[Appendix]{Artin}.
\end{proof}

%\begin{prop}\label{finitetypeprop} Let $\pi :P\rightarrow X$ be a proper surjective morphism with $P$ an algebraic space %flat over $B$.  If $\mls G'$ denotes  $\mls G\times _XP$, then the morphisms
%$$
%\pi ^*:\mls Pic_{X/B}^{\mls G}\rightarrow \mls Pic_{P/B}^{\mls G}, \ \ \ \pi ^*:\underline {\text{\rm Pic}}^{\mls %G}_{X/B}\rightarrow \underline {\text{\rm Pic}}^{\mls G}_{P/B}
%$$
%are of finite type.
%\end{prop}
%\begin{proof} It suffices to prove that the morphism $\pi ^*:\mls Pic_{X/B}^{\mls G}\rightarrow \mls Pic_{P/B}^{\mls G}$ %is of finite type. This we prove using a method similar to the one
%\end{proof}

\begin{pg}  There is a natural morphism
\begin{equation}\label{Picmap}
\underline {\text{Pic}}_{X/B}^{\mls G}\rightarrow \underline {\text{Pic}}_{X/B}, \ \ \mls L\mapsto \mls L^{\otimes n}.
\end{equation}
Denote by $\underline {\text{Pic}}_{X/B}^{\mls G}[n]$ the inverse image of the identity $e = [\mls O_X]\in \underline {\text{Pic}}_{X/B}$.  Another description of this space is as follows. The map \ref{Picmap} lifts naturally to a morphism
$$
\pi ^*:\scr Pic_{X/B}^{\mls G}\rightarrow \scr Pic_{X/B}, \ \ \mls L\mapsto \mls L^{\otimes n}.
$$
Let $e:B\rightarrow \scr Pic_{X/B}$ be the morphism corresponding to $\mls O_X$, and set
$$
\scr Pic_{X/B}^{\mls G}[n]:= \scr Pic_{X/B}^{\mls G}\times _{\scr Pic_{X/B}, e}B.
$$
The stack $\scr Pic_{X/B}^{\mls G}[n]$ classifies pairs $(\mls L, \iota )$, where $\mls L$ is a $\mls G$--twisted sheaf on $X$ and $\iota :\mls L^{\otimes n}\simeq \mls O_X$ is an isomorphism.  The space $\underline {\text{Pic}}_{X/B}^{\mls G}[n]$ is the coarse moduli space of $\scr Pic_{X/B}^{\mls G}[n]$ via the map sending $(\mls L, \iota )$ to the class of $\mls L$.  In fact, $\scr Pic_{X/B}^{\mls G}[n]$ is a $\mu _n$--gerbe over $\underline {\text{Pic}}_{X/S}^{\mls G}[n]$. 
Indeed  any point $P$ of $\underline {\text{Pic}}_{X/B}^{\mls G}[n](B)$ can fppf--locally on $B$ be represented by a $\mls G$--twisted sheaf $\mls L$ with $\mls L^{\otimes n}\simeq \mls O_X$. Conversely for any two pairs $(\mls L, \iota )$ and $(\mls L', \iota ')$ defining the same point of $\underline {\text{Pic}}_{X/B}^{\mls G}[n]$, the two $\mls G$--twisted sheaves $\mls L$ and $\mls L'$ become isomorphic after making an fppf base change on $B$.  Then $\iota $ and $\iota '$ differ by a section of $f_*\mls O_{\mls G}^* \simeq \mathbb{G}_m$. After making another fppf base change on $B$ so that this unit becomes an $n$--th power, the two pairs $(\mls L, \iota )$ and $(\mls L', \iota ')$ become isomorphic.  Moreover, this isomorphism is unique up to multiplication by an element of $\mu _n(B)$.
\end{pg}

\begin{prop} The algebraic space $\underline {\text{\rm Pic}}_{X/B}^{\mls G}[n]$ is of finite type over $B$.
\end{prop}
\begin{proof} To prove this statement it suffices by a standard limit argument to consider the case when $B$ is a noetherian scheme.  In this case by noetherian induction it suffices to find a dominant morphism $B'\rightarrow B$ such that the restriction of $\underline {\text{Pic}}_{X/B}^{\mls G}$ to $B'$ is of finite type.  We may therefore also assume that $B$ is a noetherian affine integral scheme.  By Chow's lemma  there exists a proper surjection $P\rightarrow \mls G$ with $P$ an algebraic space.  After shrinking on $B$ we may assume that $P$ is also flat over $B$.  Then $\mls G\times _{X}P$ is trivial and by \ref{modprop} the pullback map
$$
\Sec (\mls G/X)\simeq \scr Pic_{X/B}^{\mls G}[n]\rightarrow \scr Pic_{\mls G\times _XP/P}^{\mls G}[n]\simeq \Sec (\mls G\times _XP/P) 
$$
is of finite type. Consequently the map
$$
\underline {\text{\rm Pic}}^{\mls G}_{X/B}[n]\rightarrow \underline {\text{\rm Pic}}^{\mls G}_{P/B}[n]
$$
is also of finite type.  It therefore suffices to show $\underline {\text{\rm Pic}}_{P/B}^{\mls G}[n]$ is of finite type.

Now in the case when $\mls G = \cB\mu _n\times X$ there exists a globally defined $\mls G$-twisted sheaf $\mls L$ with $\mls L^{\otimes n}$ trivial.  Choose one such sheaf $\mls L$.  If $p:\mls G\rightarrow X$ denotes the projection, then for any invertible sheaf $\mls M$ on $X$ with $\mls M^{\otimes n}\simeq \mls O_X$, the sheaf $\mls L\otimes p^*\mls M$ is also a $\mls G$--twisted sheaf.  Conversely, if $\mls N$ is a $\mls G$--twisted sheaf with $\mls N^{\otimes n}$ trivial, then $\mls L\otimes \mls N^{-1}$ descends uniquely to an $n$--torsion sheaf on $X$.  Using this one sees that $\underline {\text{Pic}}_{X/B}^{\mls G}[n]$ is a trivial torsor under $\underline {\text{Pic}}_{X/B}[n]$. The proposition therefore follows from \cite[6.27 and 6.28]{Kleiman}.
\end{proof}

\begin{pg}\label{proofend} Since $\Sec (\mls G/X)\simeq \scr Pic_{X/B}^{\mls G}[n]$ is a $\mu_n$--gerbe over $\underline {\text{Pic}}_{X/B}^{\mls G}[n]$ this completes the proof of (\ref{athm} (ii)) in the special case of \ref{case1} and hence also the proof in general. \qed
\end{pg}

\end{appendix}

%%%%%%%%%%%%%%%%%%%%%%%%%%%%%%%%%%%%%%%

%\bibliographystyle{amsplain}
%\bibliography{mrabbrev,VistoliRefs}

\end{document}